\documentclass[10pt]{article}

\usepackage{amsmath,amssymb,mathrsfs,amsthm}
\usepackage{amsfonts}
\usepackage[fleqn]{mathtools}
\usepackage{dsfont}
\usepackage{latexsym}
\usepackage{soul}
\usepackage{xcolor}
\usepackage{graphicx}
\usepackage{subfig}
\usepackage{physics}
\usepackage[margin=1in]{geometry}
\usepackage{booktabs}
\usepackage{setspace}
\usepackage{float}
\usepackage{multicol}
\usepackage{enumitem}
\usepackage{graphicx}
\usepackage[dvipsnames]{xcolor}
\usepackage{ragged2e}
\usepackage[nice]{nicefrac}
\usepackage[english]{babel}
\usepackage{authblk}

\usepackage[sort&compress,numbers]{natbib}

\usepackage[title]{appendix}

\usepackage{setspace}
\newtheorem{theorem}{Theorem}[section]
\newtheorem{lemma}{Lemma}[section]

\newtheorem{corollary}{Corollary}[section]

\newtheorem{remark}{Remark}[section]
\newtheorem{definition}{Definition}[section]

\newtheorem{rule-def}[theorem]{Rule}
\theoremstyle{definition}

\numberwithin{equation}{section}

\newcommand{\wt}[1]{\widehat{#1}}
\emergencystretch=\maxdimen
\usepackage[
    colorlinks=true,
    linkcolor=magenta,    % For section refs
    citecolor=green,   % For citations
    urlcolor=blue,     % URLs
    filecolor=blue,
    linktoc=all
]{hyperref}

\let\originaleqref\eqref
\renewcommand{\eqref}[1]{\textcolor{magenta}{\originaleqref{#1}}}
\providecommand{\keywords}[1]
{
  \textbf{\textit{Keywords:}} #1
}

\providecommand{\msc}[1]
{
  \textbf{\textit{MSC codes (2020):}} #1
}

\title{On a PDE-ODE-ODE model for meniscus tissue regeneration with space-dependent, {possibly degenerate} diffusion tensor}
\author{Nishith Mohan\thanks{mohan@mathematik.uni-kl.de}}
\author{Ahmad Hussein\thanks{ahmad.hussein@math.rptu.de}}
\author{Bernd Simeon,\thanks{simeon@mathematik.uni-kl.de}}
\author{Christina Surulescu\thanks{c.surulescu@rptu.de}}
\affil{RPTU Kaiserslautern-Landau, Department of Mathematics, \authorcr Gottlieb-Daimler-Str. 48, 67663 Kaiserslautern, Germany.}
\date{\today}
\begin{document}

\makeatletter
\renewcommand\@fnsymbol[1]{}
\makeatother
\maketitle
\begin{abstract}
We investigate a tissue regeneration model that describes the coupled dynamics of human mesenchymal stem cells (MSCs),  chondrocytes, and newly formed extracellular matrix (ECM). The system features motility terms only for the mesenchymal stem cells, more precisely myopic diffusion and haptotaxis towards increasing ECM; all of those terms involve a space-dependent, possibly degenerating tensor. The chondrocytes are obtained from MSCs by differentiation and provide positive feedback to MSCs via dedifferentiation. The considered macroscopic model was obtained by upscaling from a formulation within the kinetic theory of active particles.  Under appropriate assumptions we prove global existence of classical solutions in 2D and 3D in the case of a nondegenerating tensor and global existence of weak solutions in 2D when the diffusion tensor is allowed to degenerate. Numerical simulations are performed in 2D in order to illustrate the solution behavior.
\end{abstract}

\msc{35A01 (primary), % Existence problems for PDEs: global existence, local existence, non-existence
     35K57, % Reaction-diffusion equations
     35A09, % Classical solutions to PDEs
     35B45, % A priori estimates
     35D30 % Weak solutions to PDEs
     35Q92, % PDEs in connection with biology and other natural sciences
     92C17 (secondary)  % Cell movement (chemotaxis, haptotaxis, etc.) (secondary classification)
     }
\section{Introduction and model setup}
\label{sec111}
The last two decades have witnessed a substantial development of mathematical modeling for tissue regeneration processes involved in wound healing, bone, and tissue regeneration, see e.g. \cite{Waters2021} for a comprehensive reference and \cite{Surulescu2026} for a classification of continuous approaches, along with comments on main issues related to tissue regeneration in a larger, interdisciplinary context. 

Here we focus on a macroscopic setting which couples a reactions-diffusion taxis partial differential equation (PDE) with two ordinary differential equations (ODEs) acting on the population scale of two cell phenotypes (MSCs and chondrocytes) and of ECM produced by the latter. The model,  however, takes into account the multiscality of the relevant biological and physical processes. Its deduction is very similar to that employed in \cite{Surulescu2026} and \cite{Hussein}: it starts from a formulation within the kinetic theory of active particle (KTAP) framework developed by N. Bellomo et al. (see \cite{Bellomo2017} and references therein) and performs a parabolic upscaling to the mentioned macroscopic PDE-ODE-ODE system. More precisely, we consider the space-time dynamics of MSCs performing growth, spread, and differentiation into chondrocytes, as well as dedifferentiation of the latter and their production of ECM, along with degradation of newly formed matrix by MSCs (e.g., upon secreting matrix degrading enzymes) and natural decay. Since MSCs are the only ones supposed to actively migrate, we will only perform the upscaling for this species; the dynamics of ECM and chondrocytes only involves production and decay and can be directly modeled on the macroscale.

We start at the microscopic level with the ODE system
\begin{align*}
\dot x&=v    \\
\dot v&=0\\
\dot y&=k_+\frac{\tau }{K}(R_0-y)-k_-y=:G(\tau,y)
\end{align*}
characterizing for each MSC its trajectory: $x$ denotes position, $v$ velocity, and $y$ the so-called 'activity variable' (in the KTAP) giving the amount of receptors bound to the ECM - represented here by its macroscopic density $\tau$. The second ODE assumes that there are no forces acting on the MSCs, whereas the third ODE merely describes mass action kinetics for the mentioned receptor binding, with $k_\pm$ the attachment/detachment rates, $R_0$ the total amount of receptors on a cell (assumed to be constant) and $K$ the carrying capacity of ECM. Similarly to \cite{Surulescu2026}, rescaling and denoting $z:=y^*-y$ (representing deviations of the activity variable from its steady-state $y^*$) leads to an ODE for $z$ and on the mesoscale to a kinetic transport equation (KTE) for the density $p(t,x,v,z)$ of MSCs sharing at time $t$ and position $x$ the velocity regime $v$ and deviation $z$ from the equilibrium receptor binding state:
\begin{align}\label{eq:KTE}
\partial_tp+\nabla_x\cdot (vp)+\partial_z(G(\tau,z)p)=\mathcal L[\lambda (z)]p+\varepsilon^2 \left (-\alpha_1\frac{p}{1+c_1}+\alpha _2\frac{c_2}{1+c_2}+\beta p\left (1-\frac{c_1}{K}-\frac{c_2}{K}-\frac{\tau}{K}\right )\right ),
\end{align}
where $c_1$ and $c_2$ denote the macroscopic densities of MSCs and chondrocytes, respectively; $c_1(t,x)=\int _{V\times Z}p(t,x,v,z)d(v,z)$, with $V=s\mathbb S^{n-1}$ representing the velocity space ($s$ is average speed of MSCs, assumed constant) and $Z=(y^*-1,y*)$. 
The operator $\mathcal L[\lambda (z)]$ represents cell reorientation according to the underlying anisotropy of the scaffold on which the whole dynamics takes place and is given by 
$ \mathcal L[\lambda ]p:=\lambda \left (\frac{q(x,\hat v)}{\omega }\int _Vpdv-p\right ),$
with $\hat v=v/|v|$ and $q(x,\hat v)$ representing the orientation distribution of scaffold fibers. The remaining terms on the right-hand side of \eqref{eq:KTE} describe differentiation of MSCs into chondrocytes with rate $\alpha_1$, dedifferentiation of chondrocytes back to MSCs with rate $\alpha _2$, and proliferation of MSCs, which is assumed to have rate $\beta $ and be of the logistic type - limited by the available amount of MSC, chondrocyte, and ECM densities. These processes are known to need much more time to happen than transport and reorientation; this motivates the scaling with the very small quantity $\varepsilon ^2$, which will also rescale time from mesoscopic to macroscopic, while the space variable will be rescaled only by $\varepsilon$. For more details on the operators, spaces, and their properties, we refer to \cite{Surulescu2026}. Like there, a parabolic scaling leads after analogous calculations to the macroscopic reaction-diffusion-transport PDE
\begin{align}\label{eq:macro-c1}
\partial_tc_1=\nabla \nabla :\left (\mathbb Dc_1\right )-\nabla \cdot \left (\chi (\tau)\mathbb Dc_1\nabla \tau\right )-\alpha _1\frac{c_1}{1+c_1}+\alpha _2 \frac{c_2}{1+c_2}+\beta c_1\left (1-\frac{c_1}{K}-\frac{c_2}{K}-\frac{\tau}{K}\right ),
\end{align}
where $\nabla \nabla:(\mathbb{D} c_1)=\nabla \cdot( \mathbb{D} \nabla c_1 + c_1 \nabla \cdot \mathbb{D})$ represents myopic diffusion and $\mathbb D(x):=\frac{1}{\lambda _{10}}\int _Vv\otimes v\frac{q(x,\hat v)}{\omega }dv$ is the MSC diffusion tensor carrying the anisotropy information about the fibre distribution of the scaffold, $\chi (\tau )=\frac{k_-k_+\lambda _{11}}{K(B(\tau)^2(B(\tau )+\lambda _{10})}$ the haptotactic sensitivity of MSCs towards ECM gradients\footnote{$\lambda _{10}$ and $\lambda _{11}$ are constants involved in the turning rate $\lambda (z)$, see \cite{Surulescu2026}.}, and $B(\tau )=k_+\frac{\tau}{K}+k_-$. This equation is to be supplemented with ODEs for $c_2$ and $\tau$ and the whole system endowed with initial conditions. Although the deduction of \eqref{eq:macro-c1} has been performed for $x\in \mathbb{R}^n$, the PDE can be formulated on a bounded domain $\Omega \subset \mathbb R^n$ and no-flux boundary conditions can be provided; refer again to \cite{Surulescu2026}. 

The model obtained here is in several respects simpler than that in \cite{Surulescu2026}: it involves only haptotaxis towards ECM density and does not feature the dynamics of a differentiation medium, nor that of fibre bound proteins acting as a haptotactic cue; the chondrocytes are assumed here not to diffuse, and the transition rates $\alpha_1, \alpha _2$ are here constants, whereas in \cite{Surulescu2026} they depend on the concentration of differentiation medium and on a stress-related quantity which realized there the coupling with fluid and scaffold mechanics. In this note we want to focus on the space-dependent tensor $\mathbb D$, which can be potentially degenerate. We will therefore address the global existence of the following problem:
\begin{equation}
  \begin{cases}
  \label{model}
    \partial_t c_{1} = \nabla \cdot (\mathbb{D} \nabla c_1 + c_1 \nabla \cdot \mathbb{D}) - \chi \nabla \cdot (c_1 \mathbb{D} \nabla \tau)- \alpha_1  \frac{c_{1}}{1 + c_{1}}  + \alpha_2 \frac{c_{2}}{1 + c_{2}} \\[10pt]
    \hspace{6.8cm} + \, \beta c_{1}(1 - c_{1} - c_{2} - \tau),                                                                                   &   x \in \Omega,\ t > 0,  \\[10pt]
    \partial_t c_2 =  \alpha_1 \frac{c_1}{1 + c_1} - \alpha_2 \frac{c_2}{1 + c_2},                                                         &   x \in \Omega,\ t > 0,   \\[5pt]
    \partial_t \tau = -\delta_{1} c_1 \tau - \mu \tau +  \delta_{2}\frac{c_2 }{1 + c_2},                                                                  &   x \in \Omega,\ t > 0,
  \end{cases}
\end{equation}
subject to zero-flux boundary conditions\footnote{$\nu$ denotes the outward unit normal on $\partial \Omega$}
\begin{equation}
  \label{boundary_conditions}
   (\mathbb{D} \nabla c_1) \cdot \nu = \chi (c_1 \mathbb{D} \nabla \tau) \cdot \nu - c_1 (\nabla \cdot \mathbb{D}) \cdot \nu = 0,  \quad x \in \partial \Omega, \; t > 0,
\end{equation}
and to the initial conditions:
\begin{equation}
  \label{initial_conditions}
  c_1(x, 0) = c_{10}, \quad c_2(x, 0) = c_{20}, \quad \tau(x, 0) = \tau_0, \quad x \in \Omega.
\end{equation}
Moreover, we will also perform numerical simulations. The presence of a degenerating diffusion tensor can raise serious challenges both for the analysis and the numerics, therefore, we consider this system instead of that derived in \cite{Surulescu2026}; \eqref{model} actually even replaces $\chi(\tau)$ from \eqref{eq:macro-c1} with a constant. 

The remainder of this paper is organized as follows:  Section \ref{sec:preliminaries} is setting the existence problem in its bibliographical context and enumerates the main notations and conventions used in the analysis. Section \ref{sec:global-exist} is dedicated to the proof (in 2D and 3D) of the global existence of classical solutions to the problem considered. Section \ref{sec:global-weak-deg} provides the proof for the global existence of weak solutions (in 2D) to \eqref{model} under less restrictive assumptions on the diffusion tensor $\mathbb D$, which is now allowed to infer a certain kind of degeneracy. The mathematical analysis is followed in Section \ref{sec:numerics} by numerical discretization and simulations of the solution behavior in various scenarios related to the degeneracy of $\mathbb D$. The section is completed by a discussion of the simulation results in light of the employed methods.
%%%%%%%%%%%%%%
\section{Preliminaries and notations for the analysis}\label{sec:preliminaries}
%%%%%%%%%%%%%
  The system of equations in \eqref{model} represents a haptotaxis model that incorporates non-Fickian, myopic diffusion, as considered, e.g. in \cite{Engwer2014}, see also \cite{Engwer2016,engwer2016effective}. Such and related models have also been studied in \cite{winkler2016global,heihoff2023global,Conte2021,Conte2023,CMS-26b,Corbin2020,Dietrich2022,Surulescu2026,Conte_Hiremath_Surulescu}. The mathematical analysis of \eqref{model}-\eqref{initial_conditions} builds upon the framework proposed in \cite{heihoff2023global}, which examined a haptotaxis model for cancer invasion incorporating a potentially degenerate diffusion tensor and non-Fickian diffusion. However, the complexities in our system are compounded by the presence of a production term $\frac{c_2}{1 + c_2}$ in the third equation of the model. In the cancer invasion model studied in \cite{heihoff2023global}, the component along which haptotaxis is reported only degrades over time. In contrast, in our case, the haptotactic component exhibits growth. Moreover, the component responsible for the growth of the haptotactic component itself satisfies an ordinary differential equation, further complicating the situation. 

  In both two and three-dimensional settings, we assume, as in \cite{heihoff2023global}, that the diffusion tensor $\mathbb{D}$ is uniformly positive definite. This ensures that the diffusion operator remains elliptic, allowing us to employ standard techniques from parabolic PDE theory. The global classical solvability of \eqref{model}-\eqref{initial_conditions} is derived through an appropriate change of variable, a standard approach for such haptotaxis systems.

  In two spatial dimensions, the analysis relies on constructing a quasi-energy inequality. We then establish a $ L^\infty$ estimate for the first solution component, $c_1$, using a modified Gagliardo-Nirenberg interpolation inequality \cite[Lemma A.5]{tao2014energy} together with a Moser-Alikakos\cite{alikakos1979lp} iteration argument. This estimate, combined with a suitable testing procedure, yields $W^{1,4}$ bounds for $c_2$ and $\tau$. These bounds collectively ensure that global extensibility in two dimensions, thereby completing our proof.

  In three spatial dimensions, the proof of global classical solvability relies on the coupled $L^p$-estimate approach as in \cite{pang2018global,dai2022new}. The analysis proceeds through a bootstrap of the regularity of the solution component $c_1$ from $L^1$ to $L^\infty$ via intermediate $L^2$ and $L^p \; (p > 1)$ estimates. A key ingredient is the derivation of a suitable weighted inequality, which, together with a comparison argument for an associated ODE, yields the required $L^2$ estimate under an appropriate smallness assumption on the initial data. Standard testing arguments subsequently obtain the higher-order $L^p$ estimates, while the final $L^\infty$ bound follows from the Moser-Alkalosis iteration technique.
   
  Next, we relax the assumptions on the diffusion tensor $\mathbb{D}$ by removing the requirement of uniform positive definiteness on $\overline{\Omega}$. Under this weaker set of assumptions, we establish the existence of weak solutions to \eqref{model}-\eqref{initial_conditions} in the two-dimensional setting.
   
 In the sequel, we will use the following notations and conventions:
 \begin{itemize}
	\item we abbreviate the integrals $\int_\Omega f(x) dx$ and $\int_0^T \int_\Omega f(x, t)\, dx \,dt$ $(T > 0)$ as $\int_\Omega f(x)$ and $\int_0^T \int_\Omega f(x, t)$;
	\item the numbering of the constants $ C_i, K_i > 0, i = 1, 2, 3, \ldots$ holds only within the lemma and its proof in which the constants are used. The sequence restarts once the proof is over;
    \item { if $X$ denotes a function space, then $X^*$ represents its dual;}
    \item the assumptions on the diffusion tensor and the initial data, together with the corresponding basic estimates and the main results on classical and weak solvability of \eqref{model}-\eqref{initial_conditions}, are presented separately in each section to keep the analysis self-contained.
    \item we assume the reader to be familiar with the standard Lebesgue and Sobolev spaces and their usual properties,  with the Zygmund space (see e.g., \cite[Chapter 4, Definition 6.1]{bennet}) 
    $$L\log L:=\left\{u\in L^1(\Omega)\ :\ \int _{\{u>1\}}|u|\log |u|dx<\infty \right\},$$
    { and, more generally, the so-called Orlicz space $L^\Phi$ associated to a Young function $\Phi:\mathbb R\to \mathbb R^+$ (i.e., a convex, even, and lower semicontinuous function s.t. $\Phi (0)=0$ and $\lim \limits_{x\to \infty }\Phi (x)=\infty$)
    $$L^\Phi :=\left \{u\ :\ \int _\Omega \Phi (\kappa |u|)<\infty \ \text{for some }\kappa >0\right \}\quad \text{(see e.g. \cite{book-rao})}.$$}
 \end{itemize}

%%%%%%%%%%%%%%%%%%%%%%%%%%%%%%%%%%%%%%%%%%%%%%%%%%%%%%%%%%%%%%%%%%%%%%%%%%%%%%%%%%%%%%%%%%%%%%%%%%%%%%%%%%%%%%%%%%%%%%%%%%%%%%%%%%%%%%%%%%%%%%%%%%%%%%%%%%%%%%%%%%%%%%%%%%%%%%%%%%%%%%%%%%%%%%%%%%%%%%%%%%%%%%%%%%%%%%%%%%%%%%%%%
%%%%%%%%%%%%%%%%%%%%%%%%%%%%%%%%%%%%%%%%%%%%%%%%%%%%%%%%%%%%%%%%%%%%%%%%%%%%%%%%%%%%%%%%%%%%%%%%%%%%%%%%%%%%%%%%%%%%%%%%%%%%%%%%%%%%%%%%%%%%%%%%%%%%%%%%%%%%%%%%%%%%%%%%%%%%%%%%%%%%%%%%%%%%%%%%%%%%%%%%%%%%%%%%%%%%%%%%%%%%%%%%%
%%%%%%%%%%%%%%%%%%%%%%%%%%%%%%%%%%%%%%%%%%%%%%%%%%%%%%%%%%%%%%%%%%%%%%%%%%%%%%%%%%%%%%%%%%%%%%%%%%%%%%%%%%%%%%%%%%%%%%%%%%%%%%%%%%%%%%%%%%%%%%%%%%%%%%%%%%%%%%%%%%%%%%%%%%%%%%%%%%%%%%%%%%%%%%%%%%%%%%%%%%%%%%%%%%%%%%%%%%%%%%%%%

%\section{Preliminaries}

%%%%%%%%%%%%%%%%%%%%%%%%%%%%%%%%%%%%%%%%%%%%%%%%%%%%%%%%%%%%%%%%%%%%%%%%%%%%%%%%%%%%%%%%%%%%%%%%%%%%%%%%%%%%%%%%%%%%%%%%%%%%%%%%%%%%%%%%%%%%%%%%%%%%%%%%%%%%%%%%%%%%%%%%%%%%%%%%%%%%%%%%%%%%%%%%%%%%%%%%%%%%%%%%%%%%%%%%%%%%%%%%%
%%%%%%%%%%%%%%%%%%%%%%%%%%%%%%%%%%%%%%%%%%%%%%%%%%%%%%%%%%%%%%%%%%%%%%%%%%%%%%%%%%%%%%%%%%%%%%%%%%%%%%%%%%%%%%%%%%%%%%%%%%%%%%%%%%%%%%%%%%%%%%%%%%%%%%%%%%%%%%%%%%%%%%%%%%%%%%%%%%%%%%%%%%%%%%%%%%%%%%%%%%%%%%%%%%%%%%%%%%%%%%%%%
%%%%%%%%%%%%%%%%%%%%%%%%%%%%%%%%%%%%%%%%%%%%%%%%%%%%%%%%%%%%%%%%%%%%%%%%%%%%%%%%%%%%%%%%%%%%%%%%%%%%%%%%%%%%%%%%%%%%%%%%%%%%%%%%%%%%%%%%%%%%%%%%%%%%%%%%%%%%%%%%%%%%%%%%%%%%%%%%%%%%%%%%%%%%%%%%%%%%%%%%%%%%%%%%%%%%%%%%%%%%%%%%%

\section{\texorpdfstring
{Global existence of classical solutions to \eqref{model}-\eqref{initial_conditions} under assumptions \ref{lab_ass} on $\mathbb{D}$ in two and three spatial dimensions}
{Classical solutions}}\label{sec:global-exist}

In this section, we derive our results on the global classical solvability of \eqref{model}-\eqref{initial_conditions} in two and three dimensions, assuming sufficiently regular initial data and some rather stringent conditions on $\mathbb{D}$. For the sake of convenience, we restate our problem in a more concise form:
\begin{equation}
  \begin{cases}
    \label{model1}
     \partial_t c_{1} = \nabla \cdot (\mathbb{D}\nabla c_{1} + c_{1} \nabla \cdot \mathbb{D}) - \chi \nabla \cdot (c_{1} \mathbb{D} \nabla \tau) - \alpha_1  \tfrac{c_{1}}{1 + c_{1}}  + \alpha_2 \tfrac{c_{2}}{1 + c_{2}}\\[10pt]
     \hspace{7cm} ~~~~ + \beta c_{1}(1 - c_{1} - c_{2} - \tau),                                                                                                                 &  ~~ x \in \Omega, t > 0,  \\[5pt]
    \partial_t c_{2} =  \alpha_1 \frac{c_{1}}{1 + c_{1}} - \alpha_2 \frac{c_{2}}{1 + c_{2}},                                                                                 & ~~ x \in \Omega, t > 0,            \\[5pt]
    \partial_t \tau_ = -\delta c_{1} \tau - \mu \tau + \frac{c_{2} }{1 + c_{2}},                                                                                           & ~~ x \in \Omega, t > 0,            \\[5pt]
    (\mathbb{D} \nabla c_{1}) \cdot \nu = \chi (c_{1} \mathbb{D}\nabla \tau)\cdot \nu - c_{1}(\nabla \cdot \mathbb{D}) \cdot \nu = 0,     & ~~ x \in \partial \Omega, t > 0,   \\[5pt]
    c_{1}(\cdot, 0) = c_{10}, ~~ c_{2}(\cdot, 0) = c_{20}, ~~ \tau_{}(\cdot, 0) = \tau_{0},                                                                                  & ~~ x \in \Omega.
  \end{cases}
\end{equation}

%%%%%%%%%%%%%%%%%%%%%%%%%%%%%%%%%%%%%%%%%%%%%%%%%%%%%%%%%%%%%%%%%%%%%%%%%%%%%%%%%%%%%%%%%%%%%%%%%%%%%%%%%%%%%%%%%%%%%%%%%%%%%%%%%%%%%%%%%%%%%%%%%%%%%%%%%%%%%%%%%%%%%%%%%%%%%%%%%%%%%%%%%%%%%%%%%%%%%%%%%%%%%%%%%%%%%%%%%%%%%%%%%
%%%%%%%%%%%%%%%%%%%%%%%%%%%%%%%%%%%%%%%%%%%%%%%%%%%%%%%%%%%%%%%%%%%%%%%%%%%%%%%%%%%%%%%%%%%%%%%%%%%%%%%%%%%%%%%%%%%%%%%%%%%%%%%%%%%%%%%%%%%%%%%%%%%%%%%%%%%%%%%%%%%%%%%%%%%%%%%%%%%%%%%%%%%%%%%%%%%%%%%%%%%%%%%%%%%%%%%%%%%%%%%%%

\subsection{Assumptions on \texorpdfstring{$\mathbb{D}$}{D} : classical solutions (non degenerate case)}
\label{lab_ass}
In order to show the global existence of classical solutions for \eqref{model1} the diffusion tensor $\mathbb{D} \in C^2(\overline{\Omega}; \mathbb{R}^{n \times n}) ~  (n \subseteq \mathbb{N})$ is assumed to be positive definite everywhere, with components $\mathbb{D} = (d_{ij})_{i, j = 1, \ldots, n}, ~~ d_{ij} \in C^{ {2}}(\overline{\Omega})$. Additionally, $\mathbb{D}$ satisfies  $(\nabla \cdot \mathbb{D}) \cdot \nu = 0$ on the boundary $\partial \Omega$. Moreover, $\mathbb{D}$ satisfies a uniform boundedness condition, i.e., there exist two constants $D_1, D_2 > 0$ such that for all $\xi \in \mathbb{R}^n$ and $x \in \overline{\Omega}$ the following holds true
\begin{equation}
  \label{assumption0}
  D_1 |\xi|^2 \leq \sum_{i, j = 1}^n d_{ij}(x) \xi_j \xi_i \leq D_2 |\xi|^2.
\end{equation}
Furthermore, under the assumptions mentioned above we can fix a $M > 0$ such that
\begin{equation}
  \label{assumption1}
  \|\nabla \cdot \mathbb{D}\|_{L^\infty(\Omega)} \leq M \quad \text{and} \quad \|\nabla \cdot (\nabla \cdot \mathbb{D})\|_{L^\infty(\Omega)} \leq M.
\end{equation}
A regularity property required to carry out the analysis is related with the divergence of $\mathbb{D}$ (applied column-wise) and is given as follows:
\begin{definition}[{\cite[Definition 2.3]{heihoff2023global}}]
\label{req_weakSol}
  Given that $\mathbb{D}$ is a positive semidefinite matrix, $\mathbb{D} = (\mathbb{D}_1, \cdots, \mathbb{D}_n) \in L^1(\Omega; \mathbb{R}^{n \times n})$, and $\nabla \cdot \mathbb{D} := (\nabla \cdot \mathbb{D}_1, \cdots, \nabla \cdot \mathbb{D}_n) \in L^1(\Omega; \mathbb{R}^n)$, it satisfies the following divergence estimate with exponent $\lambda\footnote{$\lambda$ here is unrelated to $\lambda(z)$ in the cell reorientation operator $\mathcal{L}[\lambda(z)]$ from Section~\ref{sec111}.} \in [\tfrac{1}{2}, \tfrac{2}{3}]$ and a constant $A > 0$:
    \begin{equation}
      \label{assumption2}
      \left|\int_\Omega (\nabla \cdot \mathbb{D}) \cdot \Phi\right| \leq A \left( \int_\Omega (\Phi \cdot \mathbb{D} \Phi)^\lambda + 1\right)
    \end{equation}
    for all $\Phi \in C^0(\overline{\Omega};\mathbb{R}^n)$.
\end{definition}
Also, we select the initial data such that
\begin{equation}
\begin{cases}
  \label{initial_condition_assumption}
  c_{10}, c_{20}, \tau_0 \in C^{2 + \vartheta}(\overline{\Omega}), ~\text{for some}~\vartheta \in (0, 1), ~\text{with}      \\[5pt]
  (\mathbb{D} \nabla c_{10}) \cdot \nu = (\mathbb{D} \nabla \tau_0) \cdot \nu = \nabla c_{20}\cdot \nu  = 0 ~\text{on}~ \partial \Omega, ~~\text{and}         \\[5pt]
  c_{10} > 0, c_{20} > 0, \tau_0 > 0 ~\text{on}~\overline{\Omega}.\quad %\text {\cred Can this be weakened to $\tau_0\ge 0$?}
\end{cases}
\end{equation}
%%%%%%%%%%%%%%%%%%%%%%%%%%%%%%%%%%%%%%%%%%%%%%%%%%%%%%%%%%%%%%%%%%%%%%%%%%%%%%%%%%%%%%%%%%%%%%%%%%%%%%%%%%%%%%%%%%%%%%%%%%%%%%%%%%%%%%%%%%%%%%%%%%%%%%%%%%%%%%%%%%%%%%%%%%%%%%%%%%%%%%%%%%%%%%%%%%%%%%%%%%%%%%%%%%%%%%%%%%%%%%%%%
%%%%%%%%%%%%%%%%%%%%%%%%%%%%%%%%%%%%%%%%%%%%%%%%%%%%%%%%%%%%%%%%%%%%%%%%%%%%%%%%%%%%%%%%%%%%%%%%%%%%%%%%%%%%%%%%%%%%%%%%%%%%%%%%%%%%%%%%%%%%%%%%%%%%%%%%%%%%%%%%%%%%%%%%%%%%%%%%%%%%%%%%%%%%%%%%%%%%%%%%%%%%%%%%%%%%%%%%%%%%%%%%%
\subsubsection{Main results: classical solutions}

The main results of this section are stated below. Their proofs will then be performed throughout the remainder of this Section \ref{sec:global-exist}. We begin with the global classical solvability result in a two-dimensional setting.
\begin{theorem}
\label{classical_main_th1}
   Let $\Omega \subseteq \mathbb{R}^2$ be a bounded domain with a smooth boundary. The parameters $\alpha_1, \alpha_2, \beta, \delta, \mu$, and $\chi$ are all positive, and $\lambda \in [\tfrac{1}{2}, \tfrac{2}{3}]$ with $\tfrac{\lambda}{1 - \lambda} \leq 2$. The diffusion tensor $\mathbb{D} \in C^2(\overline{\Omega}; \, \mathbb{R}^{2 \times 2})$ is assumed to be symmetric, positive definite, and uniformly bounded on $\overline{\Omega}$. Additionally, $\mathbb{D}$ satisfies the divergence estimate of Definition \ref{req_weakSol} with exponent $\lambda$ and  the boundary condition $(\nabla \cdot \mathbb{D}) \cdot \nu = 0$  on  $\partial \Omega$. Moreover, let $c_{10}, c_{20}, \tau_0 \in C^{2 + \vartheta}(\overline{\Omega})$, $\vartheta \in (0, 1)$, be some initial conditions with $c_{10}, c_{20}, \tau_0 > 0$ on $\overline{\Omega}$ such that $(\mathbb{D} \nabla c_{10}) \cdot \nu = (\mathbb{D} \nabla \tau_0) \cdot \nu =0$ on $\partial \Omega$ holds true. Then there exists a unique {tuple} of nonnegative functions $c_1, c_2, \tau \in C^{2, 1}(\overline{\Omega} \times [0, \infty))$ such that the triple $(c_1, c_2, \tau)$ is a global classical solution to \eqref{model1}.
\end{theorem}
%%%%%%%%%%%%%%%%%%%%%%%%%%%%%%%%%%%%%%%%%%%%%%%%%%%%%%%%%%%%%%%%%%%%%%%%%%%%%%%%%%%%%%%%%%%%%%%%%%%%%%%%%%%%%%%%%%%%%%%%%%%%%%%%%%%%%%%%%%%%%%%%%%%%%%%%%%%%%%%%%%%%%%%%%%%%%%%%%%%%%%%%%%%%%%%%%%%%%%%%%%%%%%%%%%%%%%%%%%%%%%%%%%%%%%%%%%%%%%%%%%%%%%%%%%%%%%%%%%%%%%%%%%%%%%%%%%%%%%%%%%%%%%%%%%%%%%%%%%%%%%%%%%%%%%%%%%%%%%%%%%%%%%%%%%%%%%%%%%%%%%%%%%%%%%%%%%%%%%%%%%%%%%
The next theorem establishes the global classical solvability of \eqref{model1} in three dimensions.
\begin{theorem}
    \label{3D_Class_Theorem}
  Let $\Omega \subseteq \mathbb{R}^3$ be a bounded domain with smooth boundary, and let $\alpha_1,\alpha_2,\beta,\delta,\mu,\chi>0$. Assume that the diffusion tensor $\mathbb{D}\in C^2(\overline{\Omega} \, ;\mathbb{R}^{3\times3})$ is symmetric, uniformly positive definite and satisfies $(\nabla\cdot\mathbb{D})\cdot\nu=0$ on $\partial\Omega$. Let $c_{10},c_{20},\tau_0\in C^{2+\vartheta}(\overline{\Omega})$, $\vartheta\in(0,1)$, be positive initial data satisfying
  $(\mathbb{D}\nabla c_{10})\cdot\nu = (\mathbb{D}\nabla\tau_0)\cdot\nu = 0 \; \; \text{on }\partial\Omega$.
  Define $\Pi:=\|c_{10}\|_{L^1(\Omega)},$
  and suppose that
  \begin{align*}
  \|e^{-\frac{\chi}{2}\tau_0}c_{10}\|_{L^2(\Omega)}^2 \leq s_0 \qquad \text{and} \qquad  s_0:=\sqrt{\frac{1}{3\Lambda^\ast}},
  \end{align*}
  where $\Lambda^\ast$ is given by \eqref{Lambdastar}. Then there exists a constant $\Pi^\ast>0$ satisfying
  \begin{align*}
  \Lambda_\ast C_{\mathrm{GN}}^{(3)}(\Pi^\ast)^3 +\overline{\Lambda}\Pi^\ast = \frac{2s_0}{3},
  \end{align*}
  where $\Lambda_\ast$ and $\overline{\Lambda}$ are defined in \eqref{Lambda_ast} and \eqref{ov_lambda}, respectively, and depend only on the model parameters and the three-dimensional Gagliardo-Nirenberg interpolation constant $C_{\mathrm{GN}}^{(3)}(\Omega) > 0$, such that, whenever
   \begin{align*}
   \frac{|\Omega|}{2} \left(1+\sqrt{1+\frac{4\alpha_2}{\beta}}\right) < \Pi < \Pi^\ast,
   \end{align*}
 the system \eqref{model1} admits a global classical solution.
\end{theorem}

%\noindent
\begin{remark}
{Notice that, depending mainly on $\Pi$, but also influenced by the quotient $\frac{\alpha_2}{\beta}$, the size of the spatial domain on which this solution exists might be substantially reduced. Thus, a sufficiently large mean value (over space domain) of initial MSC density is required, which essentially depends on the quotient between the rates of chondrocyte dedifferentiation and MSC proliferation.}
\end{remark}

%%%%%%%%%%%%%%%%%%%%%%%%%%%%%%%%%%%%%%%%%%%%%%%%%%%%%%%%%%%%%%%%%%%%%%%%%%%%%%%%%%%%%%%%%%%%%%%%%%%%%%%%%%%%%%%%%%%%%%%%%%%%%%%%%%%%%%%%%%%%%%%%%%%%%%%%%%%%%%%%%%%%%%%%%%%%%%%%%%%%%%%%%%%%%%%%%%%%%%%%%%%%%%%%%%%%%%%%%%%%%%%%%
%%%%%%%%%%%%%%%%%%%%%%%%%%%%%%%%%%%%%%%%%%%%%%%%%%%%%%%%%%%%%%%%%%%%%%%%%%%%%%%%%%%%%%%%%%%%%%%%%%%%%%%%%%%%%%%%%%%%%%%%%%%%%%%%%%%%%%%%%%%%%%%%%%%%%%%%%%%%%%%%%%%%%%%%%%%%%%%%%%%%%%%%%%%%%%%%%%%%%%%%%%%%%%%%%%%%%%%%%%%%%%%%%

\subsection{Local existence and some basic estimates}
On the basis of \eqref{assumption0}, \eqref{assumption1}, and \eqref{assumption2}, $\mathbb{D}$ can be regarded as a positive definite matrix-valued function. Consequently, the diffusion operator in the first equation of \eqref{model1} is strictly elliptic. As a result, much of the existence and regularity theory applicable to the Laplacian operator $\Delta$ extends to this case. We remark that the results established in this subsection are valid in both two and three spatial dimensions. 

Accordingly, the following result on local existence and uniqueness of classical solutions can be adapted from similar systems studied in \cite{tao2014energy} and \cite{pang2017global}:

\begin{lemma}
  \label{local_existence}
Let $\Omega \subseteq \mathbb{R}^n$, with $n \in \{2,3\}$, be a bounded domain with smooth boundary. Let $\alpha_1$, $\alpha_2$, $\beta$, $\delta$, $\mu$, and $\chi$ be positive constants, and assume that \eqref{assumption0}, \eqref{assumption1}, \eqref{assumption2}, and \eqref{initial_condition_assumption} are satisfied. Then there exists a $T_{\max} \in (0, \infty)$ and a uniquely determined classical solution $(c_1, c_2, \tau) \in (C^{2, 1}(\overline{\Omega} \times [0, \infty)))^3$ of \eqref{model1} such that either $T_{\max} = \infty$ or
  \begin{align}
    \label{extensibility_criterion}
    \limsup_{t \nearrow T_{\max}} \left\{\|c_1(\cdot, t)\|_{L^\infty(\Omega)} + \|\tau(\cdot, t)\|_{W^{1, 4}(\Omega)}\right\} = \infty.
  \end{align}
  Moreover, we have $c_1 > 0, c_2 > 0$ and $\tau > 0$ in $\overline{\Omega} \times (0, T_{\max})$.
\end{lemma}

%%%%%%%%%%%%%%%%%%%%%%%%%%%%%%%%%%%%%%%%%%%%%%%%%%%%%%%%%%%%%%%%%%%%%%%%%%%%%%%%%%%%%%%%%%%%%%%%%%%%%%%%%%%%%%%%%%%%%%%%%%%%%%%%%%%%%%%%%%%%%%%%%%%%%%%%%%%%%%%%%%%%%%%%%%%%%%%%%%%%%%%%%%%%%%%%%%%%%%%%%%%%%%%%%%%%%%%%%%%%%%%%%
%%%%%%%%%%%%%%%%%%%%%%%%%%%%%%%%%%%%%%%%%%%%%%%%%%%%%%%%%%%%%%%%%%%%%%%%%%%%%%%%%%%%%%%%%%%%%%%%%%%%%%%%%%%%%%%%%%%%%%%%%%%%%%%%%%%%%%%%%%%%%%%%%%%%%%%%%%%%%%%%%%%%%%%%%%%%%%%%%%%%%%%%%%%%%%%%%%%%%%%%%%%%%%%%%%%%%%%%%%%%%%%%%

\begin{lemma}
\label{lemmal1c1}
 The first solution component of \eqref{model1} satisfies
  \begin{equation}
    \label{l1c1}
    \int_\Omega c_{1}(\cdot, t) \leq \max\left\{\int_\Omega c_{10}, \frac{|\Omega|}{2} \left(1 + \sqrt{1 + \frac{4 \alpha_2}{\beta}}\right)\right\} =: M_1
  \end{equation}
  for all $t \in (0,T_{\max})$. Additionally, 
  \begin{equation}
  \label{c1_bound1}
    \int_t^{t + \gamma} \int_\Omega c_{1}^2 \leq  \frac{M_1 + \alpha_2 |\Omega|}{\beta} + M_1, \qquad \text{for any} ~~  0 < \gamma \leq \min \left\{1, \frac{T_{\max}}{2}\right\}. 
  \end{equation}
\noindent
%{\cred How do we know $T_{\max}>1$? I guess we should rather integrate on $(t,t+\gamma)$, with $\gamma $ sufficiently small. } {\color{blue} - Please check if this is fine}

\end{lemma}

\begin{proof}
  Spatially integrating the first equation of \eqref{model1} over $\Omega$ results in
  \begin{align*}
  \label{l1c11}
    \dfrac{d}{dt}\int_\Omega c_{1} = & - \alpha_1\int_\Omega \dfrac{c_{1}}{1 + c_{1}} + \alpha_2 \int_\Omega  \dfrac{c_{2}}{1 + c_{2}} + \beta \int_\Omega c_{1} - \beta \int_\Omega c_{1}^2 - \beta \int_\Omega c_{1} c_{2}  - \beta \int_\Omega c_{1} \tau
  \end{align*}
 for all $t \in (0, T_{\max})$. Using the positivity of solution components, and the fact that $\frac{c_{i}}{1 + c_{i}} < 1$ for $c_{i} \geq 0,\ i \in \{1, 2\}$, we have
  \begin{equation}
    \label{l1c12}
    \dfrac{d}{dt} \int_\Omega c_{1} \leq \alpha_2 |\Omega| + \beta \int_\Omega c_{1} - \beta \int_\Omega c_{1}^2 \quad \text{for all}~ t \in (0, T_{\max}).
  \end{equation}
  The Cauchy-Schwarz inequality ensures that $\int_\Omega c_{1}^2 \geq \frac{1}{|\Omega|} (\int_\Omega c_{1})^2$, hence
  \begin{equation}
    \label{l1c13}
    \dfrac{d}{dt} \int_\Omega c_{1} \leq \alpha_2 |\Omega| + \beta \int_\Omega c_{1} - \frac{\beta}{|\Omega|} \left(\int_\Omega c_{1}\right)^2  \quad \text{for all}~ t \in (0, T_{\max}).
  \end{equation}
 Applying an ODE comparison principle to \eqref{l1c13} results in
  \begin{equation}
  \label{l1c14}
    \|c_{1}(\cdot, t)\|_{L^1(\Omega)} \leq \max\left\{\|c_{10}\|_{L^1(\Omega)}, \frac{|\Omega|}{2} \left(1 + \sqrt{1 + \frac{4 \alpha_2}{\beta}}\right)\right\} \quad \text{for all}~ t \in (0, T_{\max, })
  \end{equation}
  thus yielding \eqref{l1c1}.
  Since $\beta > 0$ an integration of \eqref{l1c12} in view of \eqref{l1c14} over $(t, t + \gamma)$ for any $0 < \gamma \leq \min \{1, \frac{T_{\max}}{2}\}$ results in
  \begin{align*}
      \int_t^{t + \gamma} \int_\Omega c_1^2 \leq \frac{M_1 + \alpha_2 |\Omega|}{\beta} + M_1,
  \end{align*}
  thus establishing \eqref{c1_bound1}.
\end{proof}
%%%%%%%%%%%%%%%%%%%%%%%%%%%%%%%%%%%%%%%%%%%%%%%%%%%%%%%%%%%%%%%%%%%%%%%%%%%%%%%%%%%%%%%%%%%%%%%%%%%%%%%%%%%%%%%%%%%%%%%%%%%%%%%%%%%%%%%%%%%%%%%%%%%%%%%%%%%%%%%%%%%%%%%%%%%%%%%%%%%%%%%%%%%%%%%%%%%%%%%%%%%%%%%%%%%%%%%%%%%%%%%%%
%%%%%%%%%%%%%%%%%%%%%%%%%%%%%%%%%%%%%%%%%%%%%%%%%%%%%%%%%%%%%%%%%%%%%%%%%%%%%%%%%%%%%%%%%%%%%%%%%%%%%%%%%%%%%%%%%%%%%%%%%%%%%%%%%%%%%%%%%%%%%%%%%%%%%%%%%%%%%%%%%%%%%%%%%%%%%%%%%%%%%%%%%%%%%%%%%%%%%%%%%%%%%%%%%%%%%%%%%%%%%%%%%
\begin{lemma}
  \label{lemma_c2_bound}
  Let $T>0$. Then there exists a constant $M_{c_2}(T) > 0$ such that the second solution component satisfies
  \begin{equation}
    \label{c2_bound}
    \|c_2(\cdot,t)\|_{L^\infty(\Omega)}
    \leq M_{c_2}(T)
    \qquad \text{for all } t\in(0,\widehat{T}),
  \end{equation}
  where
  $\widehat{T}:=\min\{T,T_{\max}\}$.
\end{lemma}

\begin{proof}
  By the nonnegativity of $c_1$ and $c_2$ in $\overline{\Omega}\times(0,T_{\max})$, the second equation of
  \eqref{model1} yields
  $$\partial_tc_{2} \leq \alpha_1\frac{c_1}{1+c_1}  \leq \alpha_1 \qquad \text{in } \overline{\Omega}\times(0,T_{\max}).$$
%  where we have used $0\leq \frac{c_1}{1+c_1}\leq 1$ for all $c_1\geq 0$.
  Integrating the above differential inequality over $(0,t)$ gives
  $$c_2(x,t) \leq c_{20}(x)+\alpha_1t \leq \|c_{20}\|_{L^\infty(\Omega)}+\alpha_1\widehat{T} =: M_{c_2}(T)$$
  for all $x\in\overline{\Omega}$ and $t\in(0,\widehat{T})$.
\end{proof}
%%%%%%%%%%%%%%%%%%%%%%%%%%%%%%%%%%%%%%%%%%%%%%%%%%%%%%%%%%%%%%%%%%%%%%%%%%%%%%%%%%%%%%%%%%%%%%%%%%%%%%%%%%%%%%%%%%%%%%%%%%%%%%%%%%%%%%%%%%%%%%%%%%%%%%%%%%%%%%%%%%%%%%%%%%%%%%%%%%%%%%%%%%%%%%%%%%%%%%%%%%%%%%%%%%%%%%%%%%%%%%%%%
%%%%%%%%%%%%%%%%%%%%%%%%%%%%%%%%%%%%%%%%%%%%%%%%%%%%%%%%%%%%%%%%%%%%%%%%%%%%%%%%%%%%%%%%%%%%%%%%%%%%%%%%%%%%%%%%%%%%%%%%%%%%%%%%%%%%%%%%%%%%%%%%%%%%%%%%%%%%%%%%%%%%%%%%%%%%%%%%%%%%%%%%%%%%%%%%%%%%%%%%%%%%%%%%%%%%%%%%%%%%%%%%%

\begin{lemma}
  \label{lemma_tau_bound}
  For the third solution component, there exists an $M_\tau > 0$ such that
  \begin{align}
    \label{tau_bound} & \|\tau(\cdot, t)\|_{L^\infty(\Omega)} \leq \frac{r_{1\ast}}{\mu} + \|\tau_{0}\|_{L^\infty(\Omega)} =: M_\tau,\quad \text{for all}~ t \in (0, T_{\max}),
  \end{align}
   where $r_1 \equiv \frac{c_{2}}{1 + c_{2}}$ and $r_{1\ast} := \|r_1\|_{L^\infty(\Omega \times (0, \infty))}$.
\end{lemma}
\begin{proof}
  Let\begin{align*}
    \bar{\tau}(t) := \|\tau_{0 }\|_{L^\infty(\Omega)}e^{- \mu t} + \int_0^t e^{-\mu (t - s)} \|r_1(\cdot, s)\|_{L^\infty(\Omega)} ds
  \end{align*}
  for $t \in (0, \infty)$. Then $\bar{\tau} \leq M_\tau$ on $\overline{\Omega} \times (0, T_{\max})$ and the following holds true:
  \begin{equation*}
    \bar{\tau}_t + \delta c_{1} \bar{\tau} + \mu \bar{\tau} - r_1(x, t) \geq \bar{\tau}_t + \mu \bar{\tau} - \|r_1(\cdot, t)\|_{L^\infty(\Omega)} \geq 0
  \end{equation*}
  on $\overline{\Omega} \times (0, T_{\max})$, owing to the non-negativity of $c_{1}$ and $\bar \tau$. Applying an ODE comparison principle results in $\bar{\tau} \geq \tau$ on $(\Omega \times (0, T_{\max})$ and in particular $\|\tau(\cdot, t)\|_{L^\infty(\Omega)} \leq M_\tau$ for all $t \in (0, T_{\max})$. This gives \eqref{tau_bound}.
\end{proof}

%%%%%%%%%%%%%%%%%%%%%%%%%%%%%%%%%%%%%%%%%%%%%%%%%%%%%%%%%%%%%%%%%%%%%%%%%%%%%%%%%%%%%%%%%%%%%%%%%%%%%%%%%%%%%%%%%%%%%%%%%%%%%%%%%%%%%%%%%%%%%%%%%%%%%%%%%%%%%%%%%%%%%%%%%%%%%%%%%%%%%%%%%%%%%%%%%%%%%%%%%%%%%%%%%%%%%%%%%%%%%%%%%
%%%%%%%%%%%%%%%%%%%%%%%%%%%%%%%%%%%%%%%%%%%%%%%%%%%%%%%%%%%%%%%%%%%%%%%%%%%%%%%%%%%%%%%%%%%%%%%%%%%%%%%%%%%%%%%%%%%%%%%%%%%%%%%%%%%%%%%%%%%%%%%%%%%%%%%%%%%%%%%%%%%%%%%%%%%%%%%%%%%%%%%%%%%%%%%%%%%%%%%%%%%%%%%%%%%%%%%%%%%%%%%%%

\subsection{Change of variable}

Employing the strategy outlined in  \cite{tao2009global, tao2011global,tao2014energy, ke2018note, fuest2022global, giesselmann2018existence, tao2020global,li2021boundedness, surulescu2021does, fuest2023cross, heihoff2023global}, we change the variable for $c_1$ to achieve an equivalent system, which expresses the first equation of \eqref{model1} in the divergence form. To this end, we substitute
 \begin{equation}
 \label{cov1}
   z := c_1 e^{- \chi \tau} \quad \text{defined on}~ \overline{\Omega} \times [0, T_{\max}).
 \end{equation}
The associated initial data is $z_0 := c_{10} e^{- \chi \tau_0}$ defined on $\overline{\Omega}$. Using \eqref{cov1} we can rewrite \eqref{model1} in the following form
\begin{equation}
  \begin{cases}
    \label{sys1}
    \partial_t z= e^{- \chi \tau} \nabla \cdot (e^{\chi \tau} \mathbb{D} \nabla z) + e^{- \chi \tau} \nabla \cdot (z e^{\chi \tau} \nabla \cdot \mathbb{D}) - \alpha_1 \frac{z}{1 + z e^{\chi \tau}} + \alpha_2 \frac{c_2 e^{-\chi \tau}}{1 + c_2}                                                                                                                                                                \\[5pt]
    \hspace{1cm} + \beta z (1 - z e^{\chi \tau} - c_2 -\tau)+ \delta \chi e^{\chi \tau} z^2 \tau + \mu \chi \tau z - \frac{\chi z c_2}{1 + c_2},     &   x \in \Omega, t > 0,               \\[5pt]
    \partial_t c_2 = \alpha_1 \frac{z e^{\chi \tau}}{1 + z e^{\chi \tau}} - \alpha_2 \frac{c_2}{1 + c_2} ,                                           &   x \in \Omega, t > 0,               \\[5pt]
    \partial_t \tau = -\delta z  e^{\chi \tau} \tau - \mu \tau + \frac{c_2}{1 + c_2},                                                                &   x \in \Omega, t > 0,               \\[5pt]
    (\mathbb{D} \nabla z) \cdot \nu = - z (\nabla \cdot \mathbb{D}) \cdot \nu = 0,                                            &   x \in \partial \Omega, t > 0,      \\[5pt]
    z(\cdot, 0) = z_0, \quad c_2(\cdot, 0) = c_{20},     &   x \in \Omega.
   \end{cases}
\end{equation}
The rationale for using the transformation \eqref{cov1} resulting in the construction of the equivalent system in \eqref{sys1} is that it can eliminate the gradient involved with taxis terms by integrating it into the diffusion operator for $z$. As pointed out in \cite{tao2011global, tao2014energy, li2021boundedness, surulescu2021does}, the models \eqref{model1} and \eqref{sys1} are equivalent within the framework of classical solutions. Hence, such a transformation should not pose any difficulty in translating results, as it involves only a switching between the variables $z$ and $c_1$. We will start the analysis by deriving some bounds for this equivalent system.

%%%%%%%%%%%%%%%%%%%%%%%%%%%%%%%%%%%%%%%%%%%%%%%%%%%%%%%%%%%%%%%%%%%%%%%%%%%%%%%%%%%%%%%%%%%%%%%%%%%%%%%%%%%%%%%%%%%%%%%%%%%%%%%%%%%%%%%%%%%%%%%%%%%%%%%%%%%%%%%%%%%%%%%%%%%%%%%%%%%%%%%%%%%%%%%%%%%%%%%%%%%%%%%%%%%%%%%%%%%%%%%%%
%%%%%%%%%%%%%%%%%%%%%%%%%%%%%%%%%%%%%%%%%%%%%%%%%%%%%%%%%%%%%%%%%%%%%%%%%%%%%%%%%%%%%%%%%%%%%%%%%%%%%%%%%%%%%%%%%%%%%%%%%%%%%%%%%%%%%%%%%%%%%%%%%%%%%%%%%%%%%%%%%%%%%%%%%%%%%%%%%%%%%%%%%%%%%%%%%%%%%%%%%%%%%%%%%%%%%%%%%%%%%%%%%

\begin{corollary}
  \label{L1_z:cor}
  There exists a $C > 0$ such that
  \begin{equation}
    \label{zL1_bound}
    \int_\Omega z(\cdot, t) \leq C =: z^\ast \quad \text{for all}~ t \in (0, T_{\max}).
  \end{equation}
\end{corollary}
\begin{proof}
  Since $\int_\Omega z = \int_\Omega  e^{-\chi \tau} c_1\leq \int_\Omega c_1$ for all $t \in (0, T_{\max})$, \eqref{zL1_bound}  follows directly from \eqref{l1c1} in view of \eqref{tau_bound}.
\end{proof}

%%%%%%%%%%%%%%%%%%%%%%%%%%%%%%%%%%%%%%%%%%%%%%%%%%%%%%%%%%%%%%%%%%%%%%%%%%%%%%%%%%%%%%%%%%%%%%%%%%%%%%%%%%%%%%%%%%%%%%%%%%%%%%%%%%%%%%%%%%%%%%%%%%%%%%%%%%%%%%%%%%%%%%%%%%%%%%%%%%%%%%%%%%%%%%%%%%%%%%%%%%%%%%%%%%%%%%%%%%%%%%%%%
%%%%%%%%%%%%%%%%%%%%%%%%%%%%%%%%%%%%%%%%%%%%%%%%%%%%%%%%%%%%%%%%%%%%%%%%%%%%%%%%%%%%%%%%%%%%%%%%%%%%%%%%%%%%%%%%%%%%%%%%%%%%%%%%%%%%%%%%%%%%%%%%%%%%%%%%%%%%%%%%%%%%%%%%%%%%%%%%%%%%%%%%%%%%%%%%%%%%%%%%%%%%%%%%%%%%%%%%%%%%%%%%%
\subsection{\texorpdfstring
{Global classical solvability to \eqref{model1} in two dimensions}
{Classical solutions2d}}

\subsubsection{A quasi-energy inequality}

This subsection aims to derive a quasi-energy inequality for \eqref{model1}. To this end, we will begin with a standard testing procedure, through which the haptotactic term in the first equation of \eqref{model1} reduces to an inner product of gradients.

%%%%%%%%%%%%%%%%%%%%%%%%%%%%%%%%%%%%%%%%%%%%%%%%%%%%%%%%%%%%%%%%%%%%%%%%%%%%%%%%%%%%%%%%%%%%%%%%%%%%%%%%%%%%%%%%%%%%%%%%%%%%%%%%%%%%%%%%%%%%%%%%%%%%%%%%%%%%%%%%%%%%%%%%%%%%%%%%%%%%%%%%%%%%%%%%%%%%%%%%%%%%%%%%%%%%%%%%%%%%%%%%%
%%%%%%%%%%%%%%%%%%%%%%%%%%%%%%%%%%%%%%%%%%%%%%%%%%%%%%%%%%%%%%%%%%%%%%%%%%%%%%%%%%%%%%%%%%%%%%%%%%%%%%%%%%%%%%%%%%%%%%%%%%%%%%%%%%%%%%%%%%%%%%%%%%%%%%%%%%%%%%%%%%%%%%%%%%%%%%%%%%%%%%%%%%%%%%%%%%%%%%%%%%%%%%%%%%%%%%%%%%%%%%%%%

\begin{lemma}
\label{c1logbound}
There exists a $C > 0$ such that
\begin{align}
  \label{logineq0}
  \frac{d}{dt} \int_\Omega c_1 \ln c_1 + \frac{1}{4}\int_\Omega \frac{\nabla c_1 \cdot \mathbb{D} \nabla c_1}{c_1} + \frac{\beta}{4} \int_\Omega c_1^2 (\ln 2 + c_1) \leq  \chi \int_\Omega \nabla c_1 \cdot \mathbb{D} \nabla \tau  + C
\end{align}
holds for all $t > 0$.
\end{lemma}
\begin{proof}
  We test the first equation of \eqref{model1} with $\ln c_1$. This results in
  \begin{align}
    \label{logineq1}
    \nonumber & \frac{d}{dt} \int_\Omega c_1 \ln c_1 = \int_\Omega \partial_t c_{1} \ln c_1 + \int_\Omega \partial_t c_{1}                                                                                   \\[5pt]
    \nonumber & = \int_\Omega \Big\{\nabla \cdot (\mathbb{D} \nabla c_1 + c_1 \nabla \cdot \mathbb{D}) - \chi \nabla \cdot (c_1 \mathbb{D} \nabla \tau)- \alpha_1  \tfrac{c_{1}}{1 + c_{1}}  + \alpha_2 \tfrac{c_{2}}{1 + c_{2}} \\[5pt] \nonumber & + \beta c_{1}(1 - c_{1} - c_{2} - \tau)\Big\}\ln c_1 + \int_\Omega \{- \alpha_1  \tfrac{c_{1}}{1 + c_{1}}  + \alpha_2 \tfrac{c_{2}}{1 + c_{2}} + \beta c_{1}(1 - c_{1} - c_{2} - \tau)\} \\[5pt]
    \nonumber & = - \int_\Omega \frac{\nabla c_1 \cdot \mathbb{D} \nabla c_1}{c_1} - \int_\Omega \nabla c_1 \cdot (\nabla \cdot \mathbb{D}) + \chi \int_\Omega \nabla c_1 \cdot \mathbb{D} \nabla \tau - \alpha_1  \int_\Omega \tfrac{c_{1} \ln c_1}{1 + c_{1}}  + \alpha_2 \int_\Omega \tfrac{c_{2} \ln c_1 }{1 + c_{2}}  \\[5pt]
    & + \beta \int_\Omega c_{1}(1 - c_{1} - c_{2} - \tau) \ln c_1 - \alpha_1 \int_\Omega \tfrac{c_{1}}{1 + c_{1}} + \alpha_2 \int_\Omega \tfrac{c_{2}}{1 + c_{2}} + \beta\int_\Omega  c_{1}(1 - c_{1} - c_{2} - \tau)
  \end{align}
  for all $t > 0$. The first three terms on the right-hand side of \eqref{logineq1} have been obtained using integration by parts. We will first handle the term resulting from non-Fickian diffusion. From \eqref{assumption2} we have the divergence estimate with an exponent $\lambda \in [\tfrac{1}{2}, \tfrac{2}{3}]$ such that $\frac{\lambda}{1 - \lambda} \leq 2$, and a constant $A > 0$ which, in conjunction with Young's inequality, enables the following estimate:
  \begin{align}
  \label{nlogineq2}
    \nonumber \left|\int_\Omega \nabla c_1 \cdot (\nabla \cdot \mathbb{D})\right| & \leq A \int_\Omega (\nabla c_1 \cdot \mathbb{D} \nabla c_1)^\lambda + A                                     \\[5pt]
    & \leq \frac{1}{2} \int_\Omega \frac{\nabla c_1 \cdot \mathbb{D} \nabla c_1}{c_1} + 2^{\frac{\lambda}{1 - \lambda}} A^{\frac{1}{1 - \lambda}} \int_\Omega c_1^{\frac{\lambda}{1 - \lambda}}
  \end{align}
 for all $t > 0$. As we have pointed out that $\frac{\lambda}{1 - \lambda} \leq 2$ since $\lambda \in [\tfrac{1}{2}, \tfrac{2}{3}]$, we have a constant $C_1 > 0$ such that
 \begin{align}
   \label{nlogineq3}
   2^{\frac{\lambda}{1 - \lambda}} A^{\frac{1}{1 - \lambda}} c_1^{\frac{\lambda}{1 - \lambda}} - \frac{\beta}{2} c_1^2 \ln c_1 \leq C_1
 \end{align}
 for all $c_1 \geq 0$. Inserting \eqref{nlogineq2} and \eqref{nlogineq3} in \eqref{logineq1} results in
 \begin{align}
    \label{nlogineq4}
    \nonumber & \frac{d}{dt} \int_\Omega c_1 \ln c_1 \leq - \frac{1}{2} \int_\Omega \frac{\nabla c_1 \cdot \mathbb{D} \nabla c_1}{c_1}  + \chi \int_\Omega \nabla c_1 \cdot \mathbb{D} \nabla \tau - \alpha_1  \int_\Omega \tfrac{c_{1} \ln c_1}{1 + c_{1}}  + \alpha_2 \int_\Omega \tfrac{c_{2} \ln c_1 }{1 + c_{2}}                                                                        \\[5pt]
    & \nonumber - \frac{\beta}{2} \int_\Omega c_1^2 \ln c_1  + \beta \int_\Omega c_{1}(1 - c_{2} - \tau) \ln c_1 - \alpha_1 \int_\Omega \tfrac{c_{1}}{1 + c_{1}} + \alpha_2 \int_\Omega \tfrac{c_{2}}{1 + c_{2}} \\[5pt]
    & + \beta\int_\Omega  c_{1}(1 - c_{1} - c_{2} - \tau) + C_1{ |\Omega|}
 \end{align}
 for all $t > 0$. Now, by invoking \cite[Lemma 4.2]{stinner2014global}, we have a constant $C_2>0$ such that
\begin{align}
    \label{ourr}
    \frac{\beta}{2} \int_\Omega \left\{c_1 \ln c_1 - c_1^2 \ln c_1\right\} \leq - \frac{\beta}{4} \int_\Omega c_1^2 \ln(2 + c_1) + C_2,
\end{align}
for all $t > 0$. We next estimate the term $\frac{\beta}{2} \int_\Omega c_1 \ln c_1$ appearing on the right-hand side of \eqref{nlogineq4}. To this end, we use the elementary inequality $s \ln s \leq s^{\frac{3}{2}}$ for all $s \geq 0$, which yields
\begin{align*}
    \frac{\beta}{2} \int_\Omega c_1 \ln c_1 \leq \frac{\beta}{2}\int_\Omega c_1^{\frac{3}{2}},
\end{align*}
for all $t > 0$. Since the spatial dimension is $2$, the Gagliardo-Nirenberg inequality (cf. \cite[(3.19) with $q=\frac{3}{2}$]{tao2011global}), in conjunction with Young's inequality and the uniform $L^1$ bound for $c_1$ established in \eqref{l1c1}, yields
\begin{align}
    \label{gagr1}
    \nonumber \frac{\beta}{2} \int_\Omega c_1^{\frac{3}{2}} & \leq  \frac{\beta}{2}  C^{(2)}_{GN}\|\nabla \sqrt{c_1}\|_{L^2(\Omega)} \|\sqrt{c_1}\|^2_{L^2(\Omega)} + \frac{\beta}{2}  C^{(2)}_{GN} \|\sqrt{c_1}\|^3_{L^2(\Omega)}   \\[5pt]
    \nonumber & \leq \frac{\beta}{2} C^{(2)}_{GN}\left(\frac{1}{4} \int_\Omega \frac{|\nabla c_1|^2}{c_1}\right)^{\frac{1}{2}} \left(\int_\Omega c_1\right) + \frac{\mu_1}{2}  C^{(2)}_{GN} \left(\int_\Omega c_1\right)^{\frac{3}{2}} \\[5pt]
    \nonumber & \leq  \frac{\beta}{2}  C^{(2)}_{GN} M_1 C_2 \left(\frac{1}{4} \int_\Omega \frac{|\nabla c_1|^2}{c_1}\right)^{\frac{1}{2}} + \frac{\beta}{2}  C^{(2)}_{GN} M_1^{\frac{3}{2}} \\[5pt]
    & \leq  \frac{\epsilon_1}{4} \int_\Omega \frac{|\nabla c_1|^2}{c_1} + \frac{1}{4\epsilon_1}  \left(\frac{\beta}{2}  C^{(2)}_{GN} M_1 \right)^2 +  \frac{\beta}{2}  C^{(2)}_{GN} M_1^{\frac{3}{2}}
\end{align}
for all $t > 0$. By virtue of \eqref{assumption0}, we infer from \eqref{gagr1} upon setting
$\epsilon_1 = D_1$, that the following holds:
\begin{align}
    \label{gagr2}
    \nonumber \frac{\beta}{2} \int_\Omega c_1^{\frac{3}{2}} & \leq \frac{\epsilon_1}{4 D_1} \int_\Omega \frac{\nabla c_1 \cdot { D_1} \nabla c_1}{c_1} + C_3 \\[7pt]
    & = \frac{1}{4} \int_\Omega \frac{\nabla c_1 \cdot \mathbb{D} \nabla c_1}{c_1} + C_3 \quad \text{ (by setting $\epsilon_1:=1$)}
\end{align}
for all $t > 0$.
Using the non-negativity of solution components, \eqref{l1c1}, \eqref{tau_bound}, \eqref{c2_bound}, the fact that $\frac{c_i}{1 + c_i} < 1$, for $c_i \geq 0,\ i \in \{1, 2\}$, and the inequalities $s\ln s \geq -\tfrac{1}{e}$ for all $s > 0$ and $\ln s \leq \tfrac{1}{e} s$ for $s > 1$, we have
  \begin{align}
    \label{logineq3}  & - \alpha_1  \int_\Omega \frac{c_{1} \ln c_1}{1 + c_{1}} \leq \frac{\alpha_1}{e} |\Omega|,                                                                \\[5pt]
    \label{logineq4}  &  \alpha_2 \int_\Omega \frac{c_{2} }{1 + c_{2}} \ln c_1 \leq \alpha_2 \int_{\{c_1 > 1\}} \ln c_1 \leq \alpha_2 \int_{\{c_1 > 1\}} c_1 \leq \alpha_2 M_1,           \\[5pt]
    \label{logineq5}  & - \beta \int_\Omega c_1 c_2 \ln c_1 \leq \frac{\beta}{e} \int_\Omega c_2 \leq \frac{1}{e} \beta M_{\tau}|\Omega|,                                        \\[5pt]
    \label{logineq6}  & - \beta \int_\Omega c_1 \tau \ln c_1 \leq \frac{\beta}{e} \int_\Omega \tau \leq \frac{1}{e} \beta M_{\tau} |\Omega|,                                     \\[5pt]
    \label{logineq7}  & - \alpha_1 \int_\Omega \tfrac{c_{1}}{1 + c_{1}} + \alpha_2 \int_\Omega \tfrac{c_{2}}{1 + c_{2}} + \beta\int_\Omega  c_{1}(1 - c_{1} - c_{2} - \tau) \leq \alpha_2 |\Omega| + \beta M_1
  \end{align}
  for all $t > 0$. Combining \eqref{ourr} and \eqref{gagr2}-\eqref{logineq7} with \eqref{nlogineq4} yields \eqref{logineq0}.
\end{proof}

%%%%%%%%%%%%%%%%%%%%%%%%%%%%%%%%%%%%%%%%%%%%%%%%%%%%%%%%%%%%%%%%%%%%%%%%%%%%%%%%%%%%%%%%%%%%%%%%%%%%%%%%%%%%%%%%%%%%%%%%%%%%%%%%%%%%%%%%%%%%%%%%%%%%%%%%%%%%%%%%%%%%%%%%%%%%%%%%%%%%%%%%%%%%%%%%%%%%%%%%%%%%%%%%%%%%%%%%%%%%%%%%%
%%%%%%%%%%%%%%%%%%%%%%%%%%%%%%%%%%%%%%%%%%%%%%%%%%%%%%%%%%%%%%%%%%%%%%%%%%%%%%%%%%%%%%%%%%%%%%%%%%%%%%%%%%%%%%%%%%%%%%%%%%%%%%%%%%%%%%%%%%%%%%%%%%%%%%%%%%%%%%%%%%%%%%%%%%%%%%%%%%%%%%%%%%%%%%%%%%%%%%%%%%%%%%%%%%%%%%%%%%%%%%%%%
The integral on the right-hand side in \eqref{logineq0}, which originates from haptotactic cross-diffusion, can be handled by studying the time evolution of the integral $ \int_\Omega \frac{\nabla \tau \cdot \mathbb{D} \nabla \tau}{\tau}$.

\begin{lemma}
\label{tauevobound}
The following holds true
\begin{align}
  \label{tauinq0}
  \nonumber  \dfrac{1}{2} \frac{d}{dt} \int_\Omega \frac{\nabla \tau \cdot \mathbb{D} \nabla \tau}{\tau} & + \delta \int_\Omega \nabla c_1 \cdot \mathbb{D} \nabla \tau  + \frac{\delta}{2} \int_\Omega \frac{c_1 \nabla \tau \cdot \mathbb{D} \nabla \tau}{\tau} \\[5pt]
   & + \frac{\mu}{2} \int_\Omega \dfrac{\nabla \tau \cdot \mathbb{D} \nabla \tau}{\tau} \leq \frac{D_2^2}{2 D_1} \int_\Omega \frac{|\nabla c_2|^2}{c_2}
\end{align}
for all $t > 0$.
\end{lemma}
\begin{proof}
  From the third equation of \eqref{model1} we get
  \begin{align}
  \label{tauinq1}
    \nonumber & \dfrac{1}{2} \frac{d}{dt} \int_\Omega \frac{\nabla \tau \cdot \mathbb{D} \nabla \tau}{\tau} = \int_\Omega \dfrac{\nabla \tau_t \cdot \mathbb{D} \nabla \tau}{\tau} - \frac{1}{2} \int_\Omega \frac{\tau_t (\nabla \tau \cdot  \mathbb{D} \nabla \tau)}{\tau^2} \\[5pt]
    \nonumber & = - \delta \int_\Omega \nabla c_1 \cdot \mathbb{D} \nabla \tau - \delta \int_\Omega \frac{c_1 \nabla \tau \cdot \mathbb{D} \nabla \tau}{\tau} - \mu \int_\Omega \dfrac{\nabla \tau \cdot \mathbb{D} \nabla \tau}{\tau} + \int_\Omega \frac{\nabla c_2 \cdot \mathbb{D} \nabla \tau}{\tau} \cdot \frac{1}{(1 + c_2)^2} \\[5pt]
    \nonumber  & + \frac{\delta}{2} \int_\Omega \frac{c_1 \nabla \tau \cdot \mathbb{D} \nabla \tau}{\tau} + \frac{\mu}{2} \int_\Omega \frac{\nabla \tau \cdot \mathbb{D} \nabla \tau}{\tau} - \frac{1}{2} \int_\Omega \frac{\nabla \tau \cdot \mathbb{D} \nabla \tau}{\tau^2} \cdot \frac{c_2}{1 + c_2} \\[5pt]
    \nonumber  & = - \delta \int_\Omega \nabla c_1 \cdot \mathbb{D} \nabla \tau - \frac{\delta}{2} \int_\Omega \frac{c_1 \nabla \tau \cdot \mathbb{D} \nabla \tau}{\tau} - \frac{\mu}{2} \int_\Omega \dfrac{\nabla \tau \cdot \mathbb{D} \nabla \tau}{\tau} + \int_\Omega \frac{\nabla c_2 \cdot \mathbb{D} \nabla \tau}{\tau} \cdot \frac{1}{(1 + c_2)^2} \\[5pt]
    & - \frac{1}{2} \int_\Omega \frac{\nabla \tau \cdot \mathbb{D} \nabla \tau}{\tau^2} \cdot \frac{c_2}{1 + c_2}
  \end{align}
   for all $t > 0$. To handle the last two terms in \eqref{tauinq1} we will use \eqref{assumption0}, Young's inequality and the fact that $\frac{1}{(1 + s)^3} < 1$ for all $s > 0$
  % \begin{align*}
   %  & \int_\Omega \frac{\nabla c_2 \cdot \mathbb{D} \nabla \tau}{\tau} \cdot \frac{1}{(1 + c_2)^2}- \frac{1}{2} \int_\Omega \frac{\nabla \tau \cdot \mathbb{D} \nabla \tau}{\tau^2} \cdot \frac{c_2}{1 + c_2}  \\[5pt]
   %  &  \leq \frac{D_2}{2 D_1} \int_\Omega \frac{|\nabla c_2|^2}{c_2} \cdot \frac{1}{(1 + c_2)^3} + \frac{D_1}{2 D_2} \int_\Omega \frac{{\cb \mathbb D}\nabla \tau \cdot \mathbb{D} \nabla \tau}{\tau^2} \cdot \frac{c_2}{1 + c_2}  - \frac{1}{2} \int_\Omega \frac{\nabla \tau \cdot \mathbb{D} \nabla \tau}{\tau^2} \cdot \frac{c_2}{1 + c_2}  \\[5pt]
   %  & \leq \frac{D_2}{2 D_1} \int_\Omega \frac{|\nabla c_2|^2}{c_2} +  \frac{D_1}{2 D_2} \cdot D_2 \int_\Omega \frac{|\nabla \tau|^2 }{\tau^2} \cdot \frac{c_2}{1 + c_2} - \frac{D_1}{2} \int_\Omega \frac{|\nabla \tau|^2 }{\tau^2} \cdot \frac{c_2}{1 + c_2}
   %\end{align*}
  % for all $t > 0$. %{\cred How do you get rid of the blue $\mathbb D$ above? - {\color{blue} Please see below.}} Inserting \eqref{tauinq2} in \eqref{tauinq1} results in \eqref{tauinq0}.
   \begin{align}
   \label{tauinq2}
   \nonumber & \int_\Omega\frac{\nabla c_2 \cdot \mathbb{D}\nabla\tau}{\tau}\cdot\frac{1}{(1+c_2)^2}-\frac12
   \int_\Omega \frac{\nabla\tau\cdot\mathbb{D}\nabla\tau}{\tau^2}\cdot\frac{c_2}{1+c_2}\\[5pt]
   \nonumber&  \leq\int_\Omega\frac{
  (\nabla c_2\cdot\mathbb{D}\nabla c_2)^{\frac{1}{2}}(\nabla\tau\cdot\mathbb{D}\nabla\tau)^{\frac12}}{\tau}\cdot\frac{1}{(1+c_2)^2}-\frac12\int_\Omega \frac{\nabla\tau\cdot\mathbb{D}\nabla\tau}{\tau^2}\cdot\frac{c_2}{1+c_2} \\[5pt]
  \nonumber &  \leq \int_\Omega \frac{\sqrt{D_2} |\nabla c_2| (\nabla\tau\cdot\mathbb{D}\nabla\tau)^{\frac12}}{\tau} \cdot \frac{1}{(1 + c_2)^2} -\frac12\int_\Omega \frac{\nabla\tau\cdot\mathbb{D}\nabla\tau}{\tau^2}\cdot\frac{c_2}{1+c_2} \\[5pt]
  \nonumber &  \leq \frac{D_2^2}{2D_1} \int_\Omega \frac{|\nabla c_2|^2}{c_2} \cdot \frac{1}{(1 + c_2)^3} + \frac{D_1}{2 D_2} \int_\Omega \frac{\nabla \tau \cdot \mathbb{D} \nabla \tau}{\tau^2} \cdot \frac{c_2}{1 + c_2}  -\frac12\int_\Omega \frac{\nabla\tau\cdot\mathbb{D}\nabla\tau}{\tau^2}\cdot\frac{c_2}{1+c_2} \\[5pt] 
  &  \leq \frac{D_2^2}{2D_1} \int_\Omega \frac{|\nabla c_2|^2}{c_2} \cdot \frac{1}{(1 + c_2)^3} + \frac{D_1}{2 D_2} \cdot D_2\int_\Omega \frac{|\nabla \tau|^2}{\tau^2} \cdot \frac{c_2}{1 + c_2}  -\frac{D_1}{2}\int_\Omega \frac{|\nabla\tau|^2}{\tau^2}\cdot\frac{c_2}{1+c_2}
\end{align}
 for all $t > 0$. Inserting \eqref{tauinq2} in \eqref{tauinq1} results in \eqref{tauinq0}.
\end{proof}

%%%%%%%%%%%%%%%%%%%%%%%%%%%%%%%%%%%%%%%%%%%%%%%%%%%%%%%%%%%%%%%%%%%%%%%%%%%%%%%%%%%%%%%%%%%%%%%%%%%%%%%%%%%%%%%%%%%%%%%%%%%%%%%%%%%%%%%%%%%%%%%%%%%%%%%%%%%%%%%%%%%%%%%%%%%%%%%%%%%%%%%%%%%%%%%%%%%%%%%%%%%%%%%%%%%%%%%%%%%%%%%%%
%%%%%%%%%%%%%%%%%%%%%%%%%%%%%%%%%%%%%%%%%%%%%%%%%%%%%%%%%%%%%%%%%%%%%%%%%%%%%%%%%%%%%%%%%%%%%%%%%%%%%%%%%%%%%%%%%%%%%%%%%%%%%%%%%%%%%%%%%%%%%%%%%%%%%%%%%%%%%%%%%%%%%%%%%%%%%%%%%%%%%%%%%%%%%%%%%%%%%%%%%%%%%%%%%%%%%%%%%%%%%%%%%

The ill-signed integral on the right-hand side in \eqref{tauinq0}, which originates from the production term $\frac{c_2}{1 + c_2}$ in the third equation of \eqref{model1}, can be handled by studying the time evolution of the Dirichlet integral $ \int_\Omega \frac{|\nabla c_2|^2}{c_2}$.

\begin{lemma}
\label{c2evobound}
 The following holds true:
\begin{align}
  \label{c2inq0}
   \dfrac{1}{2} \frac{d}{dt} \int_\Omega \frac{|\nabla c_2|^2}{c_2} \leq \frac{\alpha_1}{2 D_1} \int_\Omega \frac{\nabla c_1 \cdot \mathbb{D} \nabla c_1}{c_1} + \frac{\alpha_2}{2} \int_\Omega \frac{|\nabla c_2|^2}{c_2}
\end{align}
for all $t > 0$.
\end{lemma}

\begin{proof}
  From the second equation of \eqref{model1} we have
  \begin{align}
  \label{c2inq1}
    \nonumber & \dfrac{1}{2} \frac{d}{dt} \int_\Omega \frac{|\nabla c_2|^2}{c_2} = \int_\Omega \frac{\nabla \partial_tc_{2} \cdot \nabla c_2}{c_2} - \frac{1}{2} \int_\Omega \frac{|\nabla c_2|^2}{c_2^2} \partial_tc_{2} \\[5pt]
    \nonumber & = \alpha_1 \int_\Omega \frac{\nabla c_1 \cdot \nabla c_2}{c_2 (1 + c_1)^2} - \alpha_2 \int_\Omega \frac{|\nabla c_2|^2}{c_2} \cdot \frac{1}{(1 + c_2)^2} - \frac{\alpha_1}{2} \int_\Omega \frac{|\nabla c_2|^2}{c_2^2} \cdot \frac{c_1}{1 + c_1} \\[5pt]
    & + \frac{\alpha_2}{2} \int_\Omega \frac{|\nabla c_2|^2}{c_2} \cdot \frac{1}{1 + c_2}
  \end{align}
   for all $t > 0$. Applying Young's inequality to the first term on the right hand side of \eqref{c2inq1}, and subsequently using \eqref{assumption0}, we obtain 
   \begin{align}
     \label{c2inq2}
     \nonumber \alpha_1 \int_\Omega \frac{\nabla c_1 \cdot \nabla c_2}{c_2 (1 + c_1)^2} & \leq \frac{\alpha_1}{2} \int_\Omega \frac{|\nabla c_1|^2}{c_1} \cdot \frac{1}{(1 + c_1)^3} + \frac{\alpha_1}{2}\int_\Omega  \frac{|\nabla c_2|^2}{c_2^2} \cdot \frac{c_1}{1 + c_1} \\[5pt]
     \nonumber & \leq \frac{\alpha_1}{2} \int_\Omega \frac{\nabla c_1 \cdot \nabla c_1}{c_1}  + \frac{\alpha_1}{2}\int_\Omega  \frac{|\nabla c_2|^2}{c_2^2} \cdot \frac{c_1}{1 + c_1} \\[5pt]
     & \leq \frac{\alpha_1}{2 D_1} \int_\Omega \frac{\nabla c_1 \cdot \mathbb{D} \nabla c_1}{c_1} + \frac{\alpha_1}{2}\int_\Omega  \frac{|\nabla c_2|^2}{c_2^2} \cdot \frac{c_1}{1 + c_1}
   \end{align}
    for all $t > 0$. Inserting \eqref{c2inq2} in \eqref{c2inq1} and dropping a negative term yields \eqref{c2inq0}.
\end{proof}

%%%%%%%%%%%%%%%%%%%%%%%%%%%%%%%%%%%%%%%%%%%%%%%%%%%%%%%%%%%%%%%%%%%%%%%%%%%%%%%%%%%%%%%%%%%%%%%%%%%%%%%%%%%%%%%%%%%%%%%%%%%%%%%%%%%%%%%%%%%%%%%%%%%%%%%%%%%%%%%%%%%%%%%%%%%%%%%%%%%%%%%%%%%%%%%%%%%%%%%%%%%%%%%%%%%%%%%%%%%%%%%%%
%%%%%%%%%%%%%%%%%%%%%%%%%%%%%%%%%%%%%%%%%%%%%%%%%%%%%%%%%%%%%%%%%%%%%%%%%%%%%%%%%%%%%%%%%%%%%%%%%%%%%%%%%%%%%%%%%%%%%%%%%%%%%%%%%%%%%%%%%%%%%%%%%%%%%%%%%%%%%%%%%%%%%%%%%%%%%%%%%%%%%%%%%%%%%%%%%%%%%%%%%%%%%%%%%%%%%%%%%%%%%%%%%

We are now in a position to develop an energy-like inequality from the outcomes of Lemmas \ref{c1logbound}, \ref{tauevobound}, and \ref{c2evobound}, expressed as follows.

\begin{lemma}
\label{entropy_lemma}
Let $T > 0$. Then there exists a $C(T) > 0$ such that
\begin{align}
  \label{entropy0}
  \nonumber & \sup_{t \in (0, \widehat{T})}  \left\{\int_\Omega  \left(c_1\ln c_1 + \frac{1}{e}\right) + \frac{\chi}{2 \delta} \int_\Omega \frac{\nabla \tau \cdot \mathbb{D} \nabla \tau}{\tau} + \frac{D_1}{8\alpha_1}  \int_\Omega \frac{|\nabla c_2|^2}{c_2}\right\}  + \frac{1}{8}  \int_0^{\widehat{T}} \int_\Omega \frac{\nabla c_1 \cdot \mathbb{D} \nabla c_1}{c_1}\\[5pt]
  & +  \frac{\beta}{4} \int_0^{\widehat{T}} \int_\Omega c_1^2 (\ln 2 + c_1) + \frac{\chi}{2} \int_0^{\widehat{T}} \int_\Omega \frac{c_1 \nabla \tau \cdot \mathbb{D} \nabla \tau}{\tau} + \frac{\mu \chi}{2 \delta} \int_0^{\widehat{T}} \int_\Omega \dfrac{\nabla \tau \cdot \mathbb{D} \nabla \tau}{\tau} \leq C(T),
\end{align}
where $\widehat{T} := \min\{T, T_{\max}\}$.
\end{lemma}

\begin{proof}
 From \eqref{logineq0}, \eqref{tauinq0} and \eqref{c2inq0} it follows that
 \begin{align*}
   \nonumber & \dfrac{d}{dt} \left\{\int_\Omega  \left(c_1\ln c_1 + \frac{1}{e}\right) + \frac{\chi}{2 \delta} \int_\Omega \frac{\nabla \tau \cdot \mathbb{D} \nabla \tau}{\tau} + \frac{D_1}{8\alpha_1}  \int_\Omega \frac{|\nabla c_2|^2}{c_2}\right\}                                                                                   \\[5pt]
   &  + \frac{1}{8}  \int_\Omega \frac{\nabla c_1 \cdot \mathbb{D} \nabla c_1}{c_1}  + \frac{\beta}{4} \int_\Omega c_1^2 (\ln 2 + c_1) + \frac{\chi}{2} \int_\Omega \frac{c_1 \nabla \tau \cdot \mathbb{D} \nabla \tau}{\tau} \\[5pt]
   &  + \frac{\mu \chi}{2 \delta} \int_\Omega \dfrac{\nabla \tau \cdot \mathbb{D} \nabla \tau}{\tau} \leq \left(\frac{\alpha_2 D_1}{8 \alpha_1} + \frac{D_2^2 \chi}{2 D_1 \delta}\right) \int_\Omega \frac{|\nabla c_2|^2}{c_2} + C_1
 \end{align*}
 for all $t > 0$. Next, we define the following two nonnegative functions on $(0, \infty)$:
 \begin{align*}
   \mathcal{F}(t) & := \int_\Omega  \left(c_1\ln c_1 + \frac{1}{e}\right) + \frac{\chi}{2 \delta} \int_\Omega \frac{\nabla \tau \cdot \mathbb{D} \nabla \tau}{\tau} + \frac{D_1}{8\alpha_1}  \int_\Omega \frac{|\nabla c_2|^2}{c_2}, \\[5pt]
   \mathcal{D}(t) & := \frac{1}{8}  \int_\Omega \frac{\nabla c_1 \cdot \mathbb{D} \nabla c_1}{c_1} + \frac{\beta}{ 4} \int_\Omega c_1^2 (\ln 2 + c_1) + \frac{\chi}{2} \int_\Omega \frac{c_1 \nabla \tau \cdot \mathbb{D} \nabla \tau}{\tau} + \frac{\mu \chi}{2\delta} \int_\Omega \dfrac{\nabla \tau \cdot \mathbb{D} \nabla \tau}{\tau}.
 \end{align*}
 Using the above functions and applying the inequality $-s \ln s \leq \frac{1}{e}$ for all $s > 0$, we arrive at the following ordinary differential inequality
 \begin{align}
 \label{entropy2}
   \nonumber \mathcal{F}'(t) + \mathcal{D}(t) & \leq \widehat{\eta} \int_\Omega \frac{|\nabla c_2|^2}{c_2} + C_1 \leq \eta \mathcal{F}(t) - \eta \int_\Omega c_1 \ln c_1 - \frac{\eta}{e} - \frac{\eta \chi}{2\delta} \int_\Omega \frac{\nabla \tau \cdot \mathbb{D} \nabla \tau}{\tau} + C_1  \\[5pt]
   & \leq \eta \mathcal{F}(t) + \frac{\eta}{e} |\Omega| + C_1
 \end{align}
 for all $t > 0$, where $\widehat{\eta} = \frac{\alpha_2 D_1}{8 \alpha_1} + \frac{D_2^2 \chi}{2 D_1 \delta}$ and $\eta = \frac{8\alpha_1}{D_1} \widehat{\eta}$. Since $\mathcal{F}(t) > 0$ and $\mathcal{D}(t) > 0$, we can apply the Gronwall's lemma to \eqref{entropy2}, thus entailing
 \begin{align}
   \label{entropy3}
   \mathcal{F}(t) \leq \mathcal{F}(0) e^{\eta t} +  \frac{\eta}{e} |\Omega| t + C_1 t \leq C_2 \quad \text{for all}~ t \in (0, {\widehat{T}}).
 \end{align}
A direct time integration of \eqref{entropy2} over the interval $(0, {\widehat{T}})$ in view of \eqref{entropy3} allows us to obtain \eqref{entropy0}.
\end{proof}
%%%%%%%%%%%%%%%%%%%%%%%%%%%%%%%%%%%%%%%%%%%%%%%%%%%%%%%%%%%%%%%%%%%%%%%%%%%%%%%%%%%%%%%%%%%%%%%%%%%%%%%%%%%%%%%%%%%%%%%%%%%%%%%%%%%%%%%%%%%%%%%%%%%%%%%%%%%%%%%%%%%%%%%%%%%%%%%%%%%%%%%%%%%%%%%%%%%%%%%%%%%%%%%%%%%%%%%%%%%%%%%%%
%%%%%%%%%%%%%%%%%%%%%%%%%%%%%%%%%%%%%%%%%%%%%%%%%%%%%%%%%%%%%%%%%%%%%%%%%%%%%%%%%%%%%%%%%%%%%%%%%%%%%%%%%%%%%%%%%%%%%%%%%%%%%%%%%%%%%%%%%%%%%%%%%%%%%%%%%%%%%%%%%%%%%%%%%%%%%%%%%%%%%%%%%%%%%%%%%%%%%%%%%%%%%%%%%%%%%%%%%%%%%%%%%

Using \eqref{entropy0} we can easily derive the following crucial estimate for the first solution component. 
\begin{lemma}
\label{entropy_lemma_ablog}
Let $T > 0$. Then there exists a $C_Z(T) > 0$ such that
\begin{align}
\label{uln_bound}
\int_\Omega z(\cdot,t)\,|\ln z(\cdot,t)| \leq C_Z(T), \quad \text{for all}~ t \in (0, \widehat{T}),
  \end{align}
where $\widehat{T} := \min\{T, T_{\max}\}$.
\end{lemma}
\begin{proof}
    Thanks to Lemma \ref{entropy_lemma} there exists a $C_1(T) > 0$ such that
    \begin{align*}
        \int_\Omega c_1 \ln c_1 \leq C_1(T), \qquad \text{for all}~ t \in (0, \widehat{T}).
    \end{align*}
    In view of \eqref{tau_bound} and \eqref{cov1} we directly have
    \begin{align*}
    \int_\Omega z \ln z = \int_\Omega e^{-\chi \tau} c_1 \ln \left(e^{-\chi v} c_1\right) & =  \int_\Omega e^{-\chi \tau} c_1 (- \chi \tau + \ln c_1)\\[5pt] & \leq \int_\Omega c_1 \ln c_1 \leq C_1(T), \qquad \text{for all}~ t \in (0, \widehat{T}).
\end{align*}
Using once again $s\ln s\geq -e^{-1}$ for all $s>0$, we obtain
 \begin{align*}
    \nonumber \int_\Omega z(\cdot, t) |\ln z(\cdot, t)| & = \int_\Omega z(\cdot, t) \ln z(\cdot, t) - 2 \int_{\{z < 1\}} z(\cdot, t) \ln z(\cdot, t)                  \\[5pt]
     & \leq \int_\Omega z(\cdot, t) \ln z(\cdot, t) + 2 e^{-1} |\Omega| \leq C_1(T) +  2 e^{-1}|\Omega | =: C_Z(T). 
\end{align*}
for all $t\in(0,\widehat{T})$.
\end{proof}

%%%%%%%%%%%%%%%%%%%%%%%%%%%%%%%%%%%%%%%%%%%%%%%%%%%%%%%%%%%%%%%%%%%%%%%%%%%%%%%%%%%%%%%%%%%%%%%%%%%%%%%%%%%%%%%%%%%%%%%%%%%%%%%%%%%%%%%%%%%%%%%%%%%%%%%%%%%%%%%%%%%%%%%%%%%%%%%%%%%%%%%%%%%%%%%%%%%%%%%%%%%%%%%%%%%%%%%%%%%%%%%%%
%%%%%%%%%%%%%%%%%%%%%%%%%%%%%%%%%%%%%%%%%%%%%%%%%%%%%%%%%%%%%%%%%%%%%%%%%%%%%%%%%%%%%%%%%%%%%%%%%%%%%%%%%%%%%%%%%%%%%%%%%%%%%%%%%%%%%%%%%%%%%%%%%%%%%%%%%%%%%%%%%%%%%%%%%%%%%%%%%%%%%%%%%%%%%%%%%%%%%%%%%%%%%%%%%%%%%%%%%%%%%%%%%

\subsubsection{\texorpdfstring{$L^\infty$}{Linfty} bounds for \texorpdfstring{$c_1$}{c1}}
The primary goal here is to develop a Moser-Alikakos-type iteration argument yielding an $L^\infty(\Omega)$ bound for the first solution component $z$. The following lemma will play a crucial role in the implication of this scheme. 

\begin{lemma}
  \label{lloglineq1c}
  Let $p > 1$ and $T > 0$. Then for all $\varrho > 0$ there exists a $C(\varrho, p, T)$ such that
  \begin{align}
  \label{llogineq0c}
    \int_\Omega z^{p + 1} \leq  \frac{p^3\varepsilon C_Z(T)}{8 D_1} \int_\Omega  z^{p - 2} \nabla z \cdot \mathbb{D} \nabla z + C_1(p, \varepsilon) \cdot \left\{M_1^{p + 1} + 1\right\}, \quad \text{for all}~ t \in (0, \widehat{T}),
  \end{align}
  where $\widehat{T} := \min\{T, T_{\max}\}$.
\end{lemma}

\begin{proof}
  The boundedness of $z$ in $L \log L$ implied by Lemmas \ref{entropy_lemma_ablog} allows us to apply the refined interpolation Gagliardo-Nirenberg inequality type \cite[Lemma 2.1]{li2021boundedness}, which can be regarded as a special case of Lemma A.5 in \cite{tao2014energy} and \eqref{assumption0} to have
  \begin{align*}
    \nonumber \int_\Omega z^{p+1} &= \|z^{\frac p2}\|_{L^{\frac{2(p+1)}p}(\Omega)}^{\frac{2(p+1)}p} \\[5pt]
   \nonumber &\leq \varepsilon\|\nabla z^{\frac{p}{2}}\|^2_{L^2(\Omega)} \cdot \int_\Omega z |\ln z^{\frac{p}{2}}| + C(p, \varepsilon) \cdot \left(\|z^{\frac{p}{2}}\|^{\frac{2(p + 1)}{p}}_{L^{\frac{2}{p}}(\Omega)} + 1\right) \\[5pt]
   \nonumber & = \frac{p\varepsilon}{2} \|\nabla z^{\frac{p}{2}}\|^2_{L^2(\Omega)} \cdot \int_\Omega z |\ln z| + C(p, \varepsilon) \cdot \left\{\left(\int_\Omega a\right)^{p + 1} + 1\right\}\\[5pt]
   \nonumber & \leq \frac{p^3\varepsilon C_Z(T)}{8} \int_\Omega z^{p - 2} |\nabla z|^2 + C_1(p, \varepsilon) \cdot \left\{M_1^{p + 1} + 1\right\} \\[5pt]
    & \leq  \frac{p^3\varepsilon C_Z(T)}{8 D_1} \int_\Omega  z^{p - 2} \nabla z \cdot \mathbb{D} \nabla z + C_1(p, \varepsilon) \cdot \left\{M_1^{p + 1} + 1\right\}
  \end{align*}
  for all $t \in (0, \widehat{T})$, thus yielding the announced estimate.
\end{proof}

%%%%%%%%%%%%%%%%%%%%%%%%%%%%%%%%%%%%%%%%%%%%%%%%%%%%%%%%%%%%%%%%%%%%%%%%%%%%%%%%%%%%%%%%%%%%%%%%%%%%%%%%%%%%%%%%%%%%%%%%%%%%%%%%%%%%%%%%%%%%%%%%%%%%%%%%%%%%%%%%%%%%%%%%%%%%%%%%%%%%%%%%%%%%%%%%%%%%%%%%%%%%%%%%%%%%%%%%%%%%%%%%%
%%%%%%%%%%%%%%%%%%%%%%%%%%%%%%%%%%%%%%%%%%%%%%%%%%%%%%%%%%%%%%%%%%%%%%%%%%%%%%%%%%%%%%%%%%%%%%%%%%%%%%%%%%%%%%%%%%%%%%%%%%%%%%%%%%%%%%%%%%%%%%%%%%%%%%%%%%%%%%%%%%%%%%%%%%%%%%%%%%%%%%%%%%%%%%%%%%%%%%%%%%%%%%%%%%%%%%%%%%%%%%%%%

We will now apply a standard testing procedure to the first equation of \eqref{sys1}, leading to an $L^p$ estimate for the solution component $z$.

\begin{lemma}
  \label{Lplemma}
  Let $T > 0$ and $p \geq 2$. Then there exists a $C(p, T) > 0$ such that
  \begin{align}
   \label{z_Lpbound}
    \int_\Omega z^p(\cdot, t) \leq C(p, T) \quad \text{for all}~ t \in (0, \widehat{T}),
  \end{align}
  where $\widehat{T} := \min\{T, T_{\max}\}$.
\end{lemma}

\begin{proof}
  Using integration by parts, we obtain the following from the first equation of \eqref{sys1}
  \begin{align}
  \label{Lpz1}
    \nonumber & \frac{d}{dt} \int_\Omega e^{\chi \tau} z^{p}  =  p\int_\Omega e^{\chi \tau} z^{p - 1} z_t + \chi \int_\Omega e^{\chi \tau} z^p \tau_t \\[5pt]
    \nonumber& = p \int_\Omega e^{\chi \tau} z^{p - 1} \Big\{e^{- \chi \tau} \nabla \cdot (e^{\chi \tau} \mathbb{D} \nabla z) + e^{- \chi \tau} \nabla \cdot (z e^{\chi \tau} \nabla \cdot \mathbb{D}) - \alpha_1 \tfrac{z}{1 + z e^{\chi \tau}} + \alpha_2 \tfrac{c_2 e^{-\chi \tau}}{1 + c_2} \\[5pt]
    \nonumber & + \beta z (1 - z e^{\chi \tau} - c_2 -\tau) + \delta \chi e^{\chi \tau} z^2 \tau + \mu \chi \tau z - \tfrac{\chi z c_2}{1 + c_2}\Big\} + \chi \int_\Omega e^{\chi \tau} z^p \Big\{-\delta z  e^{\chi \tau} \tau - \mu \tau + \tfrac{c_2}{1 + c_2}\Big\} \\[5pt]
    \nonumber & = - p (p - 1) \int_\Omega e^{\chi \tau} z^{p - 2} (\nabla z \cdot \mathbb{D} \nabla z) - p(p - 1) \int_\Omega e^{\chi \tau} z^{p - 1}  ((\nabla \cdot \mathbb{D}) \cdot \nabla z) \\[5pt]
    \nonumber & +  \int_\Omega e^{\chi \tau} z^p \left\{(p-1) \mu \chi \tau  - \tfrac{p \alpha_1}{1 + z e^{\chi \tau}} - (p-1) \chi \tfrac{c_2}{1 + c_2} \right\} + \int_\Omega p \beta e^{\chi \tau} z^p (1 - z e^{\chi \tau} - c_2 -\tau) \\[5pt]
    & + p \alpha_2 \int_\Omega \tfrac{c_2}{1 + c_2} z^{p - 1} + (p - 1) \delta \chi \int_\Omega e^{2 \chi \tau} \tau z^{p + 1}
  \end{align}
  for all $t \in (0, T_{\max})$. Using \eqref{c2_bound} and \eqref{tau_bound}, and applying Young's inequality, we obtain from \eqref{Lpz1} that
  \begin{align}
    \label{Lpz2}
    \nonumber \frac{d}{dt} \int_\Omega e^{\chi \tau} z^{p} & \leq - p (p - 1) \int_\Omega e^{\chi \tau} z^{p - 2} (\nabla z \cdot \mathbb{D} \nabla z) - p(p - 1) \int_\Omega e^{\chi \tau} z^{p - 1}  ((\nabla \cdot \mathbb{D}) \cdot \nabla z) \\[5pt]
    & + C_1 \int_\Omega z^p + C_2(p) \int_\Omega z^p + (p - 1)\delta \chi e^{2 \chi M_\tau} M_\tau \int_\Omega z^{p + 1} + C_3(p) \alpha_2^p|\Omega|
  \end{align}
  for all $t \in (0, T_{\max})$ with an adequate $C_1 > 0$. Given the assumptions \eqref{assumption0} and \eqref{assumption1} on $\mathbb{D}$ we can use Young's inequality to estimate
 \begin{align}
   \label{LPz3}
   &  \nonumber  - p (p - 1) \int_\Omega e^{\chi \tau} z^{p - 2} (\nabla z \cdot \mathbb{D} \nabla z) - p(p - 1) \int_\Omega e^{\chi \tau} z^{p - 1}  ((\nabla \cdot \mathbb{D}) \cdot \nabla z)         \\[5pt]
   &  \nonumber \leq - p(p - 1) D_1 \int_\Omega e^{\chi \tau} z^{p - 2} |\nabla z|^2 + M p (p - 1) \int_\Omega e^{\chi \tau} z^{p - 1} |\nabla z|                                                       \\[5pt]
   &  \leq - \frac{D_1 p(p - 1)}{2} \int_\Omega e^{\chi \tau} z^{p - 2} |\nabla z|^2 + \frac{p(p - 1) M^2}{2 D_1} \int_\Omega e^{\chi \tau} z^p
 \end{align}
 for all $t \in (0, T_{\max})$. Inserting \eqref{LPz3} in \eqref{Lpz2} results in
 \begin{align}
    \label{Lpz4}
    \frac{d}{dt} \int_\Omega e^{\chi \tau} z^{p} + \frac{D_1 p(p - 1)}{2} \int_\Omega e^{\chi \tau} z^{p - 2} |\nabla z|^2 & \leq C_4(p) \int_\Omega z^p + (p - 1)\delta \chi e^{2 \chi M_\tau} M_\tau \int_\Omega z^{p + 1} + C_5
  \end{align}
for all $t \in (0, T_{\max})$. Since $e^{\chi \tau} \geq 1$, we have from \eqref{Lpz4}
\begin{align}
    \label{Lpz5}
    \frac{d}{dt} \int_\Omega e^{\chi \tau} z^{p} +  \zeta(p) \int_\Omega z^{p - 2} |\nabla z|^2 & \leq C_6(p) \int_\Omega e^{\chi \tau} z^p + (p - 1)\delta \chi e^{2 \chi M_\tau} M_\tau \int_\Omega z^{p + 1} + C_5
  \end{align}
for all $t \in (0, T_{\max})$ with a suitable $\zeta(p) > 0$. Let $w(t) := \int_\Omega e^{\chi \tau(\cdot, t)} z^p(\cdot, t), t \in (0, T_{\max})$. Then $w(t)$ satisfies
\begin{align}
  \label{Lpz6}
  w'(t) + \zeta(p)  \int_\Omega  z^{p - 2} |\nabla z|^2 \leq C_6(p) w(t) + C_7(p) \cdot \left(\int_\Omega z^{p + 1}\right) + C_5
\end{align}
for all $t \in (0, T_{\max})$. Then, as in \cite[(4.12)]{surulescu2021does}, we can also obtain the following inequality by integrating \eqref{Lpz6}
\begin{align}
  \label{Lpz7}
  w(t) + \zeta(p) \int_\sigma ^t \int_\Omega z^{p - 2} |\nabla z|^2 \leq w(\sigma) e^{C_6(p) T} + C_7(p) e^{C_6(p) T} \cdot \left(\int_\sigma ^t \int_\Omega z^{p + 1}\right) + C_8(p, T)
\end{align}
for all $t \in (0, T_{\max})$ {and some $\sigma \ge 0$}. For a fixed $T > 0$, thanks to \eqref{llogineq0c} we can choose $\varepsilon :=\varepsilon(p, T) > 0$ sufficiently small such that for some $C_{10}(p, T) > 0$ 
\begin{align}
  \label{Lpz8}
   C_7(p) e^{C_6(p) T} \int_0^t \int_\Omega z^{p + 1} \leq \zeta(p) \int_0^t \int_\Omega z^{p - 2} |\nabla z|^2 + C_{9}(p, \varepsilon, T) \quad \text{for all}~ t \in (0, \widehat{T}).
\end{align}
Inserting \eqref{Lpz8} into \eqref{Lpz7} results in
\begin{align*}
  w(t) \leq w(\sigma) e^{C_6(p) t} + C_{9}(p, T) =: C_{10}(p, T) \quad \text{for all}~ t \in (0, \widehat{T}),
\end{align*}
which ensures that \eqref{z_Lpbound} holds true.
\end{proof}

%%%%%%%%%%%%%%%%%%%%%%%%%%%%%%%%%%%%%%%%%%%%%%%%%%%%%%%%%%%%%%%%%%%%%%%%%%%%%%%%%%%%%%%%%%%%%%%%%%%%%%%%%%%%%%%%%%%%%%%%%%%%%%%%%%%%%%%%%%%%%%%%%%%%%%%%%%%%%%%%%%%%%%%%%%%%%%%%%%%%%%%%%%%%%%%%%%%%%%%%%%%%%%%%%%%%%%%%%%%%%%%%%
%%%%%%%%%%%%%%%%%%%%%%%%%%%%%%%%%%%%%%%%%%%%%%%%%%%%%%%%%%%%%%%%%%%%%%%%%%%%%%%%%%%%%%%%%%%%%%%%%%%%%%%%%%%%%%%%%%%%%%%%%%%%%%%%%%%%%%%%%%%%%%%%%%%%%%%%%%%%%%%%%%%%%%%%%%%%%%%%%%%%%%%%%%%%%%%%%%%%%%%%%%%%%%%%%%%%%%%%%%%%%%%%%
We are now in a position to apply a well-established Moser-Alikakos iteration technique to the present case, which allows us to convert the estimate \eqref{z_Lpbound} into an $L^\infty$ bound.

\begin{lemma}
  \label{lemma_Linfy_bound}
  Let $T > 0$, then there exists a $C(T) > 0$ such that
  \begin{align}
    \label{L_infy_bound}
    \|z(\cdot, t)\|_{L^\infty(\Omega)} \leq C(T), ~~ \text{and, therefore} ~~ \|c_1(\cdot, t)\|_{L^\infty(\Omega)} \leq C(T)  \quad \text{for all}~ t \in (0, \widehat{T}),
  \end{align}
  where $\widehat{T} := \min\{T, T_{\max}\}$.
\end{lemma}

\begin{proof}
 The first estimate in \eqref{L_infy_bound} can be obtained by following the procedure outlined in \cite[Lemma 4.5]{surulescu2021does}. The next part follows directly from this estimate, along with \eqref{tau_bound} and \eqref{cov1}.
\end{proof}

%%%%%%%%%%%%%%%%%%%%%%%%%%%%%%%%%%%%%%%%%%%%%%%%%%%%%%%%%%%%%%%%%%%%%%%%%%%%%%%%%%%%%%%%%%%%%%%%%%%%%%%%%%%%%%%%%%%%%%%%%%%%%%%%%%%%%%%%%%%%%%%%%%%%%%%%%%%%%%%%%%%%%%%%%%%%%%%%%%%%%%%%%%%%%%%%%%%%%%%%%%%%%%%%%%%%%%%%%%%%%%%%%%%%
%%%%%%%%%%%%%%%%%%%%%%%%%%%%%%%%%%%%%%%%%%%%%%%%%%%%%%%%%%%%%%%%%%%%%%%%%%%%%%%%%%%%%%%%%%%%%%%%%%%%%%%%%%%%%%%%%%%%%%%%%%%%%%%%%%%%%%%%%%%%%%%%%%%%%%%%%%%%%%%%%%%%%%%%%%%%%%%%%%%%%%%%%%%%%%%%%%%%%%%%%%%%%%%%%%%%%%%%%%%%%%%%%%%%

%%%%%%%%%%%%%%%%%%%%%%%%%%%%%%%%%%%%%%%%%%%%%%%%%%%%%%%%%%%%%%%%%%%%%%%%%%%%%%%%%%%%%%%%%%%%%%%%%%%%%%%%%%%%%%%%%%%%%%%%%%%%%%%%%%%%%%%%%%%%%%%%%%%%%%%%%%%%%%%%%%%%%%%%%%%%%%%%%%%%%%%%%%%%%%%%%%%%%%%%%%%%%%%%%%%%%%%%%%%%%%%%%
%%%%%%%%%%%%%%%%%%%%%%%%%%%%%%%%%%%%%%%%%%%%%%%%%%%%%%%%%%%%%%%%%%%%%%%%%%%%%%%%%%%%%%%%%%%%%%%%%%%%%%%%%%%%%%%%%%%%%%%%%%%%%%%%%%%%%%%%%%%%%%%%%%%%%%%%%%%%%%%%%%%%%%%%%%%%%%%%%%%%%%%%%%%%%%%%%%%%%%%%%%%%%%%%%%%%%%%%%%%%%%%%%

\subsubsection{Gradient estimates in \texorpdfstring{$L^4$}{L4}}

As the $L^\infty$ bounds for all the solution components have been established in Lemmas \ref{lemma_c2_bound}, \ref{lemma_tau_bound} and \ref{lemma_Linfy_bound}, we are now in a position to derive appropriate bounds for the haptotactic gradient component, which will ultimately lead to the global existence of classical solutions for \eqref{sys1}. Before we start the discussion, we would like to restate that $\nabla \cdot (\mathbb{D} \nabla \cdot)$ is an elliptic operator because of the assumed positive definiteness of $\mathbb{D}$ over $\overline{\Omega}$. This structural characteristic of the involved operator is crucial in the analysis that follows.

We begin our discussion by noting the following estimate \cite[(3.9)]{heihoff2023global}, which results from the Gagliardo-Nirenberg inequality \cite[Theorem 3.1.1]{lunardi2012analytic}, along with standard regularity estimates for the elliptic operator $\nabla \cdot (\mathbb{D} \nabla \cdot)$ under zero-flux boundary conditions: {there exists a constant $K^\ast>0$ such that}
\begin{align}
\label{gag_ninnn11}
\int_{\Omega} |\nabla \psi|^4
\leq K^\ast
\left(
\int_{\Omega} |\nabla \cdot (\mathbb{D}\nabla \psi)|^2
+
\int_{\Omega} |\psi|^2
\right)
\|\psi\|_{L^\infty(\Omega)}^2
\end{align}
for all $\psi \in C^2(\overline{\Omega})$ with
$(\mathbb{D}\nabla \psi)\cdot \nu = 0$ on $\partial\Omega$. Hence, in both two and three spatial dimensions, the $L^\infty(\Omega)$ estimate for $z$ provided by \eqref{L_infy_bound} allows us to infer from the above inequality that for all $T > 0$ 
\begin{align}
  \label{gag_niren_1}
  \int_\Omega |\nabla z|^4 \leq K^\ast \widehat{C}(T) \left(\int_\Omega |\nabla \cdot (\mathbb{D} \nabla z)|^2 + 1\right)
\end{align}
for all $t \in (0, \widehat{T})$ where $\widehat{T} := \min\{T, T_{\max}\}$. In the subsequent lemma we will test the following variant of the first equation in \eqref{sys1} with $- \nabla \cdot (\mathbb{D} \nabla z)$:
\begin{align}
  \label{gradnab1}
  \nonumber \partial_t z & = \nabla \cdot (\mathbb{D} \nabla z) + \chi \nabla \tau \cdot \mathbb{D} \nabla z + \nabla \cdot (z (\nabla \cdot \mathbb{D})) + \chi z(\nabla \tau \cdot (\nabla \cdot \mathbb{D})) - \alpha_1 \tfrac{z}{1 + z e^{\chi \tau}} + \alpha_2 \tfrac{c_2 e^{-\chi \tau}}{1 + c_2} \\[5pt]
  & + \beta z (1 - z e^{\chi \tau} - c_2 -\tau)+ \delta \chi e^{\chi \tau} z^2 \tau + \mu \chi \tau z - \tfrac{\chi z c_2}{1 + c_2}.
\end{align}

%%%%%%%%%%%%%%%%%%%%%%%%%%%%%%%%%%%%%%%%%%%%%%%%%%%%%%%%%%%%%%%%%%%%%%%%%%%%%%%%%%%%%%%%%%%%%%%%%%%%%%%%%%%%%%%%%%%%%%%%%%%%%%%%%%%%%%%%%%%%%%%%%%%%%%%%%%%%%%%%%%%%%%%%%%%%%%%%%%%%%%%%%%%%%%%%%%%%%%%%%%%%%%%%%%%%%%%%%%%%%%%%%
%%%%%%%%%%%%%%%%%%%%%%%%%%%%%%%%%%%%%%%%%%%%%%%%%%%%%%%%%%%%%%%%%%%%%%%%%%%%%%%%%%%%%%%%%%%%%%%%%%%%%%%%%%%%%%%%%%%%%%%%%%%%%%%%%%%%%%%%%%%%%%%%%%%%%%%%%%%%%%%%%%%%%%%%%%%%%%%%%%%%%%%%%%%%%%%%%%%%%%%%%%%%%%%%%%%%%%%%%%%%%%%%%
\begin{lemma}
   \label{lemma_nabL4_bound}
   Let $T > 0$. Then there exists a $C(T) > 0$ such that
    \begin{align}
      \label{tau_c2ess}
      \|\nabla c_2(\cdot, t)\|_{L^4(\Omega)} \leq C(T) \quad \text{and} \quad \|\nabla \tau(\cdot, t)\|_{L^4(\Omega)} \leq C(T) \quad  \text{for all}~ t \in (0, \widehat{T}),
  \end{align}
  where $\widehat{T} := \min\{T, T_{\max}\}$.
\end{lemma}
\begin{proof}
  As mentioned above we test \eqref{gradnab1} with $- \nabla \cdot (\mathbb{D} \nabla z)$ and  employ integration by parts (noting $(\nabla \cdot \mathbb{D}) \cdot \nu = 0$ on $\partial \Omega$) and Young's inequality to have
  \begin{align}
  \label{nab_tau1}
    \nonumber & \frac{1}{2} \frac{d}{dt} \int_\Omega (\nabla z \cdot \mathbb{D} \nabla z)  = \int_\Omega \nabla z_t \cdot \mathbb{D} \nabla z   \\[5pt]
    \nonumber & = \int_\Omega \nabla (\nabla \cdot (\mathbb{D} \nabla z)) \cdot \mathbb{D} \nabla z + \chi \int_\Omega \nabla (\nabla \tau \cdot \mathbb{D} \nabla z) \cdot \mathbb{D} \nabla z + \chi \int_\Omega \nabla (z\nabla \tau \cdot (\nabla \cdot \mathbb{D})) \cdot \mathbb{D} \nabla z \\[5pt]
    \nonumber & + \int_\Omega \nabla(\nabla \cdot (z \nabla \cdot \mathbb{D})) \cdot \mathbb{D} \nabla z  + \int_\Omega \nabla \Big\{- \alpha_1 \tfrac{z}{1 + z e^{\chi \tau}} + \alpha_2 \tfrac{c_2 e^{-\chi \tau}}{1 + c_2} + \beta z (1 - z e^{\chi \tau} - c_2 -\tau) \\[5pt]
    \nonumber & \hspace{9cm} + \delta \chi e^{\chi \tau} z^2 \tau + \mu \chi \tau z - \tfrac{\chi z c_2}{1 + c_2}\Big\} \cdot \mathbb{D} \nabla z \\[5pt]
    \nonumber & \leq -  \frac{1}{2} \int_\Omega |\nabla \cdot (\mathbb{D} \nabla z)|^2 + 2 \chi^2 \int_\Omega |\nabla \tau \cdot \mathbb{D} \nabla \tau||\nabla z \cdot \mathbb{D} \nabla z| + 2 \int_\Omega |\nabla \cdot (z \nabla \cdot \mathbb{D})|^2 \\[5pt]
    \nonumber & + 2 \chi^2 \int_\Omega z^2 |\nabla \tau|^2 |\nabla \cdot \mathbb{D}|^2 + 2\int_\Omega  \Big\{- \alpha_1\tfrac{z}{1 + z e^{\chi \tau}} + \alpha_2 \tfrac{c_2 e^{-\chi \tau}}{1 + c_2} + \beta z (1 - z e^{\chi \tau} - c_2 -\tau)  \\
    &\hspace{9cm} + \delta \chi e^{\chi \tau} z^2 \tau + \mu \chi \tau z - \tfrac{\chi z c_2}{1 + c_2}\Big\}^2
  \end{align}
  for all $t > 0$. From \eqref{c2_bound}, \eqref{tau_bound},  and the first estimate in \eqref{L_infy_bound} we can have a $C_1(T) > 0$ such that
  \begin{align}
   \label{nab_tau2}
    \int_\Omega  \Big\{- \tfrac{\alpha_1 z}{1 + z e^{\chi \tau}} + \tfrac{\alpha_2 c_2 e^{-\chi \tau}}{1 + c_2} + \beta z (1 - z e^{\chi \tau} - c_2 -\tau)  + \delta \chi e^{\chi \tau} z^2 \tau + \mu \chi \tau z - \tfrac{\chi z c_2}{1 + c_2}\Big\}^2 \leq C_1(T)
  \end{align}
  for all $t \in (0, \widehat{T})$. Using the assumptions on $\mathbb{D}$ outlined in \eqref{assumption1}, applying Young's inequality and taking into account the first estimate in \eqref{L_infy_bound}, we obtain the following estimates for some constants $K_i > 0, i = 1 \cdots 6$ and for all $\rho > 0$:
  \begin{align}
    \label{nab_tau3}
     & 2 \chi^2 \int_\Omega |\nabla \tau \cdot \mathbb{D} \nabla \tau||\nabla z \cdot \mathbb{D} \nabla z| \leq \frac{K_1}{\rho} \int_\Omega |\nabla \tau|^4 + K_2 \rho \int_\Omega |\nabla z|^4,              \\[5pt]
    \label{nab_tau4}
     & 2 \int_\Omega |\nabla \cdot (z \nabla \cdot \mathbb{D})|^2 \leq 4 M^2 \left(\int_\Omega |\nabla z|^2 + \int_\Omega z^2\right) \leq K_3\rho \int_\Omega |\nabla z|^4 + \frac{K_4(T)}{\rho},                \\[5pt]
    \label{nab_tau5}
    &  2 \chi^2 \int_\Omega z^2 |\nabla \tau|^2 |\nabla \cdot \mathbb{D}|^2 \leq \chi^2 M^2 \left(\int_\Omega z^4 + \int_\Omega |\nabla \tau|^4\right) \leq K_5 \int_\Omega |\nabla \tau|^4 + K_6(T)
  \end{align}
 for all $t \in (0, \widehat{T})$. Inserting \eqref{nab_tau2}-\eqref{nab_tau5} in \eqref{nab_tau1} results in 
  \begin{align}
  \label{nab_tau6}
   \nonumber \frac{1}{2} \frac{d}{dt} \int_\Omega (\nabla z \cdot \mathbb{D} \nabla z) \leq -  \frac{1}{2} \int_\Omega |\nabla \cdot (\mathbb{D} \nabla z)|^2 & + (K_2 + K_3)\rho \int_\Omega |\nabla z|^4  \\[5pt]
   & + C_2(\rho) \int_\Omega |\nabla \tau|^4 + C_3(\rho, T)
 \end{align}
 for all $t \in (0, \widehat{T})$. Our next goal is to establish an estimate for the time derivative of the gradient terms of the haptotactic component $\tau$. To this end we calculate
 \begin{align}
   \label{nab_tau7}
   \nonumber & \frac{1}{4} \frac{d}{dt} \int_\Omega |\nabla \tau|^4 = \int_\Omega |\nabla \tau|^2 \nabla \tau \cdot \nabla \tau_t = \int_\Omega |\nabla \tau|^2 \nabla \tau \cdot \nabla \left\{-\delta z  e^{\chi \tau} \tau - \mu \tau + \frac{c_2}{1 + c_2}\right\}       \\[5pt]
   & = \int_\Omega \left\{- \delta e^{\chi \tau} z - \delta \chi e^{\chi \tau} \tau z - \mu\right\} |\nabla \tau|^4 - \delta \int_\Omega e^{\chi \tau} \tau |\nabla \tau|^2 \nabla \tau \cdot \nabla z + \int_\Omega \frac{|\nabla \tau|^2}{(1 + c_2)^2} \nabla \tau \cdot \nabla c_2
 \end{align}
 for all $t > 0$. From  \eqref{c2_bound}, \eqref{tau_bound} and the first estimate in \eqref{L_infy_bound} we can have positive constants $C_4(T), C_5(T)$ and $C_6$ such that \eqref{nab_tau7} can be expressed in the following form:
 \begin{align}
  \label{nab_tau8}
  \frac{1}{4} \frac{d}{dt} \int_\Omega |\nabla \tau|^4 \leq C_4(T) \int_\Omega |\nabla \tau|^4 + C_5(T) \int_\Omega |\nabla \tau|^2 \nabla \tau \cdot \nabla z + C_6 \int_\Omega |\nabla \tau|^2 \nabla \tau \cdot \nabla c_2
 \end{align}
 for all $t \in (0, \widehat{T})$. Applying Young's inequality to the last two terms of \eqref{nab_tau8} yields the following for some positive constants and for every $\rho > 0$:
\begin{align}
  \label{nab_tau9}
  \frac{d}{dt} \int_\Omega |\nabla \tau|^4 \leq C_7(T, \rho) \int_\Omega |\nabla \tau|^4 + C_8(T) \rho\int_\Omega |\nabla z|^4 + C_{9} \int_\Omega |\nabla c_2|^4
 \end{align}
 for all $t \in (0, \widehat{T})$. Our next task is to handle the term $\int_\Omega |\nabla c_2|^4$ on the right-hand side of \eqref{nab_tau9}. To this end, we derive an estimate for its time derivative: %of the gradient of the second solution component $c_2$:
 \begin{align}
   \label{c2nab_tau1}
   \nonumber & \frac{1}{4} \frac{d}{dt} \int_\Omega |\nabla c_2|^4 = \int_\Omega |\nabla c_2|^2 \nabla c_2 \cdot \nabla \partial_tc_{2} = \int_\Omega |\nabla c_2|^2 \nabla c_2 \cdot \nabla \left\{\alpha_1 \frac{z e^{\chi \tau}}{1 + z e^{\chi \tau}} - \alpha_2 \frac{c_2}{1 + c_2}\right\}       \\[5pt]
   & = \alpha_1 \int_\Omega \frac{e^{\chi \tau}|\nabla c_2|^2}{(1 + z e^{\chi \tau})^2} \nabla c_2 \cdot \nabla z + \alpha_1 \chi \int_\Omega \frac{ze^{\chi \tau} |\nabla c_2|^2}{(1 + ze^{\chi \tau})^2} \nabla c_2 \cdot \nabla \tau - \alpha_2 \int_\Omega \frac{|\nabla c_2|^4}{(1 + c_2)^2}
 \end{align}
 for all $t > 0$. Again, following the procedure outlined above in the derivation of \eqref{nab_tau9}, we obtain from \eqref{c2nab_tau1}, for some positive constants and for every $\rho > 0$:
 \begin{align}
  \label{c2nab_tau2}
  \frac{d}{dt} \int_\Omega |\nabla c_2|^4 \leq C_{10}(\rho) \int_\Omega |\nabla \tau|^4 + C_{11}(T) \rho \int_\Omega |\nabla z|^4 + C_{12}(T, \rho) \int_\Omega |\nabla c_2|^4
 \end{align}
 for all $t \in (0, \widehat{T})$. Combining \eqref{nab_tau6}, \eqref{nab_tau9} and \eqref{c2nab_tau2} we have 
 \begin{align}
   \label{main_tau1}
   \nonumber & \frac{d}{dt} \Big\{ \int_\Omega (\nabla z \cdot \mathbb{D} \nabla z) + \int_\Omega |\nabla c_2|^4 + \int_\Omega |\nabla \tau|^4\Big\} \\[5pt]
   \nonumber & \leq - \frac{1}{2}\int_\Omega |\nabla \cdot (\mathbb{D} \nabla z)|^2  + C_{13}(T)\rho \int_\Omega |\nabla z|^4 + C_{14}(T, \rho) \int_\Omega |\nabla c_2|^4 \\[5pt]
   & + C_{15}(T, \rho) \int_\Omega |\nabla \tau|^4 + C_{16}(T, \rho)
 \end{align}
 for all $t \in (0, \widehat{T})$. Setting $\rho := \frac{1}{2 C_{13} K^\ast \widehat{C}(T)}$ we can utilize \eqref{gag_niren_1} to have for some positive constant $C_{17} > 0$
 \begin{equation}
   \label{main_tau2}
   C_{13}(T) \rho \int_\Omega |\nabla z|^4 \leq \frac{1}{2}\int_\Omega |\nabla \cdot (\mathbb{D} \nabla z)|^2 + C_{17}
 \end{equation}
 for all $t \in (0, \widehat{T})$. Inserting \eqref{main_tau2} in \eqref{main_tau1} and utilizing the positivity of $\int_\Omega (\nabla z \cdot \mathbb{D} \nabla z)$ we get
 \begin{align}
   \label{main_tau3}
   \nonumber  \frac{d}{dt} \Big\{ \int_\Omega (\nabla z \cdot \mathbb{D} \nabla z) + \int_\Omega |\nabla c_2|^4 & + \int_\Omega |\nabla \tau|^4\Big\} \\[5pt]
   & \leq \kappa \left\{\int_\Omega (\nabla z \cdot \mathbb{D} \nabla z) + \int_\Omega |\nabla c_2|^4 + \int_\Omega |\nabla \tau|^4\right\} + C_{18}
 \end{align}
  for all $t \in (0, \widehat{T})$, where $\kappa > 0$ is a constant. An application of Gronwall's lemma to \eqref{main_tau3}  yields the desired estimates in \eqref{tau_c2ess}.
\end{proof}

%%%%%%%%%%%%%%%%%%%%%%%%%%%%%%%%%%%%%%%%%%%%%%%%%%%%%%%%%%%%%%%%%%%%%%%%%%%%%%%%%%%%%%%%%%%%%%%%%%%%%%%%%%%%%%%%%%%%%%%%%%%%%%%%%%%%%%%%%%%%%%%%%%%%%%%%%%%%%%%%%%%%%%%%%%%%%%%%%%%%%%%%%%%%%%%%%%%%%%%%%%%%%%%%%%%%%%%%%%%%%%%%%
%%%%%%%%%%%%%%%%%%%%%%%%%%%%%%%%%%%%%%%%%%%%%%%%%%%%%%%%%%%%%%%%%%%%%%%%%%%%%%%%%%%%%%%%%%%%%%%%%%%%%%%%%%%%%%%%%%%%%%%%%%%%%%%%%%%%%%%%%%%%%%%%%%%%%%%%%%%%%%%%%%%%%%%%%%%%%%%%%%%%%%%%%%%%%%%%%%%%%%%%%%%%%%%%%%%%%%%%%%%%%%%%%

\subsubsection{Global extensibility in two dimensions}

We are now in a position to prove the first main result of this section, which asserts global classical solvability for \eqref{model1} in two dimensions.
\begin{proof}[Proof of Theorem \ref{classical_main_th1}]
  Thanks to Lemma \ref{local_existence}, particularly the extensibility criterion \eqref{extensibility_criterion} therein, and combining the outcomes of Lemmas \ref{lemma_Linfy_bound} and \ref{lemma_nabL4_bound}, yields the desired result.
\end{proof}

%%%%%%%%%%%%%%%%%%%%%%%%%%%%%%%%%%%%%%%%%%%%%%%%%%%%%%%%%%%%%%%%%%%%%%%%%%%%%%%%%%%%%%%%%%%%%%%%%%%%%%%%%%%%%%%%%%%%%%%%%%%%%%%%%%%%%%%%%%%%%%%%%%%%%%%%%%%%%%%%%%%%%%%%%%%%%%%%%%%%%%%%%%%%%%%%%%%%%%%%%%%%%%%%%%%%%%%%%%%%%%%%%%%%
%%%%%%%%%%%%%%%%%%%%%%%%%%%%%%%%%%%%%%%%%%%%%%%%%%%%%%%%%%%%%%%%%%%%%%%%%%%%%%%%%%%%%%%%%%%%%%%%%%%%%%%%%%%%%%%%%%%%%%%%%%%%%%%%%%%%%%%%%%%%%%%%%%%%%%%%%%%%%%%%%%%%%%%%%%%%%%%%%%%%%%%%%%%%%%%%%%%%%%%%%%%%%%%%%%%%%%%%%%%%%%%%%%%%
\subsection{\texorpdfstring
{Global classical solvability to \eqref{model1} in three dimensions}
{Classical solutions3d}}

We now turn to the global classical solvability of \eqref{model1} in three spatial dimensions. Our approach is based on the coupled $L^p$-estimate technique developed in \cite{pang2018global,dai2022new} for haptotaxis problems without diffusion tensors. More precisely, we successively improve the regularity of $z$ along the scale
$$L^1(\Omega)\rightarrow L^2(\Omega)\rightarrow L^p(\Omega)\rightarrow L^\infty(\Omega),$$
for arbitrary $p>1$. To initiate this bootstrap procedure, we first show that the weighted energy functional
$u:=\int_\Omega e^{\chi\tau}z^2$ satisfies the autonomous differential inequality
$$u'(t)+u(t)\leq\Lambda^\ast u^3(t)+\widetilde{\Lambda} \qquad\text{for all }t\in(0,T_{\max}),$$
where $\Lambda^\ast,\widetilde{\Lambda}>0$ are appropriate constants. An application of the comparison principle for ordinary differential equations then guarantees the boundedness of $u(t)$, provided that the initial datum $c_{10}$ satisfies an appropriate smallness condition. This results in an $L^2$ estimate for $z$.

Having established an a priori $L^2(\Omega)$ estimate for $z$, we subsequently derive $L^p(\Omega)$ estimates for every $p>1$. The desired $L^\infty(\Omega)$ estimate is then obtained via the Moser-Alikakos iteration. Following this approach, we begin by establishing the fundamental a priori $L^2(\Omega)$ estimate for $z$.

\begin{lemma}
    \label{3d_lemma1}
     Let the assumptions of Theorem~\ref{3D_Class_Theorem} be satisfied. Then there exists a constant $C > 0$ such that
    \begin{align}
    \label{uwsql2}
    \|z(\cdot,t)\|_{L^2(\Omega)}
    \le C
    \qquad \text{for all } t\in(0,T_{\max}).
\end{align}
    \end{lemma}
\begin{proof}
  From the first equation in \eqref{sys1} we can directly calculate
  \begin{align}
      \label{3s22}
      \nonumber \frac{d}{dt} \int_\Omega e^{\chi \tau} z^2 & = 2 \int_\Omega e^{\chi \tau} z z_t + \chi \int_\Omega e^{\chi \tau}z^2 \tau_t  \\[5pt]
      \nonumber & = - 2 \int_\Omega e^{\chi \tau} (\nabla z \cdot \mathbb{D} \nabla z) - 2 \int_\Omega e^{\chi \tau} z ((\nabla \cdot \mathbb{D}) \cdot \nabla z) \\[5pt] 
      \nonumber & + \int_\Omega e^{\chi \tau} z^2 \left\{ \mu \chi \tau - \tfrac{c_2}{1 + c_2} \chi - \tfrac{2 \alpha_1}{1 + z e^{\chi \tau}} \right\} \\[5pt] 
      & + 2 \int_\Omega \beta e^{\chi \tau} z^2 (1 - z e^{\chi \tau} - c_2 -\tau) + 2 \alpha_2 \int_\Omega \frac{c_2}{1 + c_2} z + \delta \chi \int_\Omega e^{2 \chi \tau} z^3
  \end{align}
 for all $t > 0$. Recalling the nonnegativity of the solution components, the elementary bound $\frac{s}{1 + s}\leq 1$ for all $s\geq 0$, and \eqref{tau_bound}, we arrive at the following differential inequality:
 \begin{align}
     \label{3s45}
     \nonumber \frac{d}{dt} \int_\Omega e^{\chi \tau} z^2 & \leq - 2 \int_\Omega e^{\chi \tau} (\nabla z \cdot \mathbb{D} \nabla z) - 2 \int_\Omega e^{\chi \tau} z ((\nabla \cdot \mathbb{D}) \cdot \nabla z)  \\[5pt]
 & + \delta \chi e^{2 \chi M_\tau} \int_\Omega z^3 + 2 e^{\chi M_\tau} \left(\mu \chi M_\tau + \beta\right) \int_\Omega z^2 + 2 \alpha_2 \int_\Omega z
 \end{align}
for all $t > 0$. We first estimate the first two terms on the right-hand side of \eqref{3s45}. In view of \eqref{assumption0}-\eqref{assumption1} and Young's inequality, we obtain
\begin{align}
    \label{ewqq1}
    \nonumber & - 2 \int_\Omega e^{\chi \tau} (\nabla z \cdot \mathbb{D} \nabla z) - 2 \int_\Omega e^{\chi \tau} z ((\nabla \cdot \mathbb{D}) \cdot \nabla z) \\[5pt]
    \nonumber & \leq - 2 D_1 \int_\Omega e^{\chi \tau} |\nabla z|^2 + 2 M \int_\Omega e^{\chi \tau} z  |\nabla z| \\[5pt]
    & \leq - D_1 \int_\Omega e^{\chi \tau} |\nabla z|^2 + \frac{M^2}{D_1} \int_\Omega e^{\chi \tau} z^2
\end{align}
for all $t > 0$. Inserting \eqref{ewqq1} into \eqref{3s45} results in
\begin{align}
    \label{pl78}
    \nonumber \frac{d}{dt} \int_\Omega e^{\chi \tau} z^2 + \int_\Omega e^{\chi \tau} z^2  & +  D_1 \int_\Omega e^{\chi \tau} |\nabla z|^2 \leq \; \; \delta \chi e^{2 \chi M_\tau} \int_\Omega z^3 \\[5pt]
    & { + \left [2 e^{\chi M_\tau} \left(\mu \chi M_\tau + \beta\right ) +  e^{\chi M_\tau} \left (\tfrac{M^2}{D_1}+1 \right) \right ]}\int_\Omega z^2 + 2 \alpha_2 \int_\Omega z
\end{align}
for all $t > 0$. By means of Young's inequality, we have
\begin{align}
\label{eqaa}
    & \widehat{\Lambda}  \int_\Omega z^2 = \int_\Omega z^{\frac{3}{2}} \times \widehat{\Lambda}  z^{\frac{1}{2}}  \leq \frac{1}{2} \int_\Omega z^3 + \frac{\widehat{\Lambda}^2}{2} \int_\Omega z, 
\end{align}
for all $t > 0$ where 
\begin{align}
    \label{widehatlambda}
    \widehat{\Lambda} := { 2 e^{\chi M_\tau} \left(\mu \chi M_\tau + \beta\right ) +  e^{\chi M_\tau} \left (\tfrac{M^2}{D_1}+1 \right). }
    %2 e^{\chi M_\tau} \left(\mu \chi M_\tau + \beta + \tfrac{M^2}{D_1} + e^{\chi M_\tau}\right).
\end{align}
Inserting \eqref{eqaa} into \eqref{pl78} results in
\begin{align}
    \label{po78a}
    \frac{d}{dt} \int_\Omega e^{\chi \tau} z^2 + \int_\Omega e^{\chi \tau} z^2  + D_1 \int_\Omega e^{\chi \tau} |\nabla z|^2 \leq \left(\delta \chi e^{2 \chi M_\tau}  + \frac{1}{2}\right) \int_\Omega z^3 + \left(2 \alpha_2 + \tfrac{\widehat{\Lambda}^2}{2}\right)\int_\Omega z
\end{align}
for all $t > 0$. Invoking the Gagliardo-Nirenberg inequality  \cite[(4.8)-(4.10)]{dai2022new} with some constant $C_{\mathrm{GN}}^{(3)} = C_{\mathrm{GN}}^{(3)}(\Omega) > 0$ and Young's inequality, we have
\begin{align}
    \label{gag11}
    \nonumber \Lambda_\ast \|z\|_{L^3(\Omega)}^3 & \leq \Lambda_\ast C_{\mathrm{GN}}^{(3)} \|\nabla z\|^{\frac{3}{2}}_{L^2(\Omega)} \|z\|^{\frac{3}{2}}_{L^2(\Omega)} + \Lambda_\ast C_{\mathrm{GN}}^{(3)} \|z\|^3_{L^1(\Omega)} \\[5pt]
    & \leq D_1\|\nabla z\|^2_{L^2(\Omega)} + \frac{27}{2048 D_1} \Lambda_\ast^4 (C_{\mathrm{GN}}^{(3)})^4 \|z\|^6_{L^2(\Omega)} + \Lambda_\ast C_{\mathrm{GN}}^{(3)} \|z\|^3_{L^1(\Omega)}, 
\end{align}
for all $t > 0$, where 
\begin{align}
    \label{Lambda_ast}
    \Lambda_\ast := \left(\delta \chi e^{2 \chi M_\tau} + \frac{1}{2}\right).
\end{align}
Inserting \eqref{gag11} into \eqref{po78a} results in 
\begin{align}
    \label{plzt}
     \frac{d}{dt} \int_\Omega e^{\chi \tau} z^2 + \int_\Omega e^{\chi \tau} z^2 \leq \frac{27}{2048 D_1} \Lambda_\ast^4 (C_{\mathrm{GN}}^{(3)})^4 \left(\int_\Omega e^{\chi \tau} z^2 \right)^3 +  \Lambda_\ast C_{\mathrm{GN}}^{(3)} \|z\|^3_{L^1(\Omega)} + \overline{\Lambda} \|z\|_{L^1(\Omega)}
\end{align}
for all $t > 0$, where 
\begin{align}
    \label{ov_lambda}
    \overline{\Lambda} := \left(2 \alpha_2 + \tfrac{\widehat{\Lambda}^2}{2}\right),
\end{align}
where $\widehat{\Lambda}$ is defined in \eqref{widehatlambda}.
%\begin{align*}
   %{\cb z^\ast (t):=\|z(t)\|_ {L^1(\Omega)},}\quad   u(t):=\int_\Omega e^{\chi \tau(\cdot, t)} z^2(\cdot, t), \qquad \text{for all}~ t > 0.   
%\end{align*}
With the above notations {(also recall \eqref{zL1_bound})}, \eqref{plzt} takes the form
\begin{align}
\label{p0jh}
    u'(t) + u(t) \leq \Lambda^\ast u^3(t) +  \Lambda_\ast C_{\mathrm{GN}}^{(3)} (z^\ast)^3 + \overline{\Lambda} z^\ast
\end{align}
for all $t \in (0, T_{\max})$ with 
\begin{align}
\label{Lambdastar}
    \Lambda^\ast := \frac{27}{2048 D_1} \Lambda_\ast^4 (C_{\mathrm{GN}}^{(3)})^4.
\end{align}
Thanks to \eqref{p0jh} we have
\begin{align*}
    u'(t)  \leq \Psi(u(t)), \qquad \text{for all}~ t \in (0, T_{\max})
\end{align*}
with 
\begin{align*}
    \Psi(s) = - s + \Lambda^\ast s^3 +  \Lambda_\ast C_{\mathrm{GN}}^{(3)} (z^\ast)^3 + \overline{\Lambda} z^\ast \qquad \text{for all}~ s > 0.
\end{align*}
$\Psi(s)$ attains its minimum
\begin{align}
\label{couzt}
     \Psi(s_0) = - \frac{2 s_0}{3} +  \Lambda_\ast C_{\mathrm{GN}}^{(3)} (z^\ast)^3 + \overline{\Lambda} z^\ast 
\end{align}
at 
\begin{align*}
    s_0 := \sqrt{\frac{1}{3 \Lambda^\ast}}.
\end{align*}
We next establish the estimate $u(t)< s_0$ under a suitable smallness assumption on the initial data. Let $\Pi^\ast>0$ be a constant to be specified below, and suppose that there exists $\Pi>0$ satisfying
\begin{equation*}
\max\left\{
\|c_{10}\|_{L^1(\Omega)},
\frac{|\Omega|}{2}
\left(
1+\sqrt{1+\frac{4\alpha_2}{\beta}}
\right)
\right\}
\leq \Pi<\Pi^\ast.
\end{equation*}
Then, thanks to Lemma~\ref{lemmal1c1} we have 
\begin{equation}
\label{c1_lp}
\|c_1(\cdot,t)\|_{L^1(\Omega)}
\leq \Pi
\qquad\text{for all }t\in(0,T_{\max}).
\end{equation}
Thanks to Corollary~\ref{L1_z:cor} and \eqref{c1_lp}, we also have
\begin{align*}
\|z(\cdot,t)\|_{L^1(\Omega)} \leq \Pi
\qquad \text{for all } t\in(0,T_{\max}).
\end{align*}
We are therefore in a position to set $z^\ast < \Pi$. Consequently, we deduce from \eqref{couzt} that
\begin{align}
\label{couzt1}
     \Psi(s_0) = - \frac{2 s_0}{3} +  \Lambda_\ast C_{\mathrm{GN}}^{(3)} \Pi^3 + \overline{\Lambda} \Pi. 
\end{align}
Consider the function $F(r):=\Lambda_\ast C_{\mathrm{GN}}^{(3)}r^3+\overline{\Lambda}r, \; r \geq 0$.
Clearly, $F\in C([0,\infty))$. Moreover, $F'(r)=3\Lambda_\ast C_{\mathrm{GN}}^{(3)}r^2+\overline{\Lambda}>0$ for all $r\geq 0$, and hence $F$ is strictly increasing on $[0,\infty)$. Since $F(0)=0$ and $F(r)\to\infty$ as $r\to\infty$, the intermediate value theorem guarantees the existence of a unique constant $\Pi^\ast>0$ such that
\begin{equation*}
\Lambda_\ast C_{\mathrm{GN}}^{(3)}(\Pi^\ast)^3 +\overline{\Lambda}\Pi^\ast = \frac{2s_0}{3}.
\end{equation*}
Therefore, whenever $\Pi<\Pi^\ast$, the monotonicity of $F$ implies
\begin{align*}
\Lambda_\ast C_{\mathrm{GN}}^{(3)}\Pi^3 +\overline{\Lambda}\Pi = F(\Pi) < F(\Pi^\ast) = \frac{2s_0}{3}.
\end{align*}
Combining this with \eqref{couzt1}, we obtain $\Psi(s_0) < 0$. Since $\Psi(s_0) < 0$, a comparison argument applied to $u'(t)\leq \Psi(u(t))$ shows that $u(t) < s_0$ for all $t\in(0,T_{\max})$, provided that $u(0)\leq s_0$. The latter follows from the assumption $\|e^{- \chi \tau_0/2}c_{10}\|_{L^2(\Omega)}^2\leq s_0$. Therefore, in view of the above arguments, we have
\begin{align}
\label{l21}
\int_\Omega e^{\chi\tau}z^2 < s_0
\qquad \text{for all } t\in(0,T_{\max}).
\end{align}
Finally, the desired estimate \eqref{uwsql2} follows immediately from \eqref{l21}.
\end{proof}
%%%%%%%%%%%%%%%%%%%%%%%%%%%%%%%%%%%%%%%%%%%%%%%%%%%%%%%%%%%%%%%%%%%%%%%%%%%%%%%%%%%%%%%%%%%%%%%%%%%%%%%%%%%%%%%%%%%%%%%%%%%%%%%%%%%%%%%%%%%%%%%%%%%%%%%%%%%%%%%%%%%%%%%%%%%%%%%%%%%%%%%%%%%%%%%%%%%%%%%%%%%%%%%%%%%%%%%%%%%%%%%%%%%%%%%%%%%%%%%%%%%%%%%%%%%%%%%%%%%%%%%%%%%%%%%%%%%%%%%%%%%%%%%%%%%%%%%%%%%%%%%%%%%%%%%%%%%%%%%%%%%%%%%%%%%%%%%%%%%%%%%%%%%%%%%%%%%%%%%%%%%%%%%%%%%%%%%%%%%%%%%%%%%%%%%%%%%%%%%%%%%%%%%%%%%%%%%%%%%%%%%%%%%%%%%%%%%%%%%%%%%%%%%%%%%%%%%%%%%%%%%%%%%%%%%%%%%%%%%%%%%%%%%

With the estimate established in Lemma~\ref{3d_lemma1} at hand, we proceed to derive $L^p$-bounds $(p>1)$ for $z$ by employing a standard testing argument together with suitable Gagliardo-Nirenberg inequalities.

\begin{lemma}
    \label{3d_lemma2}
     Let the assumptions of Theorem~\ref{3D_Class_Theorem} be satisfied. Then for all $p > 1$ there exists a constant $C(T)>0$ such that
    \begin{align}
    \label{ulp}
    \|z(\cdot,t)\|_{L^p(\Omega)}
    \le C(T)
    \qquad \text{for all } t\in(0, T_{\max}).
\end{align}
    \end{lemma}
\begin{proof}
 Following the argument leading to \eqref{Lpz4} in Lemma \ref{Lplemma}, we obtain
 \begin{align}
    \label{Lpz43d}
    \nonumber \frac{d}{dt} \int_\Omega e^{\chi \tau} z^{p} + \int_\Omega e^{\chi \tau} z^{p} & + \frac{D_1 p(p - 1)}{2} \int_\Omega e^{\chi \tau} z^{p - 2} |\nabla z|^2 \\[5pt]
    & \leq C_1(p)  \int_\Omega z^p + (p - 1)\delta \chi e^{2 \chi M_\tau} M_\tau \int_\Omega z^{p + 1} + C_2
  \end{align}
for all $t \in (0, T_{\max})$. By means of Young's inequality, we directly infer
\begin{align}
    \label{Lpz434}
    \nonumber \frac{d}{dt} \int_\Omega e^{\chi \tau} z^{p} + \int_\Omega e^{\chi \tau} z^{p} & + \frac{2 D_1 (p - 1)}{p} \int_\Omega |\nabla z^{\frac p2}|\\[5pt]
    & \leq \underbrace{(C_1(p) + (p - 1)\delta \chi e^{2 \chi M_\tau} M_\tau)}_{=:\Gamma}\int_\Omega z^{p + 1} + C_3(p)
  \end{align}
for all $t \in (0, T_{\max})$. Invoking the Gagliardo-Nirenberg inequality from \cite[(4.24)]{dai2022new}, together with the $L^2$ estimate in \eqref{uwsql2}, we obtain
\begin{align}
\label{gag_3d33}
\nonumber \Gamma \|z\|_{L^{p+1}(\Omega)}^{p+1}
& = \Gamma \|z^{\frac p2}\|_{L^{\frac{2(p+1)}{p}}(\Omega)}^{\frac{2(p+1)}{p}} \\[5pt]
\nonumber &\leq \Gamma C_{\mathrm{GN}}^{(3)}\left(
\|\nabla z^{\frac p2}\|_{L^2(\Omega)}^{\frac{2\theta(p+1)}{p}}
\|z^{\frac p2}\|_{L^{\frac{4}{p}}(\Omega)}^{\frac{2(1-\theta)(p+1)}{p}}
+\|z^{\frac p2}\|_{L^{\frac{4}{p}}(\Omega)}^{\frac{2(p+1)}{p}}
\right) \\[5pt]
&\leq C_{4}\left(\|\nabla z^{\frac p2}\|_{L^2(\Omega)}^{\frac{2\theta(p+1)}{p}}+1\right),
\end{align}
where
$$
\theta=\frac{p(3p-3)}{(p+1)(3p-2)}\in(0,1),
$$
due to $p>1$, and 
$$
\theta\frac{2(p+1)}{p}=2-\frac{2}{3p-2}\in(0,2),
$$
for any $p>1$. Applying Young's inequality to \eqref{gag_3d33}, we obtain
\begin{align}
    \label{you1}
     \Gamma \int_\Omega z^{p+1} \leq  \frac{2D_1(p-1)}{p}\int_\Omega |\nabla z^{\frac p2}|^2 + \frac{1}{D_1}C_5(p), 
\end{align}
for all $t \in (0, T_{\max})$. Inserting \eqref{you1} into \eqref{Lpz434} results in 
\begin{align}
    \label{Lpz435}
    \frac{d}{dt} \int_\Omega e^{\chi \tau} z^{p} + \int_\Omega e^{\chi \tau} z^{p} \leq C_6(p)
  \end{align}
for all $t \in (0, T_{\max})$. Applying Gronwall's lemma to \eqref{Lpz435}, we obtain
\begin{align}
\label{pmain1}
\int_\Omega e^{\chi \tau}z^p
\leq C_7(p, T)
\end{align}
for all $t\in(0, T_{\max})$. The desired estimate \eqref{ulp} then follows immediately from \eqref{pmain1} and \eqref{tau_bound}.
\end{proof}

%%%%%%%%%%%%%%%%%%%%%%%%%%%%%%%%%%%%%%%%%%%%%%%%%%%%%%%%%%%%%%%%%%%%%%%%%%%%%%%%%%%%%%%%%%%%%%%%%%%%%%%%%%%%%%%%%%%%%%%%%%%%%%%%%%%%%%%%%%%%%%%%%%%%%%%%%%%%%%%%%%%%%%%%%%%%%%%%%%%%%%%%%%%%%%%%%%%%%%%%%%%%%%%%%%%%%%%%%%%%%%%%%%%%%%%%%%%%%%%%%%%%%%%%%%%%%%%%%%%%%%%%%%%%%%%%%%%%%%%%%%%%%%%%%%%%%%%%%%%%%%%%%%%%%%%%%%%%%%%%%%%%%%%%%%%%%%%%%%%%%%%%%%%%%%%%%%%%%%%%%%%%%%%%%%%%%%%%%%%%%%%%%%%%%%%%%%%%%%%%%%%%%%%%%%%%%%%%%%%%
Having established $L^p(\Omega)$ estimates for $z$ for every $p>1$, we are now in a position to derive an $L^\infty(\Omega)$ bound for $z$. This is achieved by applying the Moser-Alikakos iteration technique (cf. \cite[Lemma 4.4]{dai2020global}).
\begin{lemma}
    \label{linf3d}
     Let the assumptions of Theorem~\ref{3D_Class_Theorem} be satisfied. Then there exists a constant $C(T)>0$ such that
    \begin{align}
    \label{ulinf}
    \|z(\cdot,t)\|_{L^\infty(\Omega)}
    \le C(T)
    \qquad \text{for all } t\in(0, T_{\max}).
\end{align}
\end{lemma}
%%%%%%%%%%%%%%%%%%%%%%%%%%%%%%%%%%%%%%%%%%%%%%%%%%%%%%%%%%%%%%%%%%%%%%%%%%%%%%%%%%%%%%%%%%%%%%%%%%%%%%%%%%%%%%%%%%%%%%%%%%%%%%%%%%%%%%%%%%%%%%%%%%%%%%%%%%%%%%%%%%%%%%%%%%%%%%%%%%%%%%%%%%%%%%%%%%%%%%%%%%%%%%%%%%%%%%%%%%%%%%%%%%%%%%%%%%%%%%%%%%%%%%%%%%%%%%%%%%%%%%%%%%%%%%%%%%%%%%%%%%%%%%%%%%%%%%%%%%%%%%%%%%%%%%%%%%%%%%%%%%%%%%%%%%%%%%%%%%%%%%%%%%%%%%%%%%%%%%%%%%%%%%%%%%%%%%%%%%%%%%%%%%%%%%%%%%%%%%%%%%%%%%%%%%%%%%%%%%%%
Since all three solution components are bounded in $L^\infty(\Omega)$ and the Gagliardo-Nirenberg inequality \eqref{gag_ninnn11} remains valid in three dimensions, the estimate \eqref{gag_niren_1} continues to hold. Consequently, the argument leading to \eqref{lemma_nabL4_bound} carries over to the three-dimensional setting, yielding, for every $T>0$,
\begin{align}
\label{tau_c3ess}
\|\nabla \tau(\cdot,t)\|_{L^4(\Omega)}
\leq C
\qquad \text{for all } t\in(0, \widehat{T}),
\end{align}
 where $\widehat{T} := \min\{T, T_{\max}\}$.
%%%%%%%%%%%%%%%%%%%%%%%%%%%%%%%%%%%%%%%%%%%%%%%%%%%%%%%%%%%%%%%%%%%%%%%%%%%%%%%%%%%%%%%%%%%%%%%%%%%%%%%%%%%%%%%%%%%%%%%%%%%%%%%%%%%%%%%%%%%%%%%%%%%%%%%%%%%%%%%%%%%%%%%%%%%%%%%%%%%%%%%%%%%%%%%%%%%%%%%%%%%%%%%%%%%%%%%%%%%%%%%%%%%%%%%%%%%%%%%%%%%%%%%%%%%%%%%%%%%%%%%%%%%%%%%%%%%%%%%%%%%%%%%%%%%%%%%%%%%%%%%%%%%%%%%%%%%%%%%%%%%%%%%%%%%%%%%%%%%%%%%%%%%%%%%%%%%%%%%%%%%%%%%%%%%%%%%%%%%%%%%%%%%%%%%%%%%%%%%%%%%%%%%%%%%%%%%%%%%%

\subsubsection{Global extensibility in three dimensions}

We are now in a position to prove this section's second main result, asserting global classical solvability for \eqref{model1} in three dimensions.
\begin{proof}[Proof of Theorem \ref{3D_Class_Theorem}]
  Thanks to Lemma \ref{local_existence}, particularly the extensibility criterion \eqref{extensibility_criterion} therein, combining \eqref{ulinf} and \eqref{tau_c3ess}, yields the desired outcome.
\end{proof}

%%%%%%%%%%%%%%%%%%%%%%%%%%%%%%%%%%%%%%%%%%%%%%%%%%%%%%%%%%%%%%%%%%%%%%%%%%%%%%%%%%%%%%%%%%%%%%%%%%%%%%%%%%%%%%%%%%%%%%%%%%%%%%%%%%%%%%%%%%%%%%%%%%%%%%%%%%%%%%%%%%%%%%%%%%%%%%%%%%%%%%%%%%%%%%%%%%%%%%%%%%%%%%%%%%%%%%%%%%%%%%%%%%%%
%%%%%%%%%%%%%%%%%%%%%%%%%%%%%%%%%%%%%%%%%%%%%%%%%%%%%%%%%%%%%%%%%%%%%%%%%%%%%%%%%%%%%%%%%%%%%%%%%%%%%%%%%%%%%%%%%%%%%%%%%%%%%%%%%%%%%%%%%%%%%%%%%%%%%%%%%%%%%%%%%%%%%%%%%%%%%%%%%%%%%%%%%%%%%%%%%%%%%%%%%%%%%%%%%%%%%%%%%%%%%%%%%%%%
%%%%%%%%%%%%%%%%%%%%%%%%%%%%%%%%%%%%%%%%%%%%%%%%%%%%%%%%%%%%%%%%%%%%%%%%%%%%%%%%%%%%%%%%%%%%%%%%%%%%%%%%%%%%%%%%%%%%%%%%%%%%%%%%%%%%%%%%%%%%%%%%%%%%%%%%%%%%%%%%%%%%%%%%%%%%%%%%%%%%%%%%%%%%%%%%%%%%%%%%%%%%%%%%%%%%%%%%%%%%%%%%%%%%

\section{\texorpdfstring
{Global existence of weak solutions to \eqref{model}-\eqref{initial_conditions} under assumptions \ref{weak_Dass} on $\mathbb{D}$ in two spatial dimensions}
{weak solutions}}\label{sec:global-weak-deg}

In the previous section, we established the global existence of classical solutions for \eqref{model}, \eqref{boundary_conditions} and \eqref{initial_conditions} under certain stringent conditions on $\mathbb{D}$ in two dimensions. Here, our goal is to construct weak solutions for the same problem under comparatively weaker assumptions on the diffusion tensor. {The proof proceeds along the path in \cite{heihoff2023global}, with necessary modifications required by the increased complexity of our PDE-ODE-ODE system.} We would like to mention that the notion of function spaces, notably $W_{\mathbb{D}}^{1, 2}(\Omega)$ and $W^{1, p}_{\textup{div}}(\Omega;\mathbb{R}^{n \times n})$ used in this section have been discussed in \cite[Definition 2.8]{heihoff2023global}; for a detailed coverage of these spaces refer to \cite{sawyer2010degenerate}.

\subsection{Assumptions on \texorpdfstring{$\mathbb{D}$}{D} : weak solutions (degenerate case)}
\label{weak_Dass}
We begin by presenting a key result, laying the foundation for the assumptions that follow:

\begin{definition}[{\cite[Definition 2.10]{heihoff2023global}}]
  \label{the_embedding}
   Let $\Omega \subseteq \mathbb{R}^n, n \subseteq \mathbb{N}$, be a bounded domain with a smooth boundary. We say a positive semidefinite $\mathbb{D} \in C^0(\overline{\Omega}; \mathbb{R}^{n \times n})$ allows for a compact $L^1(\Omega)$ embedding if $W^{1,2}_{\mathbb{D}}(\Omega)$ embeds compactly into $L^1(\Omega)$.
\end{definition}

Then we fix $\Omega \subseteq \mathbb{R}^n ~(n = 2)$ as a smooth bounded domain. As mentioned earlier all the parameters $\alpha_1, \alpha_2, \beta, \delta, \mu$ and $\chi$ are positive.  We also take some a.e. nonnegative initial values $c_{10} \in L^{v[\ln(v)]_+}(\Omega)$ and $c_{20}, \tau_0 \in L^\infty(\Omega)$ with $\sqrt{c_{20}}, \sqrt{\tau_0} \in W^{1, 2}(\Omega)$, where $L^{v[\ln(v)]_+}(\Omega)$ is the standard Orlicz space associated with the function $v \mapsto v[\ln(v)]_+$.

We fix $\mathbb{D} \in W^{1, 2}_{\textup{div}}(\Omega; \mathbb{R}^{n \times n}) \cap C^0(\overline{\Omega}; \mathbb{R}^{n \times n})$ with the following properties
\begin{itemize}
  \item $\mathbb{D}$ is positive semidefinite everywhere.
  \item The divergence estimates of Definition \ref{req_weakSol} hold true.
  \item A positive semidefinite $\mathbb{D} \in C^0(\overline{\Omega}; \mathbb{R}^{n \times n})$ allows for a compact $L^1(\Omega)$ embedding as per Definition \ref{the_embedding}.
\end{itemize}
Using the divergence estimate of Definition \ref{req_weakSol} we can also have
\begin{align*}
  \mathbb{D} \in W^{1, q}_{\textup{div}}(\Omega; \mathbb{R}^{n \times n}) \subseteq W^{1, 2}_{\textup{div}}(\Omega; \mathbb{R}^{n \times n})
\end{align*}
with $q := \tfrac{2 \lambda}{2 \lambda - 1}$. We will also invoke the following lemma which fixes an approximate family $(\mathbb{D}_\varepsilon)_{\varepsilon \in (0, 1)}$:

\begin{lemma}[{\cite[Lemma 4.1]{heihoff2023global}}]
 \label{eplemma1}
  There exists an approximate family $(\mathbb{D}_\varepsilon)_{\varepsilon \in (0,1)} \subseteq C^2(\Omega; \mathbb{R}^{n \times n})$ with $\mathbb{D}_\varepsilon$ positive definite on $\overline{\Omega}$, $(\nabla \cdot \mathbb{D}_\varepsilon) \cdot \nu = 0$ on $\partial \Omega$ for all $\varepsilon \in (0,1)$ and
  \begin{align}
    \label{epassumption1}
    \mathbb{D}_\varepsilon \to \mathbb{D} \quad \text{in} \quad W^{1,q}_{\textup{div}}(\Omega; \mathbb{R}^{n \times n}) \cap C^0(\overline{\Omega}; \mathbb{R}^{n \times n}) \quad \text{as} \quad \varepsilon \searrow 0.
  \end{align}
  We can further choose this family in such a way as to ensure that
  \begin{align}
    \label{epassumption2}
    \left| \int_{\Omega} (\nabla \cdot \mathbb{D}_\varepsilon) \cdot \Phi \right| \leq B \left( \int_{\Omega} (\Phi \cdot \mathbb{D}_\varepsilon \Phi)^\lambda + 1 \right)
  \end{align}
  with $B := A + 1$ and with exponent $\lambda \in [\tfrac{1}{2}, \tfrac{2}{3}]$ such that $\tfrac{\lambda}{1 - \lambda} \leq 2$, and
  \begin{align}
     \label{epassumption3}
     \mathbb{D} + \varepsilon \leq \mathbb{D}_\varepsilon \leq \mathbb{D} + 3\varepsilon
  \end{align}
  for all $\Phi \in C^0(\overline{\Omega}; \mathbb{R}^n)$ and $\varepsilon \in (0,1)$.
\end{lemma}
Here the inequality of tensors is meant in the usual sense: for $\mathbb A,\ \mathbb B\in \mathbb R^{n\times n}$,  $\mathbb A\le \mathbb B$ iff $\mathbb B-\mathbb A$ is positive semidefinite.

%%%%%%%%%%%%%%%%%%%%%%%%%%%%%%%%%%%%%%%%%%%%%%%%%%%%%%%%%%%%%%%%%%%%%%%%%%%%%%%%%%%%%%%%%%%%%%%%%%%%%%%%%%%%%%%%%%%%%%%%%%%%%%%%%%%%%%%%%%%%%%%%%%%%%%%%%%%%%%%%%%%%%%%%%%%%%%%%%%%%%%%%%%%%%%%%%%%%%%%%%%%%%%%%%%%%%%%%%%%%%%%%%
%%%%%%%%%%%%%%%%%%%%%%%%%%%%%%%%%%%%%%%%%%%%%%%%%%%%%%%%%%%%%%%%%%%%%%%%%%%%%%%%%%%%%%%%%%%%%%%%%%%%%%%%%%%%%%%%%%%%%%%%%%%%%%%%%%%%%%%%%%%%%%%%%%%%%%%%%%%%%%%%%%%%%%%%%%%%%%%%%%%%%%%%%%%%%%%%%%%%%%%%%%%%%%%%%%%%%%%%%%%%%%%%%
\subsection{Main result: weak solutions}
As the assumptions on the diffusion tensor required for constructing weak solutions have been established, we now introduce the main result of this section. Notably, in the context of our intended compactness properties for the first solution component, it is crucial to express its gradient in the following manner, as in \cite{stinner2014global, stinner2016global}:
\begin{equation}
\label{c1_nab_strut1}
\nabla c_1 = 2 \sqrt{1 + c_1} \nabla \sqrt{1 + c_1}.
\end{equation}
With \eqref{c1_nab_strut1} we now present the following definition of weak solutions for \eqref{model}, \eqref{boundary_conditions} and \eqref{initial_conditions}:
\begin{definition}
    \label{def_weak_sol}
    Let $\Omega \subseteq \mathbb{R}^2$ be a bounded domain with a smooth boundary and assume that the parameters $\alpha_1, \alpha_2, \beta, \delta, \mu$ and $\chi$ are all positive. Let $\mathbb{D} \in W^{1, q}_{\textup{div}}(\Omega; \mathbb{R}^{2 \times 2}) \cap C^0(\overline{\Omega}; \mathbb{R}^{2 \times 2})$, $q \in (1, \infty)$ and $p := \max\{2, \frac{q}{q - 1}\}$. Furthermore, let $c_{10}, c_{20}, \tau_0 \in L^1(\Omega)$ be some initial data. We then call a triple  of nonnegative functions $(c_1, c_2, \tau)$ a weak solution to problem \eqref{model}, \eqref{boundary_conditions}, and \eqref{initial_conditions}, if
     \begin{equation*}
  \begin{cases}
     c_1 \in L^2_{\textup{loc}}(\overline{\Omega} \times (0, \infty)) \quad \text{with} \quad \sqrt{1 + c_1} \in L^{\frac{4}{3}}_{\textup{loc}}([0, \infty); W_{\mathbb{D}}^{1, \frac{4}{3}}(\Omega)), \\[5pt]
     c_2 \in L^\infty(\Omega \times (0, \infty)) \cap L^2_{\textup{loc}}([0, \infty); W^{1, 2}(\Omega)),                                                                                          \\[5pt]
    \tau \in L^\infty(\Omega \times (0, \infty)) \cap L^2_{\textup{loc}}([0, \infty); W_{\mathbb{D}}^{1, 2}(\Omega))
  \end{cases}
  \end{equation*}
  and the triple satisfies
  \begin{align}
    \label{weak_c1}
    \nonumber & \int_0^\infty \int_\Omega c_1 \partial_t \varphi + \int_\Omega c_{10} \varphi(\cdot, 0) = 2 \int_0^\infty \int_\Omega \sqrt{1 + c_1} \nabla \sqrt{1 + c_1} \cdot \mathbb{D} \nabla \varphi \\[5pt]
    \nonumber & + \int_0^\infty \int_\Omega c_1 (\nabla \cdot \mathbb{D}) \cdot \nabla \varphi - \chi \int_0^\infty \int_\Omega c_1 \nabla \tau \cdot \mathbb{D} \nabla \varphi + \alpha_1  \int_0^\infty \int_\Omega \tfrac{c_{1}}{1 + c_{1}} \varphi  \\[5pt]
    & - \alpha_2  \int_0^\infty \int_\Omega \tfrac{c_{2}}{1 + c_{2}}\varphi - \beta \int_0^\infty \int_\Omega c_{1}(1 - c_{1} - c_{2} - \tau)\varphi,
  \end{align}
  and
  \begin{align}
    \label{weak_c2}
    \int_0^\infty \int_\Omega c_2 \partial_t \varphi + \int_\Omega c_{20} \varphi(\cdot, 0) = -\alpha_1 \int_0^\infty \int_\Omega \tfrac{c_1}{1 + c_1} \varphi + \alpha_2 \int_0^\infty \int_\Omega \tfrac{c_2}{1 + c_2}\varphi,
  \end{align}
  as well as
  \begin{align}
    \label{weak_tau}
    \int_0^\infty \int_\Omega \tau \partial_t \varphi + \int_\Omega \tau_{0} \varphi(\cdot, 0) =  \delta \int_0^\infty \int_\Omega c_1 \tau \varphi + \mu \int_0^\infty \int_\Omega \tau \varphi - \int_0^\infty \int_\Omega \tfrac{c_2 }{1 + c_2}\varphi,
  \end{align}
  for all $\varphi \in C_c^\infty(\overline{\Omega} \times [0, \infty))$.
\end{definition}

The main result of this section is stated as follows:
\begin{theorem}
    \label{main_th_weak_sol}
     Let $\Omega \subseteq \mathbb{R}^2$ be a bounded domain with a smooth boundary and assume that the parameters $\alpha_1, \alpha_2, \beta, \delta, \mu$ and $\chi$ are all positive and $\lambda \in [\tfrac{1}{2}, \tfrac{2}{3}]$ with $\tfrac{\lambda}{1 - \lambda} \leq 2$ and $\mathbb{D} \in W_{\textup{div}}^{1, 2}(\Omega; \mathbb{R}^{2 \times 2}) \cap C^0(\overline{\Omega}; \mathbb{R}^{2 \times 2})$ be symmetric and positive semidefinite everywhere. Moreover, let $\mathbb{D}$ satisfy the divergence estimate of Definition \ref{req_weakSol} with exponent $\lambda$ and allow for a compact $L^1(\Omega)$ embedding of Definition \ref{the_embedding}. Furthermore, let $c_{10} \in L^{v[\ln v]_+}(\Omega)$ and $c_{20}, \tau_0 \in L^\infty(\Omega)$ be some initial data with $\sqrt{c_{20}}, \sqrt{\tau_0} \in W^{1, 2}(\Omega)$ and $c_{10} > 0, c_{20} > 0, \tau_0 > 0$ almost everywhere. %Here $L^{v[\ln v]_+}(\Omega)$ is the standard Zygmund space associated with the functions $v \mapsto v[\ln (v)]_+$. 
     Then the problem \eqref{model}, \eqref{boundary_conditions}, and \eqref{initial_conditions} admits a global weak solution in the sense of Definition \ref{def_weak_sol}.
\end{theorem}
%%%%%%%%%%%%%%%%%%%%%%%%%%%%%%%%%%%%%%%%%%%%%%%%%%%%%%%%%%%%%%%%%%%%%%%%%%%%%%%%%%%%%%%%%%%%%%%%%%%%%%%%%%%%%%%%%%%%%%%%%%%%%%%%%%%%%%%%%%%%%%%%%%%%%%%%%%%%%%%%%%%%%%%%%%%%%%%%%%%%%%%%%%%%%%%%%%%%%%%%%%%%%%%%%%%%%%%%%%%%%%%%%
%%%%%%%%%%%%%%%%%%%%%%%%%%%%%%%%%%%%%%%%%%%%%%%%%%%%%%%%%%%%%%%%%%%%%%%%%%%%%%%%%%%%%%%%%%%%%%%%%%%%%%%%%%%%%%%%%%%%%%%%%%%%%%%%%%%%%%%%%%%%%%%%%%%%%%%%%%%%%%%%%%%%%%%%%%%%%%%%%%%%%%%%%%%%%%%%%%%%%%%%%%%%%%%%%%%%%%%%%%%%%%%%%

\subsection{A family of approximate problems}
Having fixed an approximate family $(\mathbb{D}_\varepsilon)_{\varepsilon \in (0, 1)}$ for $\mathbb{D}$, as well as the associated constant $B$ as described in Lemma \ref{eplemma1}, we will now proceed to construct a family of approximate problems for our original problem in \eqref{model}, \eqref{boundary_conditions} and \eqref{initial_conditions}.

First, we fix some families $(c_{10\varepsilon})_{\varepsilon \in (0, 1)}, ~ (c_{20\varepsilon}')_{\varepsilon \in (0, 1)}, ~ (\tau_{0\varepsilon}')_{\varepsilon \in (0, 1)} \in C^3(\overline{\Omega})$ of positive functions with 
\begin{align*}
    (\mathbb{D}_\varepsilon \nabla c_{10\varepsilon}) \cdot \nu = (\mathbb{D}_\varepsilon \nabla \tau_{0\varepsilon}') \cdot \nu  = 0 \quad \text{on}~ \partial \Omega
\end{align*}
and
\begin{align*}
    \|c_{20\varepsilon}'\|_{L^\infty(\Omega)} \leq \|\sqrt{c_{20}}\|_{L^\infty(\Omega)} + 1, \qquad \|\tau_{0\varepsilon}'\|_{L^\infty(\Omega)} \leq \|\sqrt{\tau_{0}}\|_{L^\infty(\Omega)} + 1
\end{align*}
for all $\varepsilon \in (0, 1)$ as well as
\begin{align*}
  c_{10\varepsilon} \to c_{10} \quad \text{in}~ L^{v [\ln(v)]_+}(\Omega), \quad c_{20\varepsilon}' \to \sqrt{c_{20}} \quad \text{in}~ W^{1, 2}(\Omega), \quad \tau_{0\varepsilon}' \to \sqrt{\tau_{0}} \quad \text{in}~ W^{1, 2}(\Omega) 
\end{align*}
as $\varepsilon \searrow 0$. 

Now, we let $c_{20\varepsilon} := (c_{20\varepsilon}')^2 \in C^3(\overline{\Omega})$ and $\tau_{0\varepsilon} := (\tau_{0\varepsilon}')^2 \in C^3(\overline{\Omega})$ for all $\varepsilon \in (0, 1)$. Due to the properties established for the family $(c_{20\varepsilon}')_{\varepsilon \in (0, 1)}$ and $(\tau_{0\varepsilon}')_{\varepsilon \in (0, 1)}$, we can clearly have $c_{20\varepsilon} > 0$ and $\tau_{0\varepsilon} > 0$ on $\overline{\Omega}$, $(\mathbb{D}_\varepsilon \nabla \tau_{0\varepsilon}) \cdot \nu = 0$ on $\partial \Omega$ and
\begin{align*}
  & c_{20\varepsilon} \to c_{20} \quad ~\text{in}~ L^p(\Omega)~\text{for all}~ p \in [1, \infty), \quad \sqrt{c_{20\varepsilon}} \to \sqrt{c_{20}} \quad \text{in}~ W^{1,2}(\Omega), \\[5pt]
  & \tau_{0\varepsilon} \to \tau_{0} \quad ~\text{in}~ L^p(\Omega)~\text{for all}~ p \in [1, \infty), \quad \sqrt{\tau_{0\varepsilon}} \to \sqrt{\tau_{0}} \quad \text{in}~ W^{1,2}(\Omega)
\end{align*}
as $\varepsilon \searrow 0$. From the above assumptions, we can find a constant $M > 0$ such that
\begin{align}
  \label{weak_assup1}
  \|\nabla \cdot \mathbb{D}_\varepsilon\|_{L^2(\Omega)} \leq M, \quad  \|\mathbb{D}_\varepsilon\|_{L^\infty(\Omega)} \leq M
\end{align}
and
\begin{equation}
\begin{aligned}
  \label{weak_assup2}
  \int_\Omega c_{10\varepsilon} \leq M, \quad & \int_\Omega c_{10\varepsilon} \ln c_{10\varepsilon}  \leq M, \quad \|c_{20\varepsilon}\|_{L^\infty(\Omega)} \leq M, \quad \int_\Omega \frac{|\nabla c_{20\varepsilon}|^2}{c_{20\varepsilon}} \leq M , \\[5pt]
  & \|\tau_{0\varepsilon}\|_{L^\infty(\Omega)} \leq M, \quad \int_\Omega \frac{\nabla \tau_{0\varepsilon} \cdot \mathbb{D} \nabla \tau_{0\varepsilon}}{\tau_{0\varepsilon}} \leq M
\end{aligned}
\end{equation}
for all $\varepsilon \in (0, 1)$.

We are now in a position to introduce the following family of approximate problems for \eqref{model}, \eqref{boundary_conditions} and \eqref{initial_conditions}:
\begin{equation}
  \begin{cases}
    \label{app_prob}
     \partial_t c_{1\varepsilon} = \nabla \cdot (\mathbb{D}_\varepsilon \nabla c_{1\varepsilon} + c_{1\varepsilon} \nabla \cdot \mathbb{D}_\varepsilon) - \chi \nabla \cdot (c_{1\varepsilon} \mathbb{D}_\varepsilon \nabla \tau_\varepsilon)\\[10pt]
     \hspace{1cm} - \alpha_1  \tfrac{c_{1\varepsilon}}{1 + c_{1\varepsilon}}  + \alpha_2 \tfrac{c_{2\varepsilon}}{1 + c_{2\varepsilon}} + \beta c_{1\varepsilon}(1 - c_{1\varepsilon} - c_{2\varepsilon} - \tau_\varepsilon),                                                                   &  ~~ x \in \Omega, t > 0,  \\[7pt]
    \partial_t c_{2\varepsilon} =  \alpha_1 \frac{c_{1\varepsilon}}{1 + c_{1\varepsilon}} - \alpha_2 \frac{c_{2\varepsilon}}{1 + c_{2\varepsilon}},                        & ~~ x \in \Omega, t > 0,   \\[7pt]
    \partial_t \tau_\varepsilon = -\delta c_{1\varepsilon} \tau_\varepsilon - \mu \tau_\varepsilon + \frac{c_{2\varepsilon} }{1 + c_{2\varepsilon}},                       & ~~ x \in \Omega, t > 0,   \\[7pt]
    (\mathbb{D}_\varepsilon \nabla c_{1\varepsilon}) \cdot \nu = \chi (c_{1\varepsilon} \mathbb{D}_\varepsilon \nabla \tau_\varepsilon)  \cdot \nu- c_{1\varepsilon}(\nabla \cdot \mathbb{D}_\varepsilon) \cdot \nu = 0,                                                                                                        & ~~ x \in \partial \Omega, t > 0,   \\[5pt]
    c_{1\varepsilon}(\cdot, 0) = c_{10\varepsilon}, ~~ c_{2\varepsilon}(\cdot, 0) = c_{20\varepsilon}, ~~ \tau_{\varepsilon}(\cdot, 0) = \tau_{0\varepsilon},              & ~~ x \in \Omega,
  \end{cases}
\end{equation}
for all $\varepsilon \in (0, 1)$.
%%%%%%%%%%%%%%%%%%%%%%%%%%%%%%%%%%%%%%%%%%%%%%%%%%%%%%%%%%%%%%%%%%%%%%%%%%%%%%%%%%%%%%%%%%%%%%%%%%%%%%%%%%%%%%%%%%%%%%%%%%%%%%%%%%%%%%%%%%%%%%%%%%%%%%%%%%%%%%%%%%%%%%%%%%%%%%%%%%%%%%%%%%%%%%%%%%%%%%%%%%%%%%%%%%%%%%%%%%%%%%%%%
%%%%%%%%%%%%%%%%%%%%%%%%%%%%%%%%%%%%%%%%%%%%%%%%%%%%%%%%%%%%%%%%%%%%%%%%%%%%%%%%%%%%%%%%%%%%%%%%%%%%%%%%%%%%%%%%%%%%%%%%%%%%%%%%%%%%%%%%%%%%%%%%%%%%%%%%%%%%%%%%%%%%%%%%%%%%%%%%%%%%%%%%%%%%%%%%%%%%%%%%%%%%%%%%%%%%%%%%%%%%%%%%%

\subsection{Global existence for approximate problems and other results}

\begin{lemma}
  \label{global_existence_app}
For all $\varepsilon \in (0, 1)$ there exists a global classical solution to \eqref{app_prob}.
\end{lemma}

\begin{proof}
  The classical existence theory developed in Theorem \ref{classical_main_th1} can be applied in this context for any $\varepsilon \in (0, 1)$. Thus, a global classical solution $(c_{1\varepsilon}, c_{2\varepsilon}, \tau_\varepsilon)$ to the regularized system \eqref{app_prob} exists for all $\varepsilon \in (0, 1)$.
\end{proof}

From Lemma \ref{global_existence_app}, it is evident that the estimates derived for \eqref{model1} in Lemmas \ref{lemmal1c1}, \ref{lemma_c2_bound} and \ref{lemma_tau_bound} also hold for \eqref{app_prob} for all $\varepsilon \in (0, 1)$, that is:
\begin{lemma}
\label{lbound_comp}
  Let $T > 0$. Then there exist positive constants $\widehat{C}, \widehat{M}_1, \widehat{M}_{c_2}$ and $\widehat{M}_\tau$ such that for all $\varepsilon \in (0, 1)$
  \begin{align}
    \label{bound_comp1}
    & \int_\Omega c_{1\varepsilon}(\cdot, t) \leq \widehat{M}_1, \quad \text{for all}~ t \in (0, T_{\max, \varepsilon}), \\[5pt]
    \label{bound_comp2}
    & \|c_{2\varepsilon}(\cdot, t)\|_{L^\infty(\Omega)} \leq \widehat{M}_{c_2} \quad \text{for all}~ t \in (0, \widehat{T}_\varepsilon),   \text{ and} \\[5pt]
      \label{bound_comp3}
    &  \|\tau_{\varepsilon}(\cdot, t)\|_{L^\infty(\Omega)} \leq \widehat{M}_\tau \quad \text{for all}~ t \in (0, T_{\max,\varepsilon}),
  \end{align}
  where $\widehat{T}_\varepsilon := \min\{T, T_{\max, \varepsilon}\}$. Also, we have
\begin{equation}
\label{bound_comp4}
\int_t^{t+\varsigma}\int_\Omega c_{1\varepsilon}^2
\leq \widehat{C}
\end{equation}
for all $t\geq0$ and $0<\varsigma\leq1$ such that $t+\varsigma\leq T_{\max,\varepsilon}$.

\end{lemma}

Moreover, the energy-type inequality established in Lemma \ref{entropy_lemma} is also valid in this context. We restate it in the following lemma.

\begin{lemma}
\label{entropy_weak:lemma}
Let $T > 0$. Then there exists a  $C(T) > 0$ such that for all $\varepsilon \in (0, 1)$ any solution $(c_{1\varepsilon}, c_{2\varepsilon}, \tau_\varepsilon)$ to \eqref{app_prob} satisfies
\begin{align}
\label{weentropy0}
&\sup_{t\in(0,T)}
\left\{\int_\Omega\left(c_{1\varepsilon}\ln c_{1\varepsilon}+\frac{1}{e}\right)
+\frac{\chi}{2\delta}\int_\Omega
\frac{\nabla\tau_\varepsilon\cdot\mathbb{D}_\varepsilon\nabla\tau_\varepsilon}
{\tau_\varepsilon}
+\frac{D_1}{8\alpha_1}\int_\Omega
\frac{|\nabla c_{2\varepsilon}|^2}{c_{2\varepsilon}}
\right\}
\nonumber\\[7pt]
&\quad
+\frac18\int_0^T\int_\Omega
\frac{\nabla c_{1\varepsilon}\cdot\mathbb{D}_\varepsilon\nabla c_{1\varepsilon}}
{c_{1\varepsilon}}
+\frac{\beta}{4}\int_0^T\int_\Omega
c_{1\varepsilon}^2(\ln2+c_{1\varepsilon})
\nonumber\\[7pt]
&\quad
+\frac{\chi}{2}\int_0^T\int_\Omega
\frac{c_{1\varepsilon}\nabla\tau_\varepsilon\cdot\mathbb{D}_\varepsilon\nabla\tau_\varepsilon}
{\tau_\varepsilon}
+\frac{\mu\chi}{2\delta}\int_0^T\int_\Omega
\frac{\nabla\tau_\varepsilon\cdot\mathbb{D}_\varepsilon\nabla\tau_\varepsilon}
{\tau_\varepsilon}
\le C(T).
\end{align}
\end{lemma}

%%%%%%%%%%%%%%%%%%%%%%%%%%%%%%%%%%%%%%%%%%%%%%%%%%%%%%%%%%%%%%%%%%%%%%%%%%%%%%%%%%%%%%%%%%%%%%%%%%%%%%%%%%%%%%%%%%%%%%%%%%%%%%%%%%%%%%%%%%%%%%%%%%%%%%%%%%%%%%%%%%%%%%%%%%%%%%%%%%%%%%%%%%%%%%%%%%%%%%%%%%%%%%%%%%%%%%%%%%%%%%%%%
%%%%%%%%%%%%%%%%%%%%%%%%%%%%%%%%%%%%%%%%%%%%%%%%%%%%%%%%%%%%%%%%%%%%%%%%%%%%%%%%%%%%%%%%%%%%%%%%%%%%%%%%%%%%%%%%%%%%%%%%%%%%%%%%%%%%%%%%%%%%%%%%%%%%%%%%%%%%%%%%%%%%%%%%%%%%%%%%%%%%%%%%%%%%%%%%%%%%%%%%%%%%%%%%%%%%%%%%%%%%%%%%%

\subsection{Some \texorpdfstring{$\varepsilon$}{epsilon}-independent estimates}

Apart from the bounds in Lemma \ref{lbound_comp}, the limit procedure required for constructing weak solutions to \eqref{app_prob} requires some additional estimates. These will be derived from the energy-type functional of Lemma \ref{entropy_weak:lemma}.

\subsubsection{Estimates for \texorpdfstring{ $c_{2\varepsilon}$ }{c2} and  \texorpdfstring{ $\tau_{\varepsilon}$ }{c2}}

\begin{lemma}
\label{lemma_epsilon1}
  Let $T > 0$. Then, there exists a $C(T) > 0$ such that for all $\varepsilon \in (0,1)$, the following hold true
  \begin{align}
    \label{nabla_c21} \int_\Omega \|c_{2\varepsilon}(\cdot, t)\|^2_{W^{1,2}(\Omega)} \leq C(T), \quad \text{for all}~ t \in (0, T),
   \end{align}
   and
   \begin{align}
    \label{nabla_tau_1} \int_0^T \|\tau_\varepsilon(\cdot, s)\|^2_{W_{\mathbb{D}}^{1, 2}(\Omega)} ds ~\leq~ \int_0^T \|\tau_\varepsilon(\cdot, s)\|^2_{W_{\mathbb{D}_\varepsilon}^{1, 2}(\Omega)} ds \leq C(T).
    \end{align}
  \end{lemma}

\begin{proof}
  Fix $T > 0$. The estimate in \eqref{nabla_c21} can be easily obtained by using \eqref{weentropy0} and  \eqref{bound_comp2} in the following manner:
  \begin{equation*}
    \int_\Omega |\nabla c_{2\varepsilon}|^2 = \int_\Omega \frac{|\nabla c_{2\varepsilon}|^2}{c_{2\varepsilon}}c_{2\varepsilon} \leq \widehat{M}_{c_2} \int_\Omega \frac{|\nabla c_{2\varepsilon}|^2}{c_{2\varepsilon}} \leq C(T), \quad \text{for all}~ t \in (0, T).
  \end{equation*}
  To verify \eqref{nabla_tau_1}, by using \eqref{weentropy0} and \eqref{bound_comp3}, we can easily obtain a constant $M_{2\tau} > 0$ such that the following holds true
  \begin{equation}
  \label{nab_tau_bound11}
    \int_\Omega \nabla \tau_\varepsilon \cdot \mathbb{D}_\varepsilon \nabla \tau_\varepsilon \leq \widehat{M}_\tau \int_\Omega \frac{\nabla \tau_\varepsilon \cdot \mathbb{D}_\varepsilon \nabla \tau_\varepsilon}{\tau_\varepsilon} =: M_{2\tau}
  \end{equation}
  for all $t \in (0, T)$ and all $\varepsilon \in (0,1)$. From \eqref{epassumption3}, we know that $\mathbb{D} \leq \mathbb{D}_\varepsilon$ for all $\varepsilon \in (0,1)$. Then, combining \eqref{weentropy0} with \eqref{bound_comp3} and the fact that $\mathbb{D} \leq \mathbb{D}_\varepsilon$ for all $\varepsilon \in (0,1)$ yields \eqref{nabla_tau_1}.
\end{proof}

Next, we derive uniform integrability estimates for the time derivatives of $c_{2\varepsilon}$ and $\tau_\varepsilon$.

\begin{lemma}
  \label{time_det_c2tau}
    Let $T > 0$. Then, there exists a $C(T) > 0$ such that for all $\varepsilon \in (0,1)$, the following hold true:
    \begin{align}
      \label{time_tau1}
      \int_0^T \|\partial_t c_{2\varepsilon}(\cdot, t)\|_{(W^{3,2}(\Omega))^\ast}dt \leq C(T) \quad \text{and} \quad \int_0^T \|\partial_t \tau_{\varepsilon}(\cdot, t)\|_{(W^{3,2}(\Omega))^\ast}dt \leq C(T).
    \end{align}
\end{lemma}

\begin{proof}
  Fix $T > 0$. As the diffusion terms are absent in the second and third equations of \eqref{app_prob}, the bounds in \eqref{time_tau1} follow directly from the uniform boundedness of $\frac{c_{i\varepsilon}}{1 + c_{i\varepsilon}} < 1$, $c_{i\varepsilon} \geq 0, i \in \{1, 2\}$, the bounds in \eqref{bound_comp2}-\eqref{bound_comp3}, and the fact that $W^{3,2}(\Omega)$ embeds continuously into $L^\infty(\Omega)$.
\end{proof}

%%%%%%%%%%%%%%%%%%%%%%%%%%%%%%%%%%%%%%%%%%%%%%%%%%%%%%%%%%%%%%%%%%%%%%%%%%%%%%%%%%%%%%%%%%%%%%%%%%%%%%%%%%%%%%%%%%%%%%%%%%%%%%%%%%%%%%%%%%%%%%%%%%%%%%%%%%%%%%%%%%%%%%%%%%%%%%%%%%%%%%%%%%%%%%%%%%%%%%%%%%%%%%%%%%%%%%%%%%%%%%%%%
%%%%%%%%%%%%%%%%%%%%%%%%%%%%%%%%%%%%%%%%%%%%%%%%%%%%%%%%%%%%%%%%%%%%%%%%%%%%%%%%%%%%%%%%%%%%%%%%%%%%%%%%%%%%%%%%%%%%%%%%%%%%%%%%%%%%%%%%%%%%%%%%%%%%%%%%%%%%%%%%%%%%%%%%%%%%%%%%%%%%%%%%%%%%%%%%%%%%%%%%%%%%%%%%%%%%%%%%%%%%%%%%%

\subsubsection{Estimates for \texorpdfstring{$\sqrt{1 + c_{1\varepsilon}}$}{c1epsilon} }

%As we have obtained the necessary estimates for constructing weak solutions for the second and third solution components, we now turn our attention to the first.
We now turn our attention towards obtaining some estimates for the approximation of the first solution component.

\begin{lemma}
  \label{lemma_sqrt1}
  Let $T > 0$. Then there exists a $C(T) > 0$ such that for all $\varepsilon \in (0,1)$ %, the following hold true:
  \begin{align}
  \label{sqrt1}
    \int_0^T \left\|\sqrt{1 + c_{1 \varepsilon}(\cdot, t)}\right\|_{W^{1, 2}_{\mathbb{D}}(\Omega)}^{2}dt \leq C(T),
   \end{align}
   and
   \begin{align}
   \label{sqrt2}
    \int_0^T \left\|\sqrt{1 + c_{1 \varepsilon}(\cdot, t)}\right\|_{W^{1, \frac{4}{3}}_{\mathbb{D}}(\Omega)}^{\frac{4}{3}}dt \leq  \int_0^T \left\|\sqrt{1 + c_{1 \varepsilon}(\cdot, t)}\right\|_{W^{1, \frac{4}{3}}_{\mathbb{D}_\varepsilon}(\Omega)}^{\frac{4}{3}}dt \leq C(T).
  \end{align}
\end{lemma}

\begin{proof}
  Fix $T > 0$. Then from \eqref{weentropy0} we have for all $\varepsilon \in (0,1)$
  \begin{align}
  \label{sqrt3}
     \int_0^T \int_\Omega \frac{\nabla c_{1\varepsilon} \cdot \mathbb{D}_\varepsilon \nabla c_{1\varepsilon}}{c_{1\varepsilon}} \leq C_1.
  \end{align}
 Using \eqref{sqrt3} we get
 \begin{align}
 \label{sqrt4}
   \nonumber \int_0^T \int_\Omega \nabla \sqrt{1 + c_{1\varepsilon}} \cdot \mathbb{D}_\varepsilon \nabla \sqrt{1 + c_{1\varepsilon}} & = \frac{1}{4} \int_0^T \int_\Omega \frac{\nabla c_{1\varepsilon} \cdot \mathbb{D}_\varepsilon \nabla c_{1\varepsilon}}{1 + c_{1\varepsilon}} \\[5pt]
   & \leq \frac{1}{4} \int_0^T \int_\Omega \frac{\nabla c_{1\varepsilon} \cdot \mathbb{D}_\varepsilon \nabla c_{1\varepsilon}}{c_{1\varepsilon}} \leq \frac{1}{4} C_1,
 \end{align}
 thus yielding \eqref{sqrt1}. By an application of Young's inequality, we obtain
 \begin{align}
 \label{sqrt5}
   \nonumber \int_0^T \int_\Omega (\nabla \sqrt{1 + c_{1\varepsilon}} \cdot \mathbb{D}_\varepsilon \nabla \sqrt{1 + c_{1\varepsilon}})^{\frac{2}{3}} & = 2^{-\frac{4}{3}} \int_0^T \int_\Omega \frac{(\nabla c_{1\varepsilon} \cdot \mathbb{D}_\varepsilon \nabla c_{1\varepsilon})^{\frac{2}{3}}}{(1 + c_{1\varepsilon})^{\frac{2}{3}}} \\[5pt]
   \nonumber & \leq 2^{-\frac{4}{3}} \frac{2}{3}\int_0^T \int_\Omega \frac{\nabla c_{1\varepsilon} \cdot \mathbb{D}_\varepsilon \nabla c_{1\varepsilon}}{1 + c_{1\varepsilon}} + 2^{-\frac{4}{3}}3^{-1} |\Omega|T \\[5pt]
   \nonumber &  \leq 2^{-\frac{1}{3}}3^{-1} \int_0^T \int_\Omega \frac{\nabla c_{1\varepsilon} \cdot \mathbb{D}_\varepsilon \nabla c_{1\varepsilon}}{c_{1\varepsilon}} + 2^{-\frac{4}{3}}3^{-1} |\Omega|T \\[5pt]
   & \leq 2^{-\frac{1}{3}}3^{-1} C_1 + 2^{-\frac{2}{3}}3^{-1} |\Omega|T.
 \end{align}
 From \eqref{epassumption3}, we know that $\mathbb{D} \leq \mathbb{D}_\varepsilon$ for all $\varepsilon \in (0,1)$. Combining this with \eqref{sqrt5} yields the estimate in \eqref{sqrt2}.
\end{proof}

%%%%%%%%%%%%%%%%%%%%%%%%%%%%%%%%%%%%%%%%%%%%%%%%%%%%%%%%%%%%%%%%%%%%%%%%%%%%%%%%%%%%%%%%%%%%%%%%%%%%%%%%%%%%%%%%%%%%%%%%%%%%%%%%%%%%%%%%%%%%%%%%%%%%%%%%%%%%%%%%%%%%%%%%%%%%%%%%%%%%%%%%%%%%%%%%%%%%%%%%%%%%%%%%%%%%%%%%%%%%%%%%%
Next, to establish the strong compactness property for $\sqrt{1 + c_{1\varepsilon}(\cdot, t)}$, we need to appropriately control its time derivative. The following lemma presents this result.

\begin{lemma}
  \label{time_sqrt0}
  Let $T > 0$. Then, there exists a $C(T) > 0$ such that for all $\varepsilon \in (0,1)$, the following holds true
  \begin{align}
    \label{tim_sq1}
    \int_0^T \left\|\partial_t \sqrt{1 + c_{1\varepsilon}(\cdot, t)}\right\|_{(W^{3, 2}(\Omega))^\ast}dt \leq C(T).
  \end{align}

\end{lemma}

\begin{proof}
  Fix $T > 0$. Using $\partial_t \sqrt{1 + c_{1\varepsilon}} = \frac{\partial_t c_{1\varepsilon}}{2\sqrt{1 + c_{1\varepsilon}}}$ in $\Omega \times (0, \infty)$ and an arbitrary test function $\varphi \in C^\infty(\overline{\Omega})$, we derive from the first equation of \eqref{app_prob}
  \begin{align}
    \label{tim_c1}
    \nonumber & 2 \int_\Omega \partial_t \sqrt{1 + c_{1\varepsilon}} \cdot \varphi = \int_\Omega \left\{- \mathbb{D}_\varepsilon \nabla c_{1\varepsilon} - c_{1\varepsilon} \nabla \cdot \mathbb{D}_\varepsilon + \chi c_{1\varepsilon} \mathbb{D}_\varepsilon \nabla \tau_\varepsilon\right\} \cdot \nabla \frac{\varphi}{\sqrt{1 + c_{1\varepsilon}}}   \\[5pt]
    \nonumber & + \int_\Omega \left\{- \alpha_1  \tfrac{c_{1\varepsilon}}{1 + c_{1\varepsilon}}  + \alpha_2 \tfrac{c_{2\varepsilon}}{1 + c_{2\varepsilon}} + \beta c_{1\varepsilon}(1 - c_{1\varepsilon} - c_{2\varepsilon} - \tau_\varepsilon)\right\} \frac{\varphi}{\sqrt{1 + c_{1\varepsilon}}} \\[5pt]
    \nonumber &=  \frac{1}{2} \int_\Omega \frac{\nabla c_{1 \varepsilon} \cdot \mathbb{D}_\varepsilon \nabla c_{1 \varepsilon}}{\sqrt{1 + c_{1\varepsilon}}^3} \varphi - \int_\Omega \frac{\nabla c_{1 \varepsilon} \cdot \mathbb{D}_\varepsilon \nabla \varphi}{\sqrt{1 + c_{1\varepsilon}}} + \frac{1}{2} \int_\Omega \frac{c_{1\varepsilon} ((\nabla \cdot \mathbb{D}_\varepsilon) \cdot \nabla c_{1\varepsilon})}{\sqrt{1 + c_{1\varepsilon}}^3}\varphi - \int_\Omega \frac{c_{1\varepsilon}(\nabla \cdot \mathbb{D}_\varepsilon) \cdot \nabla \varphi}{\sqrt{1 + c_{1\varepsilon}}} \\[5pt]
    \nonumber & - \frac{\chi}{2} \int_\Omega \frac{c_{1\varepsilon} \nabla c_{1\varepsilon} \cdot \mathbb{D}_\varepsilon \nabla \tau_\varepsilon}{\sqrt{1 + c_{1\varepsilon}}^3}\varphi + \chi \int_\Omega \frac{c_{1\varepsilon} \nabla \tau_\varepsilon \cdot \mathbb{D}_\varepsilon \nabla \varphi}{\sqrt{1 + c_{1\varepsilon}}} - \alpha_1  \int_\Omega \frac{c_{1\varepsilon}}{\sqrt{1 + c_{1\varepsilon}}^3} \varphi + \alpha_2 \int_\Omega \frac{c_{2\varepsilon}}{1 + c_{2\varepsilon}}  \frac{\varphi}{\sqrt{1 + c_{1\varepsilon}}} \\[5pt]
    & + \beta \int_\Omega c_{1\varepsilon}(1 - c_{1\varepsilon} - c_{2\varepsilon} - \tau_\varepsilon) \frac{\varphi}{\sqrt{1 + c_{1\varepsilon}}}
  \end{align}
 for all $t \in (0, T)$ and all $\varepsilon \in (0, 1)$. We will now make several straightforward estimations involving the inequalities of Cauchy and Young, the boundedness of functions $\tfrac{c_{i \varepsilon}}{1 + c_{i \varepsilon}} < 1$ and $\tfrac{1}{1 + c_{i \varepsilon}} < 1$ for $c_{i\varepsilon} \geq 0,\ i \in \{1, 2\}$, \eqref{epassumption2}, \eqref{weak_assup1}, \eqref{bound_comp1}-\eqref{bound_comp4} and \eqref{nab_tau_bound11}:
 \begin{align}
 \label{tim_c3}
   \left|\frac{1}{2} \int_\Omega \frac{\nabla c_{1 \varepsilon} \cdot \mathbb{D}_\varepsilon \nabla c_{1 \varepsilon}}{\sqrt{1 + c_{1\varepsilon}}^3} \varphi\right| \leq \left(\frac{1}{2} \int_\Omega \frac{\nabla c_{1 \varepsilon} \cdot \mathbb{D}_\varepsilon \nabla c_{1 \varepsilon}}{c_{1\varepsilon}} \right) \cdot \|\varphi\|_{L^\infty(\Omega)}
 \end{align}
 and
\begin{align}
  \label{tim_c41}
  \nonumber \left|\int_\Omega \frac{\nabla c_{1 \varepsilon} \cdot \mathbb{D}_\varepsilon \nabla \varphi}{\sqrt{1 + c_{1\varepsilon}}}\right| & \leq \left(\int_\Omega \nabla \varphi \cdot \mathbb{D}_\varepsilon \nabla \varphi\right)^{\frac{1}{2}} \cdot \left(\int_\Omega \frac{\nabla c_{1 \varepsilon} \cdot \mathbb{D}_\varepsilon \nabla c_{1 \varepsilon}}{1 + c_{1\varepsilon}}\right)^{\frac{1}{2}}   \\[5pt]
  \nonumber & \leq \frac{1}{2} \int_\Omega \nabla \varphi \cdot \mathbb{D}_\varepsilon \nabla \varphi + \frac{1}{2} \int_\Omega \frac{\nabla c_{1 \varepsilon} \cdot \mathbb{D}_\varepsilon \nabla c_{1 \varepsilon}}{1 + c_{1\varepsilon}}\\[5pt]
  \nonumber & \leq \frac{1}{2} \int_\Omega \nabla \varphi \cdot \mathbb{D}_\varepsilon \nabla \varphi + \frac{1}{2} \int_\Omega \frac{\nabla c_{1 \varepsilon} \cdot \mathbb{D}_\varepsilon \nabla c_{1 \varepsilon}}{c_{1\varepsilon}} \\[5pt]
  & \leq \frac{1}{2} \|\nabla \varphi\|^2_{L^\infty(\Omega)} \|\mathbb{D}_\varepsilon\|_{L^\infty(\Omega)} |\Omega| + \frac{1}{2} \int_\Omega \frac{\nabla c_{1 \varepsilon} \cdot \mathbb{D}_\varepsilon \nabla c_{1 \varepsilon}}{c_{1\varepsilon}}
\end{align}
and
 \begin{align}
   \label{tim_c2}
   \nonumber & \left|\frac{1}{2}\int_\Omega \frac{c_{1\varepsilon} ((\nabla \cdot \mathbb{D}_\varepsilon) \cdot \nabla c_{1\varepsilon})}{\sqrt{1 + c_{1\varepsilon}}^3}\varphi\right| = \left|\int_\Omega\frac{c_{1\varepsilon}}{1+c_{1\varepsilon}}(\nabla\cdot\mathbb{D}_\varepsilon)\cdot\nabla\sqrt{1+c_{1\varepsilon}}\,\varphi\right| \\[5pt]
   \nonumber& = \left|\int_\Omega\frac{c_{1\varepsilon}}{1+c_{1\varepsilon}}(\nabla\cdot\mathbb{D}_\varepsilon)\cdot\nabla\left(\sqrt{1+c_{1\varepsilon}}\varphi\right) - \int_\Omega
\frac{c_{1\varepsilon}}{\sqrt{1+c_{1\varepsilon}}}(\nabla\cdot\mathbb{D}_\varepsilon)\cdot\nabla\varphi\right| \\[5pt]
 \nonumber & \leq \int_\Omega \left|\frac{c_{1\varepsilon}}{1+c_{1\varepsilon}}(\nabla\cdot\mathbb{D}_\varepsilon)\cdot
\nabla\left(\sqrt{1+c_{1\varepsilon}}\varphi\right)\right| +
\int_\Omega\frac{c_{1\varepsilon}}{\sqrt{1+c_{1\varepsilon}}}\left|(\nabla\cdot\mathbb{D}_\varepsilon)\cdot\nabla\varphi
\right|\\[5pt]
\nonumber &  \leq \int_\Omega|\left(\nabla\cdot\mathbb{D}_\varepsilon)\cdot\nabla\left(\sqrt{1+c_{1\varepsilon}}\varphi\right)
\right| +\|\nabla\cdot\mathbb{D}_\varepsilon\|_{L^\infty(\Omega)}\|\nabla\varphi\|_{L^\infty(\Omega)} \int_\Omega
\frac{c_{1\varepsilon}}{\sqrt{1+c_{1\varepsilon}}}\\[5pt]
\nonumber & \leq 
B\int_\Omega\left(\nabla\left(\sqrt{1+c_{1\varepsilon}}\varphi\right)\cdot\mathbb{D}_\varepsilon\nabla\left(\sqrt{1+c_{1\varepsilon}}\varphi\right)
\right)^\lambda+B 
%{\cb \frac{|\Omega|}{2D_1}\int_\Omega \nabla\left(\sqrt{1+c_{1\varepsilon}}\varphi\right)\cdot\mathbb{D}_\varepsilon\nabla\left(\sqrt{1+c_{1\varepsilon}}\varphi\right)+\frac{M^2}{2}}
+\|\nabla\cdot\mathbb{D}_\varepsilon\|_{L^\infty(\Omega)}\|\nabla\varphi\|_{L^\infty(\Omega)}
|\Omega|^{1/2}\left(\int_\Omega c_{1\varepsilon}\right)^{1/2}\\[5pt]
\nonumber &  \leq B\int_\Omega\nabla\left(\sqrt{1+c_{1\varepsilon}}\varphi\right)\cdot\mathbb{D}_\varepsilon\nabla\left(\sqrt{1+c_{1\varepsilon}}\varphi\right)+B(|\Omega|+1) 
%{\cb \frac{|\Omega|}{2D_1}\int_\Omega \nabla\left(\sqrt{1+c_{1\varepsilon}}\varphi\right)\cdot\mathbb{D}_\varepsilon\nabla\left(\sqrt{1+c_{1\varepsilon}}\varphi\right)+\frac{M^2}{2}}
+\|\nabla\cdot\mathbb{D}_\varepsilon\|_{L^\infty(\Omega)}\|\nabla\varphi\|_{L^\infty(\Omega)}|\Omega|^{1/2}\widehat{M}_1^{1/2} \\[5pt]
\nonumber & \leq \frac{B}{2}\int_\Omega|\varphi|^2\frac{\nabla c_{1\varepsilon}\cdot\mathbb{D}_\varepsilon\nabla c_{1\varepsilon}}{1+c_{1\varepsilon}}
+2B\int_\Omega(1+c_{1\varepsilon})\nabla\varphi\cdot\mathbb{D}_\varepsilon\nabla\varphi+B(|\Omega|+1)\\[5pt] 
\nonumber&+\|\nabla\cdot\mathbb{D}_\varepsilon\|_{L^\infty(\Omega)}\|\nabla\varphi\|_{L^\infty(\Omega)}|\Omega|^{1/2}\widehat{M}_1^{1/2}\\[5pt]
\nonumber & \leq \frac{B}{2}\int_\Omega|\varphi|^2\frac{\nabla c_{1\varepsilon}\cdot\mathbb{D}_\varepsilon\nabla c_{1\varepsilon}}{c_{1\varepsilon}}+2B\|\mathbb{D}_\varepsilon\|_{L^\infty(\Omega)}\|\nabla\varphi\|_{L^\infty(\Omega)}^2
\left(|\Omega|+\int_\Omega c_{1\varepsilon}\right)+B(|\Omega|+1) \\[5pt]
\nonumber & +\|\nabla\cdot\mathbb{D}_\varepsilon\|_{L^\infty(\Omega)}\|\nabla\varphi\|_{L^\infty(\Omega)}
|\Omega|^{1/2}\widehat{M}_1^{1/2}\\[5pt]
& \leq\frac{B\|\varphi\|_{L^\infty(\Omega)}^2}{2}\int_\Omega\frac{\nabla c_{1\varepsilon}\cdot\mathbb{D}_\varepsilon
\nabla c_{1\varepsilon}}{c_{1\varepsilon}}+C_1
 \end{align}
and
\begin{align}
  \label{tim_c5}
  \nonumber \left|\int_\Omega \frac{c_{1\varepsilon}(\nabla \cdot \mathbb{D}_\varepsilon) \cdot \nabla \varphi}{\sqrt{1 + c_{1\varepsilon}}}\right| & \leq \left\{\frac{1}{2} \int_\Omega \frac{c_{1\varepsilon}^2 }{1 + c_{1\varepsilon}} + \frac{1}{2} \int_\Omega |\nabla \cdot \mathbb{D}_\varepsilon|^2\right\}\cdot \|\nabla \varphi\|_{L^\infty(\Omega)}  \\[5pt]
  \nonumber & \leq \left\{\frac{1}{2} \int_\Omega c_{1\varepsilon} + \frac{1}{2} \int_\Omega |\nabla \cdot \mathbb{D}_\varepsilon|^2\right\}\cdot \|\nabla \varphi\|_{L^\infty(\Omega)} \\[5pt]
  & \leq (\tfrac{\widehat{M}_1}{2} + \tfrac{M^2}{2}) \cdot \|\nabla \varphi\|_{L^\infty(\Omega)}
\end{align}
and
\begin{align}
  \label{time_c6}
  \nonumber \left|\frac{\chi}{2} \int_\Omega \frac{c_{1\varepsilon} \nabla c_{1\varepsilon} \cdot \mathbb{D}_\varepsilon \nabla \tau_\varepsilon}{\sqrt{1 + c_{1\varepsilon}}^3}\varphi\right| & \leq \frac{\chi}{2} \left\{\int_\Omega \frac{\nabla c_{1\varepsilon} \cdot \mathbb{D}_\varepsilon \nabla c_{1\varepsilon}}{1 + c_{1\varepsilon}} \right\}^{\frac{1}{2}} \cdot \left\{\int_\Omega \left(\frac{c_{1\varepsilon}}{1 + c_{1\varepsilon}}\right)^2 \nabla \tau_\varepsilon \cdot \mathbb{D}_\varepsilon \nabla \tau_\varepsilon\right\}^{\frac{1}{2}}\|\varphi\|_{L^\infty(\Omega)}  \\[5pt]
  \nonumber & \leq \frac{\chi}{4} \left\{\int_\Omega \frac{\nabla c_{1\varepsilon} \cdot \mathbb{D}_\varepsilon \nabla c_{1\varepsilon}}{1 + c_{1\varepsilon}} + \int_\Omega \left(\frac{c_{1\varepsilon}}{1 + c_{1\varepsilon}}\right)^2 \nabla \tau_\varepsilon \cdot \mathbb{D}_\varepsilon \nabla \tau_\varepsilon\right\} \|\varphi\|_{L^\infty(\Omega)} \\[5pt]
  & \nonumber \leq \frac{\chi}{4} \left\{\int_\Omega \frac{\nabla c_{1\varepsilon} \cdot \mathbb{D}_\varepsilon \nabla c_{1\varepsilon}}{c_{1\varepsilon}} + \int_\Omega \nabla \tau_\varepsilon \cdot \mathbb{D}_\varepsilon \nabla \tau_\varepsilon\right\}\|\varphi\|_{L^\infty(\Omega)}\\[5pt]
  & \leq \frac{\chi}{4} \left\{\int_\Omega \frac{\nabla c_{1\varepsilon} \cdot \mathbb{D}_\varepsilon \nabla c_{1\varepsilon}}{c_{1\varepsilon}} +  M_{2\tau} \right\}\|\varphi\|_{L^\infty(\Omega)}
\end{align}
and
\begin{align}
  \label{tim_c7}
  \nonumber \left|\chi \int_\Omega \frac{c_{1\varepsilon} \nabla \tau_\varepsilon \cdot \mathbb{D}_\varepsilon \nabla \varphi}{\sqrt{1 + c_{1\varepsilon}}}\right| & \leq \chi \left\{\int_\Omega \frac{c_{1\varepsilon}^2}{1 + c_{1\varepsilon}}\right\}^{\frac{1}{2}} \cdot \left\{\int_\Omega \nabla \tau_\varepsilon \cdot \mathbb{D}_\varepsilon \nabla \tau_{\varepsilon}\right\}^{\frac{1}{2}} \|\nabla \varphi\|_{L^\infty(\Omega)} \\[5pt]
  \nonumber & \leq \frac{\chi}{2} \left\{\int_\Omega c_{1\varepsilon} + \int_\Omega \nabla \tau_\varepsilon \cdot \mathbb{D}_\varepsilon \nabla \tau_{\varepsilon}\right\}\|\nabla \varphi\|_{L^\infty(\Omega)} \\[5pt]
  & \leq \frac{\chi}{2} \left\{\widehat{M}_1 +  M_{2\tau}\right\}\|\nabla \varphi\|_{L^\infty(\Omega)}
\end{align}
as well as
\begin{align}
  \label{tim_c8}
  & \left|\alpha_1  \int_\Omega \frac{c_{1\varepsilon}}{\sqrt{1 + c_{1\varepsilon}}^3} \varphi\right| \leq \frac{\alpha_1}{2} \left\{\int_\Omega c_{1\varepsilon}^2 + |\Omega|\right\}\|\varphi\|_{L^\infty(\Omega)}, \\[5pt]
  \label{tim_c21}
  & \left|\alpha_2 \int_\Omega \frac{c_{2\varepsilon}}{1 + c_{2\varepsilon}}  \frac{\varphi}{\sqrt{1 + c_{1\varepsilon}}} \right| \leq \alpha_2 \|\varphi\|_{L^\infty(\Omega)},
\end{align}
and
\begin{align}
  \label{tim_c22}
  & \left| \beta \int_\Omega c_{1\varepsilon}(1 - c_{1\varepsilon} - c_{2\varepsilon} - \tau_\varepsilon) \frac{\varphi}{\sqrt{1 + c_{1\varepsilon}}}\right| \leq \beta \left\{(\widehat{M}_{c_2} + \widehat{M}_\tau + 1)\widehat{M}_1 + \int_\Omega  c_{1\varepsilon}^2\right\} \|\varphi\|_{L^\infty(\Omega)}
\end{align}
Inserting \eqref{tim_c3}-\eqref{tim_c22} in \eqref{tim_c1} results in
\begin{align}
\label{tim_fi1}
  2\left| \int_\Omega \partial_t \sqrt{1 + c_{1\varepsilon}} \cdot \varphi \right| \leq f_\varepsilon(t) \cdot \|\varphi\|_{W^{1, \infty}(\Omega)} \quad \text{for all}~ \varphi \in C^\infty(\overline{\Omega})
\end{align}
where,
\begin{align*}
  f_\varepsilon(t) := C_2(T)\left\{ \int_\Omega \frac{\nabla c_{1\varepsilon} \cdot \mathbb{D}_\varepsilon \nabla c_{1\varepsilon}}{c_{1\varepsilon}} + \int_\Omega c_{1\varepsilon}^2 + 1\right\}
\end{align*}
for all $t \in (0, T)$ and all $\varepsilon \in (0, 1)$ with an adequate $C_2(T) > 0$. As $W^{3, 2}(\Omega)$ is continuously embedded into $W^{1, \infty}(\Omega)$, \eqref{tim_fi1} ensures that there exists a constant $C_3 > 0$ such that
\begin{align*}
  \left\|\partial_t \sqrt{1 + c_{1\varepsilon}(\cdot, t)}\right\|_{(W^{3, 2}(\Omega))^\ast} = \sup_{\substack{\varphi \in C^\infty(\Omega),\\ \|\varphi\|_{W^{3, 2}(\Omega)} \leq 1}} \left|\int_\Omega \partial_t \sqrt{1 + c_{1\varepsilon}(\cdot, t)} \cdot \varphi\right| \leq C_3 f_\varepsilon(t)
\end{align*}
for all $t \in (0, T)$ and all $\varepsilon \in (0, 1)$. A combination of the estmates in \eqref{bound_comp4} and \eqref{weentropy0} then yields a $C_4(T) > 0$ such that
\begin{align*}
  \int_0^T f_\varepsilon(t) \leq C_4(T) \quad \text{for all}~ \varepsilon \in (0, 1),
\end{align*}
which establishes \eqref{tim_sq1}. %{\cred Looks like you apply \eqref{bound_comp4} rather than \eqref{bound_comp1}. And there you deal with $T_{\max,\varepsilon}$; how does the latter relate to your $T$?} {\color{blue} Please see the modified lemma. }

\end{proof}

%%%%%%%%%%%%%%%%%%%%%%%%%%%%%%%%%%%%%%%%%%%%%%%%%%%%%%%%%%%%%%%%%%%%%%%%%%%%%%%%%%%%%%%%%%%%%%%%%%%%%%%%%%%%%%%%%%%%%%%%%%%%%%%%%%%%%%%%%%%%%%%%%%%%%%%%%%%%%%%%%%%%%%%%%%%%%%%%%%%%%%%%%%%%%%%%%%%%%%%%%%%%%%%%%%%%%%%%%%%%%%%%%
%%%%%%%%%%%%%%%%%%%%%%%%%%%%%%%%%%%%%%%%%%%%%%%%%%%%%%%%%%%%%%%%%%%%%%%%%%%%%%%%%%%%%%%%%%%%%%%%%%%%%%%%%%%%%%%%%%%%%%%%%%%%%%%%%%%%%%%%%%%%%%%%%%%%%%%%%%%%%%%%%%%%%%%%%%%%%%%%%%%%%%%%%%%%%%%%%%%%%%%%%%%%%%%%%%%%%%%%%%%%%%%%%

\subsection{Construction of weak solutions}

We are now able to construct the solution candidates by employing various compact embedding arguments, allowing us to obtain them as the limit of the sequence of approximate solutions, thereby establishing the proof of our main result for this section. %, Theorem \ref{main_th_weak_sol}.

\begin{proof}[Proof of Theorem \ref{main_th_weak_sol}]

  We first claim that there exist a null sequence $(\varepsilon_j)_{j \in \mathbb{N}} \subseteq (0, 1)$ and a.e. non-negative functions
  \begin{align*}
  \begin{cases}
    c_1 \in L^2_{\textup{loc}}(\overline{\Omega} \times (0, \infty)) \quad \text{with} \quad \sqrt{1 + c_1} \in L^{\frac{4}{3}}_{\textup{loc}}([0, \infty); W_{\mathbb{D}}^{1, \frac{4}{3}}(\Omega)), \\[5pt]
    c_2 \in L^\infty(\Omega \times (0, \infty)) \cap L^2_{\textup{loc}}([0, \infty); W^{1, 2}(\Omega)),                                                                                          \\[5pt]
    \tau \in L^\infty(\Omega \times (0, \infty)) \cap L^2_{\textup{loc}}([0, \infty); W_{\mathbb{D}}^{1, 2}(\Omega))
  \end{cases}
  \end{align*}
  such that
  \begin{align}
    \label{weak_conv1} & c_{1\varepsilon} \to c_1     \quad \text{in}~ L^2_{\textup{loc}}(\overline{\Omega} \times [0, \infty)) ~\text{and a.e. in}~ \Omega \times [0, \infty),                                \\[5pt]
    \label{weak_conv2} & c_{1\varepsilon}^2 \to c_1^2 \quad \text{in}~ L^1_{\textup{loc}}(\overline{\Omega} \times [0, \infty)) ~\text{and a.e. in}~ \Omega \times [0, \infty),                                \\[5pt]
    \label{weak_conv3} & \sqrt{1 + c_{1\varepsilon}} \rightharpoonup \sqrt{1 + c_1} \quad \text{in} ~ L^{\frac{4}{3}}_{\textup{loc}}([0, \infty); W_{\mathbb{D}}^{1, \frac{4}{3}}(\Omega)),              \\[5pt]
    \label{weak_conv4} & c_{2\varepsilon} \to c_2 \quad \text{in}~ L^p_{\textup{loc}}(\overline{\Omega} \times [0, \infty)) ~\text{for all}~ p \in [1, \infty) ~\text{and a.e. in}~ \Omega \times [0, \infty),                \\[5pt]
    \label{weak_conv5} & c_{2\varepsilon} \rightharpoonup c_2 \quad \text{in}~ L^2_{\textup{loc}}([0, \infty); W^{1, 2}(\Omega)),                                                                         \\[5pt]
    \label{weak_conv6} & \tau_{\varepsilon} \to \tau \quad ~\text{in}~ L^p_{\textup{loc}}(\overline{\Omega} \times [0, \infty)) ~\text{for all}~ p \in [1, \infty) ~\text{and a.e. in}~ \Omega \times [0, \infty),             \\[5pt]
    \label{weak_conv7} & \tau_{\varepsilon} \rightharpoonup \tau \quad ~\text{in}~ L^2_{\textup{loc}}([0, \infty); W^{1, 2}_{\mathbb{D}}(\Omega))
  \end{align}
 as $\varepsilon = \varepsilon_j \searrow 0$. The validity of \eqref{weak_conv1}-\eqref{weak_conv7} can be established in the following manner:

 Collecting the estimates provided by Lemmas \ref{lemma_epsilon1} and \ref{lemma_sqrt1} we conclude that $(\sqrt{1 + c_{1\varepsilon}})_{\varepsilon \in (0, 1)}$ is bounded in $L^{\frac43}_{\textup{loc}}
([0,\infty);W_{\mathbb{D}}^{1,\frac43}(\Omega))$ and $(\tau_\varepsilon)_{\varepsilon \in (0, 1)}$ is bounded in $L^2_{\textup{loc}}([0, \infty); W_{\mathbb{D}}^{1,2}(\Omega))$, while $(c_{2\varepsilon})_{\varepsilon \in (0, 1)}$ is bounded in $L^2_{\textup{loc}}([0, \infty); W^{1,2}(\Omega))$. Furthermore, by Lemmas \ref{time_det_c2tau} and \ref{time_sqrt0}, the families $(\partial_t\sqrt{1 + c_{1\varepsilon}})_{\varepsilon \in (0, 1)}, (\partial_t c_{2\varepsilon})_{\varepsilon \in (0, 1)}$ and $(\partial_t\tau_{\varepsilon})_{\varepsilon \in (0, 1)}$ are bounded in $L^1_{\textup{loc}}([0, \infty); (W^{3, 2}(\Omega))^\ast)$. With the above information at hand, we can now apply the
Aubin--Lions lemma \cite[p.~215]{temam2024navier} to the above
mentioned families using the corresponding triples of embedded
spaces
\begin{align*}
W_{\mathbb{D}}^{1,\frac43}(\Omega) \subseteq L^1(\Omega) \subseteq (W^{3,2}(\Omega))^\ast,
\end{align*}
for the first solution component,
\begin{align*}
W^{1,2}(\Omega) \subseteq L^1(\Omega) \subseteq (W^{3,2}(\Omega))^\ast,
\end{align*}
for the second solution component, and
\begin{align*}
W_{\mathbb{D}}^{1,2}(\Omega) \subseteq L^1(\Omega) \subseteq (W^{3,2}(\Omega))^\ast,
\end{align*}
for the third solution component. Then applying the Aubin-Lions lemma countably infinitely many times
on the time intervals $[0,T]$, $T\in\mathbb N$, combined with an
extension and diagonal sequence argument, we can construct a null
sequence $(\varepsilon_j)_{j\in\mathbb N}\subseteq(0,1)$ and
functions $\sqrt{1+\widetilde c_1}$, $c_2$ and $\tau$ such that
\begin{align*}
\sqrt{1+c_{1\varepsilon}} &\to\sqrt{1+\widetilde c_1} \quad\text{in}~ L_{\textup{loc}}^{\frac43}([0,\infty);L^1(\Omega)),\\[5pt]
c_{2\varepsilon} & \to c_2
\quad\text{and}\quad \tau_\varepsilon\to\tau \quad\text{in}~
L_{\textup{loc}}^2([0,\infty);L^1(\Omega))
\end{align*}
as $\varepsilon = \varepsilon_j \searrow 0$.

 Moreover, the sequence above can be constructed in such a way that it ensures $\sqrt{1 + c_{1\varepsilon}} \to \sqrt{1 + \tilde{c}_1}$, $c_{2\varepsilon} \to c_2$ and $\tau_\varepsilon \to \tau$ pointwise almost everywhere as $\varepsilon = \varepsilon_j \searrow 0$ by passing to a further subsequence if necessary.

 Furthermore, by \eqref{weak_assup2} and Lemma \ref{lbound_comp}, the families $(c_{2\varepsilon})_{\varepsilon \in (0,1)}$ and $(\tau_{\varepsilon})_{\varepsilon \in (0,1)}$ are bounded in $L^\infty(\Omega \times (0, \infty))$. Consequently, the convergence properties established above directly imply \eqref{weak_conv4} and \eqref{weak_conv6}, ensuring that $c_2$ and $\tau$ remain nonnegative a.e and that $c_2, \tau \in L^\infty(\Omega \times (0, \infty))$.

 Also, the pointwise convergence everywhere for the constructed sequence ensures that
 \begin{align*}
   c_{1\varepsilon} \to \tilde{c}_1 \quad {\textup{and}} \quad c_{1\varepsilon}^2 \to \tilde{c}_{1}^2 \quad \text{a.e. pointwise}%\quad \text{\cred shouldn't the limits be $\tilde c_1, \tilde c_1^2$?} \text{{\color{blue} Yes, done.}}
 \end{align*}
 as $\varepsilon = \varepsilon_j \searrow 0$. This establishes the nonnegativity of $\tilde c_1$ as well. From Lemma \ref{entropy_weak:lemma} we can have a $\widehat{C}_1(T) > 0$ such that 
 \begin{align*}
     \int_0^T\int_\Omega c_{1\varepsilon}^2(\ln2+c_{1\varepsilon}) \leq \widehat{C}_1(T). 
 \end{align*} 
 If $c_{1\varepsilon} \geq 1$ we have 
 \begin{align*}
     c_{1\varepsilon}^2 |\ln c_{1\varepsilon}| = c_{1\varepsilon}^2 \ln c_{1\varepsilon} \leq c_{1\varepsilon}^2 \ln (2 + c_{1\varepsilon}).
 \end{align*}
 On the other hand , since $\displaystyle \sup_{0<s\leq 1}s^2|\ln s| = \tfrac{1}{2e}$, we obtain $c_{1\varepsilon}^2|\ln c_{1\varepsilon}|\leq \frac{1}{2e}$ whenever $0<c_{1\varepsilon}<1$. Consequently,
\begin{align*}
    \nonumber \int_0^T\int_\Omega c_{1\varepsilon}^2|\ln c_{1\varepsilon}| & \leq \int_0^T\int_\Omega c_{1\varepsilon}^2(\ln2+c_{1\varepsilon}) +\frac{|\Omega|T}{2e}  \\[5pt]
    &\leq \widehat{C}_1(T)+\frac{|\Omega|T}{2e} =: C(T).
\end{align*}
Hence, from above we conclude that there exists a $C(T) > 0$ such that
\begin{align*}
    \int_0^T\int_\Omega
    c_{1\varepsilon}^2|\ln c_{1\varepsilon}| \leq C(T)
\end{align*}
for all $\varepsilon\in(0,1)$.
%{\cred How do you get this from there? Unless it is not the same $C(T)$ I do not see it holding for all $c_{1\varepsilon}\ge 0$.} {\color{blue}Please see the above proof.} 

We can now use Vitali's theorem together with the de la Vall\'ee Poussin criterion for uniform integrability to obtain
\begin{align*}
c_{1\varepsilon}^2\to c_1^2
\quad\text{in}~
L^1(\Omega\times(0,T)).
\end{align*}
In particular,
\begin{align*}
\|c_{1\varepsilon}-c_1\|_{L^2(\Omega\times(0,T))}^2
&\leq
\int_0^T\int_\Omega
|c_{1\varepsilon}^2-c_1^2|
\to0,
\end{align*}
which establishes \eqref{weak_conv1} and \eqref{weak_conv2}.

The convergence properties in \eqref{weak_conv3}, \eqref{weak_conv5}, and \eqref{weak_conv7} then follow from a standard extraction procedure since the corresponding family of functions is bounded in the respective spaces as established in Lemmas \ref{lemma_epsilon1} and \ref{lemma_sqrt1}.

%%%%%%%%%%%%%%%%%%%%%%%%%%%%%%%%%%%%%%%%%%%%%%%%%%%%%%%%%%%%%%%%%%%%%%%%%%%%%%%%%%%%%%%%%%%%%%%%%%%%%%%%%%%%%%%%%%%%%%%%%%%%%%%%%%%%%%%%%%%%%%%%%%%%%%%%%%%%%%%%%%%%%%%%%%%%%%%%%%%%%%%%%%%%%%%%%%%%%%%%%%%%%%%%%%%%%%%%%%%%%%%%%
%%%%%%%%%%%%%%%%%%%%%%%%%%%%%%%%%%%%%%%%%%%%%%%%%%%%%%%%%%%%%%%%%%%%%%%%%%%%%%%%%%%%%%%%%%%%%%%%%%%%%%%%%%%%%%%%%%%%%%%%%%%%%%%%%%%%%%%%%%%%%%%%%%%%%%%%%%%%%%%%%%%%%%%%%%%%%%%%%%%%%%%%%%%%%%%%%%%%%%%%%%%%%%%%%%%%%%%%%%%%%%%%%
As a final step, we will derive the following convergence results:
  \begin{align}
    \label{c1_convergence1}   & \int_0^T \int_\Omega \sqrt{1 + c_{1\varepsilon}} \nabla \sqrt{1 + c_{1\varepsilon}} \cdot \mathbb{D}_\varepsilon \nabla \varphi \to \int_0^T \int_\Omega \sqrt{1 + c_{1}} \nabla \sqrt{1 + c_{1}} \cdot \mathbb{D} \nabla \varphi  ~~\text{as}~~ \varepsilon = \varepsilon_j \searrow 0
    \end{align}
    and
    \begin{align}
    \label{tau_convergence1}  & \int_0^T \int_\Omega c_{1\varepsilon} \nabla \tau_\varepsilon \cdot \mathbb{D}_\varepsilon \nabla \varphi \to \int_0^T \int_\Omega c_1 \nabla \tau \cdot \mathbb{D} \nabla \varphi ~~\text{as}~~ \varepsilon = \varepsilon_j \searrow 0
  \end{align}
  for all $\varphi \in C_c^\infty(\overline{\Omega} \times [0, \infty))$.

 We first fix $\varphi \in C_c^\infty(\overline{\Omega} \times [0, \infty))$ and $T > 0$ such that $\textup{supp} ~ \varphi \subseteq \overline{\Omega} \times [0, T)$. Then from {the estimates in the proofs of} Lemmas \ref{lemma_epsilon1} and \ref{lemma_sqrt1} we can have a $C_1 > 0$ such that
\begin{align}
\label{convergence22}
  \int_0^T \int_\Omega (\nabla \sqrt{1 + c_{1\varepsilon}} \cdot \mathbb{D}_\varepsilon \nabla \sqrt{1 + c_{1\varepsilon}})^{\frac{2}{3}} \leq C_1 \quad \textup{and} \quad \int_0^T \int_\Omega \nabla \tau_\varepsilon \cdot \mathbb{D}_\varepsilon \nabla \tau_\varepsilon \leq C_1
\end{align}
for all $\varepsilon \in (0, 1)$. Thereafter, using \eqref{epassumption3} and \eqref{convergence22},
\begin{align}
\label{convergence23}
  \nonumber \left\|\nabla\sqrt{1+c_{1\varepsilon}}
\right\|_{L^{\frac43}(\Omega\times(0,T))} & = \left(\int_0^T\int_\Omega|\nabla\sqrt{1+c_{1\varepsilon}}|^{\frac43}
   \right)^{\frac34}     \\[5pt]
 & \leq\left(\frac{1}{\varepsilon^{\frac{2}{3}}}\int_0^T \int_\Omega (\nabla \sqrt{1 + c_{1\varepsilon}} \cdot \mathbb{D}_\varepsilon \nabla \sqrt{1 + c_{1\varepsilon}})^{\frac{2}{3}}\right)^{\frac{3}{4}} \leq  C_1^{\frac{3}{4}}\varepsilon^{-\frac{1}{2}}
\end{align}
and
\begin{align}
\label{convergence24}
  \|\nabla \tau_\varepsilon\|_{L^2(\Omega \times (0, T))} & = \left(\int_0^T \int_\Omega \nabla \tau_\varepsilon \cdot \nabla \tau_\varepsilon\right)^{\frac{1}{2}}  \leq \left(\frac{1}{\varepsilon}\int_0^T \int_\Omega \nabla \tau_\varepsilon \cdot \mathbb{D}_\varepsilon \nabla \tau_\varepsilon\right)^{\frac{1}{2}} \leq  \sqrt{C_1}\varepsilon^{-\frac{1}{2}}
\end{align}
for all $\varepsilon \in (0, 1)$. Furthermore, by \eqref{bound_comp4} and \eqref{weak_conv2}, there exists a $C_2>0$ such that
\begin{align*}
\left\|\sqrt{1+c_{1\varepsilon}}\right\|_{L^4(\Omega\times(0,T))}^4
&=\int_0^T\int_\Omega(1+c_{1\varepsilon})^2
\leq 2|\Omega|T + 2\int_0^T\int_\Omega c_{1\varepsilon}^2 =: C_2
\end{align*}
for all $\varepsilon\in(0,1)$. From \eqref{epassumption3} and \eqref{convergence23} we can then have
\begin{align*}
&2\left|\int_0^T\int_\Omega\sqrt{1+c_{1\varepsilon}}\nabla\sqrt{1+c_{1\varepsilon}}\cdot\mathbb{D}_\varepsilon\nabla\varphi-\int_0^T\int_\Omega\sqrt{1+c_1}\nabla\sqrt{1+c_1}\cdot\mathbb{D}\nabla\varphi\right|\\[5pt]
&\leq2\left|\int_0^T\int_\Omega\sqrt{1+c_{1\varepsilon}}\nabla\sqrt{1+c_{1\varepsilon}}\cdot(\mathbb{D}_\varepsilon-\mathbb{D})\nabla\varphi\right|\\[5pt]
&+2\left|\int_0^T\int_\Omega\left(\sqrt{1+c_{1\varepsilon}}\nabla\sqrt{1+c_{1\varepsilon}} -\sqrt{1+c_1}\nabla\sqrt{1+c_1}\right)\cdot\mathbb{D}\nabla\varphi\right|
\\[5pt]
&\leq2\|\mathbb{D}_\varepsilon-\mathbb{D}\|_{L^\infty(\Omega)}\|\nabla\varphi\|_{L^\infty(\Omega\times(0,T))}
\int_0^T\int_\Omega\sqrt{1+c_{1\varepsilon}}\left|\nabla\sqrt{1+c_{1\varepsilon}}\right|\\[5pt]
&\quad+2\left|\int_0^T\int_\Omega\left(\sqrt{1+c_{1\varepsilon}}\nabla\sqrt{1+c_{1\varepsilon}}-
\sqrt{1+c_1}\nabla\sqrt{1+c_1}\right)\cdot\mathbb{D}\nabla\varphi\right|\\[5pt]
&\leq2C_1^{ \frac{9}{16}}C_2^{\frac{1}{16}}\varepsilon^{-{ \frac38}}\|\mathbb{D}_\varepsilon-\mathbb{D}\|_{L^\infty(\Omega)}\|\nabla\varphi\|_{L^\infty(\Omega\times(0,T))}
\\[5pt]
&\quad+2\left|\int_0^T\int_\Omega\left(\sqrt{1+c_{1\varepsilon}}\nabla\sqrt{1+c_{1\varepsilon}}-\sqrt{1+c_1}\nabla\sqrt{1+c_1}\right)\cdot\mathbb{D}\nabla\varphi\right|
\\[5pt]
&\leq { 6\,C_1^{\frac{9}{16}}C_2^{\frac{1}{16}} \varepsilon^{\frac{5}{8}}}\,\|\nabla\varphi\|_{L^\infty(\Omega\times(0,T))}\\
&\quad+2\left|\int_0^T\int_\Omega\left(\sqrt{1+c_{1\varepsilon}}\nabla\sqrt{1+c_{1\varepsilon}}-\sqrt{1+c_1}\nabla\sqrt{1+c_1}\right)\cdot\mathbb{D}\nabla\varphi\right|,
\end{align*}
for all $\varepsilon \in (0, 1)$. To deal with the last term, we first observe from \eqref{weak_conv1}
that
\begin{align*}
\left\|\sqrt{1+c_{1\varepsilon}}-\sqrt{1+c_1}\right\|_{L^4(\Omega\times(0,T))}^4
\leq\int_0^T\int_\Omega|c_{1\varepsilon}-c_1|^2\to0
\end{align*}
as $\varepsilon=\varepsilon_j\searrow0$. Hence
\begin{align}
\label{sqrt_strong4}
\sqrt{1+c_{1\varepsilon}}\to\sqrt{1+c_1}\quad\text{in}~L^4(\Omega\times(0,T)).
\end{align}
Combining \eqref{sqrt_strong4} with the weak convergence
\eqref{weak_conv3}, we obtain
\begin{align*}
\int_0^T\int_\Omega\sqrt{1+c_{1\varepsilon}}\nabla\sqrt{1+c_{1\varepsilon}}\cdot\mathbb{D}\nabla\varphi\to\int_0^T\int_\Omega\sqrt{1+c_1}\nabla\sqrt{1+c_1}\cdot\mathbb{D}\nabla\varphi
\end{align*}
as $\varepsilon=\varepsilon_j\searrow0$. Together with \eqref{epassumption3}, this establishes \eqref{c1_convergence1}. Similarly, from \eqref{epassumption3} and \eqref{convergence24} we have
\begin{align*}
  & \left|\int_0^T \int_\Omega c_{1\varepsilon} \nabla \tau_\varepsilon \cdot \mathbb{D}_\varepsilon \nabla \varphi - \int_0^T \int_\Omega c_1 \nabla \tau \cdot \mathbb{D} \nabla \varphi\right|   \\[5pt]
  &\leq \left|\int_0^T \int_\Omega c_{1\varepsilon} \nabla \tau_\varepsilon \cdot \mathbb{D}_\varepsilon \nabla \varphi -  \int_0^T \int_\Omega c_{1} \nabla \tau_\varepsilon \cdot \mathbb{D}_\varepsilon \nabla \varphi\right|    \\[5pt]
  &  + \left| \int_0^T \int_\Omega c_{1} \nabla \tau_\varepsilon \cdot \mathbb{D}_\varepsilon \nabla \varphi - \int_0^T \int_\Omega c_{1} \nabla \tau_\varepsilon \cdot \mathbb{D}\nabla \varphi\right| \\[5pt]
  & + \left| \int_0^T \int_\Omega c_{1} \nabla \tau_\varepsilon \cdot \mathbb{D}\nabla \varphi - \int_0^T \int_\Omega c_{1} \nabla \tau \cdot \mathbb{D}\nabla \varphi\right| \\[5pt]
  & \leq \left(\int_0^T \int_\Omega \nabla \tau_\varepsilon \cdot \mathbb{D}_\varepsilon \nabla \tau_\varepsilon\right)^{\frac{1}{2}} \left(\int_0^T \int_\Omega (c_1 - c_{1\varepsilon})^2 (\nabla \varphi \cdot \mathbb{D}_\varepsilon \nabla \varphi)\right)^{\frac{1}{2}} \\[5pt]
  & + \|c_{1}\|_{L^2(\Omega \times (0, T))} \|\nabla \tau_\varepsilon\|_{L^2(\Omega \times (0, T))}\|\mathbb{D}_\varepsilon - \mathbb{D}\|_{L^\infty(\Omega)} \|\nabla \varphi\|_{L^\infty(\Omega \times (0, T))}   \\[5pt]
  & + \left| \int_0^T \int_\Omega c_{1} \nabla \tau_\varepsilon \cdot \mathbb{D}\nabla \varphi - \int_0^T \int_\Omega c_{1} \nabla \tau \cdot \mathbb{D}\nabla \varphi\right| \\[5pt]
  & \leq \sqrt{C_1} (\|\mathbb{D}_\varepsilon\|_{L^\infty(\Omega)})^{\frac{1}{2}}\|\nabla \varphi\|_{L^\infty(\Omega)}\|c_1 - c_{1\varepsilon}\|_{L^2(\Omega \times (0, T))} \\[5pt]
  & + 3 \sqrt{C_1} \sqrt{\varepsilon} \|c_1\|_{L^2(\Omega \times (0, T))} \|\nabla \varphi\|_{L^\infty(\Omega \times (0, T))} \\[5pt]
  & + \left| \int_0^T \int_\Omega c_{1} \nabla \tau_\varepsilon \cdot \mathbb{D}\nabla \varphi - \int_0^T \int_\Omega c_{1} \nabla \tau \cdot \mathbb{D}\nabla \varphi\right|
\end{align*}
for all $\varepsilon \in (0, 1)$. The first term on the right-hand side converges to zero by \eqref{weak_conv1} and the uniform boundedness of $\mathbb{D}_\varepsilon$, while the second term converges to zero as $\varepsilon=\varepsilon_j\searrow0$. Moreover, since
\begin{align*}
c_1\nabla\varphi\in L^2(\Omega\times(0,T))
\end{align*}
and $\mathbb{D}$ is bounded, the weak convergence \eqref{weak_conv7} yields
\begin{align*}
\int_0^T\int_\Omega
c_1(\nabla\tau_\varepsilon-\nabla\tau) \cdot\mathbb{D}\nabla\varphi \to0.
\end{align*}
This establishes \eqref{tau_convergence1}. Now, testing the first equation of \eqref{app_prob} with $\varphi$ yields
\begin{align}
    \label{weak_c11}
    \nonumber & \int_0^\infty \int_\Omega c_{1\varepsilon} \partial_t \varphi \ { +} \int_\Omega c_{10\varepsilon} \varphi(\cdot, 0) = 2 \int_0^\infty \int_\Omega \sqrt{1 + c_{1\varepsilon}} \nabla \sqrt{1 + c_{1\varepsilon}} \cdot \mathbb{D}_\varepsilon \nabla \varphi\\[5pt]
    \nonumber & + \int_0^\infty \int_\Omega c_{1\varepsilon} (\nabla \cdot \mathbb{D}_\varepsilon) \cdot \nabla \varphi - \chi \int_0^\infty \int_\Omega c_{1\varepsilon} \nabla \tau_\varepsilon \cdot \mathbb{D}_\varepsilon \nabla \varphi + \alpha_1  \int_0^\infty \int_\Omega \tfrac{c_{1\varepsilon}}{1 + c_{1\varepsilon}} \varphi   \\[5pt]
    & - \alpha_2  \int_0^\infty \int_\Omega \tfrac{c_{2\varepsilon}}{1 + c_{2\varepsilon}}\varphi - \beta \int_0^\infty \int_\Omega c_{1\varepsilon}(1 - c_{1\varepsilon} - c_{2\varepsilon} - \tau_\varepsilon)\varphi
  \end{align}
for all $\varepsilon \in (0, 1)$. Here by \eqref{c1_convergence1} and \eqref{tau_convergence1} we have %{\cred where do we need \eqref{weak_assup2} here?} {\color{blue}\eqref{weak_assup2} had the boundedness information, but we can exclude it so removed.}
\begin{align}
  \label{weak_c12}
  \int_0^\infty \int_\Omega c_{1\varepsilon} \partial_t \varphi { \ +} \int_\Omega c_{10\varepsilon} \varphi(\cdot, 0) \to \int_0^\infty \int_\Omega c_{1} \partial_t \varphi { \ +} \int_\Omega c_{10} \varphi(\cdot, 0)
\end{align}
and
\begin{align*}
  2 \int_0^\infty \int_\Omega \sqrt{1 + c_{1\varepsilon}} \nabla \sqrt{1 + c_{1\varepsilon}} \cdot \mathbb{D}_\varepsilon \nabla \varphi \to 2 \int_0^\infty \int_\Omega \sqrt{1 + c_{1}} \nabla \sqrt{1 + c_{1}} \cdot \mathbb{D} \nabla \varphi
\end{align*}
and
\begin{align*}
   \chi \int_0^\infty \int_\Omega c_{1\varepsilon} \nabla \tau_\varepsilon \cdot \mathbb{D}_\varepsilon \nabla \varphi \to \chi \int_0^\infty \int_\Omega c_1 \nabla \tau \cdot \mathbb{D} \nabla \varphi
\end{align*}
as $\varepsilon = \varepsilon_j \searrow 0$, whereas by Lemma \ref{eplemma1} and \eqref{weak_conv1}
\begin{align}
  \label{weak_c14}
  \int_0^\infty \int_\Omega c_{1\varepsilon} (\nabla \cdot \mathbb{D}_\varepsilon) \cdot \nabla \varphi \to \int_0^\infty \int_\Omega c_{1} (\nabla \cdot \mathbb{D}) \cdot \nabla \varphi
\end{align}
as $\varepsilon = \varepsilon_j \searrow 0$. Furthermore, since the function $s\mapsto\frac{s}{1+s}$
is continuous and bounded on $[0,\infty)$, the convergence properties \eqref{weak_conv1} and \eqref{weak_conv4}, together with the dominated convergence theorem, yield
\begin{align}
  \label{weak_c15}
  \alpha_1  \int_0^\infty \int_\Omega \tfrac{c_{1\varepsilon}}{1 + c_{1\varepsilon}} \varphi  - \alpha_2  \int_0^\infty \int_\Omega \tfrac{c_{2\varepsilon}}{1 + c_{2\varepsilon}}\varphi \to \alpha_1  \int_0^\infty \int_\Omega \tfrac{c_{1}}{1 + c_{1}} \varphi  - \alpha_2  \int_0^\infty \int_\Omega \tfrac{c_{2}}{1 + c_{2}}\varphi
\end{align}
and
\begin{align}
  \label{weak_c16}
  \beta \int_0^\infty \int_\Omega c_{1\varepsilon}(1 - c_{1\varepsilon} - c_{2\varepsilon} - \tau_\varepsilon)\varphi \to \beta \int_0^\infty \int_\Omega c_{1}(1 - c_{1} - c_{2} - \tau)\varphi
\end{align}
as $\varepsilon = \varepsilon_j \searrow 0$. The convergence results \eqref{weak_c12}-\eqref{weak_c16} allow us to pass to the limit $\varepsilon = \varepsilon_j \searrow 0$ in \eqref{weak_c11} to obtain
\begin{align}
    \label{12weak_c1}
    \nonumber & \int_0^\infty \int_\Omega c_1 \partial_t \varphi { \ +} \int_\Omega c_{10} \varphi(\cdot, 0) = 2 \int_0^\infty \int_\Omega \sqrt{1 + c_1} \nabla \sqrt{1 + c_1} \cdot \mathbb{D} \nabla \varphi \\[5pt]
    \nonumber & + \int_0^\infty \int_\Omega c_1 (\nabla \cdot \mathbb{D}) \cdot \nabla \varphi - \chi \int_0^\infty \int_\Omega c_1 \nabla \tau \cdot \mathbb{D} \nabla \varphi + \alpha_1  \int_0^\infty \int_\Omega \tfrac{c_{1}}{1 + c_{1}} \varphi  \\[5pt]
    & - \alpha_2  \int_0^\infty \int_\Omega \tfrac{c_{2}}{1 + c_{2}}\varphi - \beta \int_0^\infty \int_\Omega c_{1}(1 - c_{1} - c_{2} - \tau)\varphi.
  \end{align}
Similarly, by testing the second and third equations in \eqref{app_prob} against $\varphi$, then combining \eqref{weak_assup2}, \eqref{weak_conv1}, \eqref{weak_conv2}, \eqref{weak_conv4}, and \eqref{weak_conv6}, followed by passing to the limit $\varepsilon = \varepsilon_j \searrow 0$, we deduce:
\begin{align}
    \label{12weak_c2}
    \int_0^\infty \int_\Omega c_2 \partial_t \varphi + \int_\Omega c_{20} \varphi(\cdot, 0) = -\alpha_1 \int_0^\infty \int_\Omega \tfrac{c_1}{1 + c_1} \varphi + \alpha_2 \int_0^\infty \int_\Omega \tfrac{c_2}{1 + c_2}\varphi,
  \end{align}
  and
  \begin{align}
    \label{12weak_tau}
    \int_0^\infty \int_\Omega \tau \partial_t \varphi + \int_\Omega \tau_{0} \varphi(\cdot, 0) =  \delta \int_0^\infty \int_\Omega c_1 \tau \varphi + \mu \int_0^\infty \int_\Omega \tau \varphi - \int_0^\infty \int_\Omega \tfrac{c_2 }{1 + c_2}\varphi.
  \end{align}
  Collecting \eqref{12weak_c1}, \eqref{12weak_c2}, and \eqref{12weak_tau} leads to \eqref{weak_c1}-\eqref{weak_tau}, thereby completing the proof.
\end{proof}
\section{Numerical simulations}\label{sec:numerics} 
For the numerical treatment, we consider the following nondimensionalized model:
\begin{equation}\label{model_num}
  \begin{cases}
    \partial_t c_{1} = {\nabla \nabla : (\mathbb Dc_1)} - \chi \nabla \cdot (c_1 \mathbb{D} \nabla \tau)- \alpha_1  \frac{c_{1}}{1 + c_{1}}  + \alpha_2 \frac{c_{2}}{1 + c_{2}} \\[10pt]
    \hspace{8cm}+  c_{1}\left(1 - c_{1} - c_{2} - \tau\right),                                                                                    &   x \in \Omega, t > 0,  \\[10pt]
    \partial_t c_2 =  \alpha_1 \frac{c_1}{1 + c_1} - \alpha_2 \frac{c_2}{1 + c_2},                                                                 &   x \in \Omega, t > 0,   \\[5pt]
    \partial_t \tau = -\delta_1 c_1 \tau - \mu \tau +  \delta_2 \frac{c_2}{1 + c_2},                                                               &   x \in \Omega, t > 0,
  \end{cases}
\end{equation}
equipped with the zero-flux boundary conditions \eqref{boundary_conditions} on $\partial\Omega \times (0,T)$ and initial conditions \eqref{initial_conditions}. The parameters in \eqref{model_num} are originally denoted by $\widehat{U}$, where $U$ represents a parameter origanally appearing in \eqref{model} respectively. For notational convenience, we omit the hats throughout. Details of the nondimensionalized parameters are provided in Table~\ref{tab:params}.

Numerically approximating system \eqref{model_num} with a degenerate diffusion tensor poses several structural challenges. Where $\mathbb{D}$ has vanishing or arbitrarily small eigenvalues, uniform parabolicity is lost, weakening the regularizing effect of diffusion in the corresponding directions. Sharp spatial variations in $\mathbb{D}$ can further complicate the discretization and contribute to steep solution gradients or sharp propagation fronts. In particular, the advective contribution $\nabla \cdot (c_1 \nabla \cdot \mathbb{D})$ to the myopic diffusion $\nabla \nabla : (\mathbb Dc_1)$ and the haptotaxis $-\chi \nabla \cdot (c_1 \mathbb{D} \nabla \tau)$ can be combined into a drift contribution in divergence form $\nabla \cdot (c_1 \mathbf{w})$ where
$\mathbf{w} := \nabla \cdot \mathbb{D} - \chi \mathbb{D} \nabla \tau$, denotes the auxiliary drift vector field. The drift contribution involving $\nabla \cdot \mathbb{D}$ may become particularly pronounced where $\mathbb{D}$ varies steeply. Moreover, because $c_2$ and $\tau$ satisfy local ODEs without spatial diffusion, perturbations originating in $c_1$ are not spatially regularized by diffusion in these components. Numerical oscillations or spurious leakage of $c_1$ across degenerate interfaces can therefore propagate to $c_2$ and $\tau$ through the coupling terms, be amplified by the reaction dynamics, or contaminate the computed haptotactic gradient $\nabla \tau$.

Let $\mathcal{T}_h$ be a conforming triangulation of the computational domain $\Omega = [-L_{\text{dom}}, L_{\text{dom}}]^2$. For the spatial approximation of the cell densities $c_1, c_2$ and the ECM $\tau$, we employ standard continuous piecewise linear Lagrange finite elements:
\begin{equation}
  V_h := \left\{ v_h \in C^0(\overline{\Omega}) \;:\; v_h|_K \in \mathbb{P}_1(K) \quad \text{for all}~ K \in \mathcal{T}_h \right\} \subset H^1(\Omega).
\end{equation}
To accurately capture the mixed transport and gradient effects, the auxiliary drift field $w$ is discretized using the lowest-order Raviart--Thomas finite element space:
\begin{equation}
  W_h := \left\{ w_h \in H(\mathrm{div}; \Omega) \;:\; w_h|_K \in [\mathbb{P}_{0}(K)]^2 + \mathbf{x}\,\mathbb{P}_{0}(K) \quad \text{for all}~ K \in \mathcal{T}_h \right\}.
\end{equation}
Here, $[\mathbb{P}_0(K)]^2 = \mathbb{P}_0(K) \times \mathbb{P}_0(K)$ denotes the space of piecewise constant 2D vectors on $K$, and $\mathbf{x} = (x_1, x_2)^\top \in K \subset \mathbb{R}^2$ is the spatial position vector.
For the temporal discretization, we divide $[0,T]$ into uniform time steps $t_n = n \Delta t$ ($n=0,1,\dots,N$) with step size $\Delta t > 0$. Multiplying \eqref{model_num} by test functions $\varphi_{1,h}, \varphi_{2,h}, \varphi_{3,h} \in V_h$, integrating by parts, and applying a backward Euler scheme combined with a semi-implicit linearization and mass lumping, with the aim of improving positivity preservation and numerical robustness near degenerate interfaces, we obtain the following fully discrete formulation:

\begin{definition}
Given $c_{1,h}^n, c_{2,h}^n, \tau_h^n \in V_h$ at time step $t_n$, find $c_{1,h}^{n+1}, c_{2,h}^{n+1}, \tau_h^{n+1} \in V_h$ such that for all $\varphi_{1,h}, \varphi_{2,h}, \varphi_{3,h} \in V_h$:
\begin{subequations}\label{fully_discrete}
\begin{align}
\left(\frac{c_{1,h}^{n+1} - c_{1,h}^n}{\Delta t}, \varphi_{1,h}\right)_h
+ a_0(c_{1,h}^{n+1}, \varphi_{1,h})
+ \int_\Omega c_{1,h}^{n+1} (w_h^n \cdot \nabla \varphi_{1,h})\,dx
&= \ell_1(c_{1,h}^{n+1}, c_{1,h}^n, c_{2,h}^n, \tau_h^n; \varphi_{1,h}),
\\
\left(\frac{c_{2,h}^{n+1} - c_{2,h}^n}{\Delta t}, \varphi_{2,h}\right)_h
&= \ell_2(c_{1,h}^n, c_{2,h}^{n+1}, c_{2,h}^n; \varphi_{2,h}),
\\
\left(\frac{\tau_h^{n+1} - \tau_h^n}{\Delta t}, \varphi_{3,h}\right)_h
&= \ell_3(c_{1,h}^n, c_{2,h}^n, \tau_h^{n+1}; \varphi_{3,h}),
\end{align}
\end{subequations}
where the Fickian diffusion bilinear form is given by
\begin{equation}
  a_0(c_{1,h}^{n+1}, \varphi_{1,h}) := \int_\Omega \mathbb{D}\nabla c_{1,h}^{n+1} \cdot \nabla \varphi_{1,h}\,dx,
\end{equation}
and the discrete drift field $w_h^n \in W_h$ is updated at each time step via projection into $W_h$:
\begin{equation}
  w_h^n = \Pi_{W_h}\big(\nabla \cdot \mathbb{D} - \chi \mathbb{D} \nabla \tau_h^n\big).
\end{equation}
The semi-implicitly linearized reaction functionals (evaluated using mass lumping) are defined as
\begin{align}
\ell_1(c_{1,h}^{n+1}, c_{1,h}^n, c_{2,h}^n, \tau_h^n; \varphi_{1,h}) 
&:= \left( -\alpha_1 \frac{c_{1,h}^{n+1}}{1+c_{1,h}^n} + \alpha_2 \frac{c_{2,h}^n}{1+c_{2,h}^n} + c_{1,h}^{n+1}(1 - c_{1,h}^n - c_{2,h}^n - \tau_h^n), \, \varphi_{1,h} \right)_h, \\
\ell_2(c_{1,h}^n, c_{2,h}^{n+1}, c_{2,h}^n; \varphi_{2,h}) 
&:= \left( \alpha_1 \frac{c_{1,h}^n}{1+c_{1,h}^n} - \alpha_2 \frac{c_{2,h}^{n+1}}{1+c_{2,h}^n}, \, \varphi_{2,h} \right)_h, \\
\ell_3(c_{1,h}^n, c_{2,h}^n, \tau_h^{n+1}; \varphi_{3,h}) 
&:= \left( -\delta_1 c_{1,h}^n \tau_h^{n+1} - \mu \tau_h^{n+1} + \delta_2 \frac{c_{2,h}^n}{1+c_{2,h}^n}, \, \varphi_{3,h} \right)_h.
\end{align}
\end{definition}
Here, $(\cdot,\cdot)_h$ denotes the lumped $L^2(\Omega)$-inner product, evaluated using the vertex quadrature rule, which yields a diagonal lumped mass matrix:
\begin{equation}
  (u, v)_h := \sum_{K_A \in \mathcal{T}_h} \frac{|K_A|}{3} \sum_{i=1}^3 u(\mathbf{a}_{K_A,i})\,v(\mathbf{a}_{K,i}) 
  \qquad \text{for all } u, v \in C^0(\overline{\Omega}),
\end{equation}
where $|K_A|$ denotes the area of the triangle $K_A \in \mathcal{T}_h$ and $\left\{\mathbf{a}_{K_A,1}, \mathbf{a}_{K_A,2}, \mathbf{a}_{K_A,3}\right\}$ are its vertices. Mass lumping has been employed in finite element schemes for degenerate parabolic equations; see \cite{MizutaniSaitoSuzuki2005}. A related nodal quadrature-based inner product is used in finite element approximations of degenerate Cahn-Hilliard systems; see \cite{BarrettBloweyGarcke2001}.

\subsection{Diffusion tensor scenarios and initial conditions}

For the dimensionless diffusion tensor on $\Omega = [-L_{\text{dom}}, L_{\text{dom}}]^2$, we define
\begin{equation}
\mathbb{D}(x,y) = a_1 \Big( D_{\min}\mathbf{I} + g(x,y)\mathbb{A} \Big),
\end{equation}
where $a_1 = 0.20$ denotes the tissue porosity/tortuosity prefactor, $D_{\min} \ge 0$ acts as a background diffusion coefficient, and the constant anisotropy matrix $\mathbb{A}$ is given by
\begin{equation}
\mathbb{A} = 
\begin{pmatrix}
1 & \epsilon\\
\epsilon & \alpha
\end{pmatrix},
\qquad \alpha > 0, \quad \epsilon^2 < \alpha,
\end{equation}
with $\alpha = 0.8$ and $\epsilon = 0.4$ in all simulations. 

To investigate how microstructural tissue properties and varying motility constraints affect cell migration, tissue regeneration, and haptotactic signaling, we introduce and compare three diffusion-tensor scenarios. These configurations allow us to contrast classical diffusion with constrained fronts, as well as complex heterogeneous tissue architectures containing physical barriers:

\begin{itemize}
    \item \textbf{Case 1 (Non-degenerate anisotropic diffusion tensor):} We set a strictly positive background diffusivity $D_{\min} = 0.10$ and modulate the anisotropic part by a radially symmetric $C^2$-function with compact support:
    \begin{equation}\label{eq:g_radial}
    g(x,y) = 
    \begin{cases} 
    D_{\max} \left( 1 - \left( \frac{r}{R} \right)^2 \right)^3 & \text{if } r < R, \\
    0 & \text{if } r \ge R,
    \end{cases}
    \end{equation}
    with radial distance $r = \sqrt{(x-\mu_x)^2 + (y-\mu_y)^2}$ centered at $(\mu_x, \mu_y) = (0,0)$, amplitude $D_{\max} = 1.0$, and cutoff radius $R = 2.5$. Since $D_{\min} > 0$ holds uniformly across the domain, $\mathbb{D}(x,y)$ remains strictly positive definite, preventing any parabolic degeneracy and ensuring smooth, unhindered spreading across the entire spatial domain $\overline{\Omega}$.

    \item \textbf{Case 2 (Radially degenerate diffusion tensor with compact support):} Utilizing the same underlying radial profile $g(x,y)$ from \eqref{eq:g_radial} with $R = 2.5$, we effectively set the background diffusion strictly to zero ($D_{\min} := 0$). Consequently, the diffusion tensor degenerates identically to $\mathbb{D}(x,y) \equiv \mathbb O_2$ in the outer ring region $r \ge R$. This configuration models a localized active tissue zone surrounded by an impermeable, non-diffusive matrix boundary.
    \item \textbf{Case 3 (Strongly degenerate multi-hole diffusion tensor):} To capture highly heterogeneous biological microenvironments, we set $D_{\min} := 0$ and construct a synthetic test profile featuring multiple internal impenetrable obstacles where $\mathbb{D}(x,y) \equiv \mathbb{O}_2$. These are interspersed with interconnected migration channels characterized by uniform baseline motility. Specifically, the spatial modulation is governed by
\begin{equation}\label{eq:g_multihole}
g(x,y) = D_{\max} \, E(r^2) \, H(x,y),
\end{equation}
where the flat-top radial envelope $E(r^2)$ maintains constant baseline diffusion inside $r \le R_{\text{flat}} = 2.0$ and vanishes smoothly at the outer boundary $R = 2.85$:
\begin{equation}
E(r^2) = \big( 1 - s_{\text{ext}}(r^2)^2 \big)^3, 
\qquad 
s_{\text{ext}}(r^2) = 
\begin{cases}
0 & \text{if } r^2 \le R_{\text{flat}}^2, \\
\dfrac{r^2 - R_{\text{flat}}^2}{R^2 - R_{\text{flat}}^2} & \text{if } R_{\text{flat}}^2 < r^2 < R^2, \\
1 & \text{if } r^2 \ge R^2,
\end{cases}
\end{equation}
while the periodic perforation function $H(x,y)$ incorporates internal structural holes via a smooth truncated trigonometric mapping:
\begin{equation}
H(x,y) = u(x,y)^2 \big( 3 - 2u(x,y) \big), 
\qquad 
u(x,y) = \min\!\left( 1, \, \max\!\left( 0, \, \frac{\cos(\omega x)\cos(\omega y) - r_0}{w} \right) \right),
\end{equation}
parameterized by obstacle frequency $\omega = 2.5$, cutoff threshold $r_0 = -0.35$, and transition width $w = 0.35$.
    
\end{itemize}

\begin{figure}[htbp]
    \centering
    \subfloat[Non-degenerate diffusion tensor profile for Case 1 ($D_{\min}=0.10$).]{%
        \includegraphics[width=0.82\linewidth]{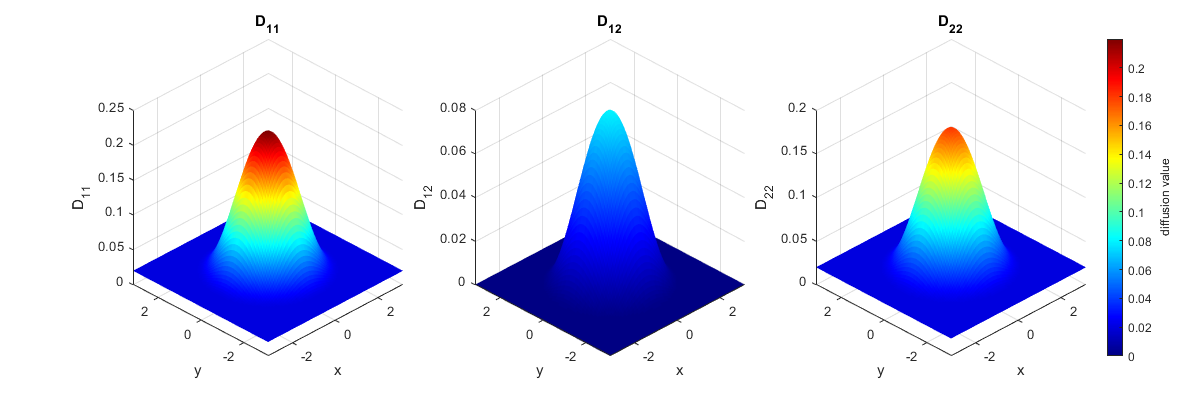}
        \label{fig:diff_nondeg}}

    \vspace{0.3cm}

    \subfloat[Radially degenerate diffusion tensor profile for Case 2 ($D_{\min}=0$).]{%
        \includegraphics[width=0.82\linewidth]{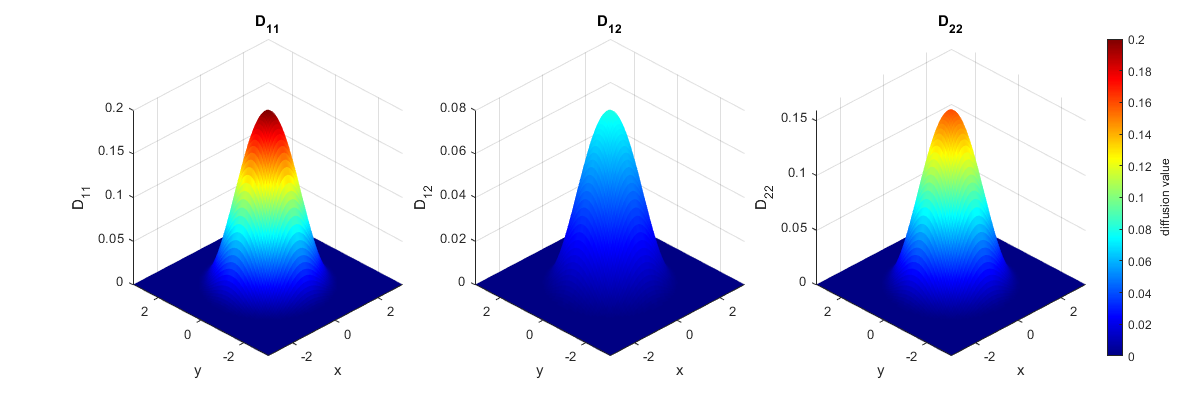}
        \label{fig:diff_deg}}

    \vspace{0.3cm}

    \subfloat[Strongly degenerate multi-hole diffusion tensor profile for Case 3 ($D_{\min}=0$).]{%
        \includegraphics[width=0.82\linewidth]{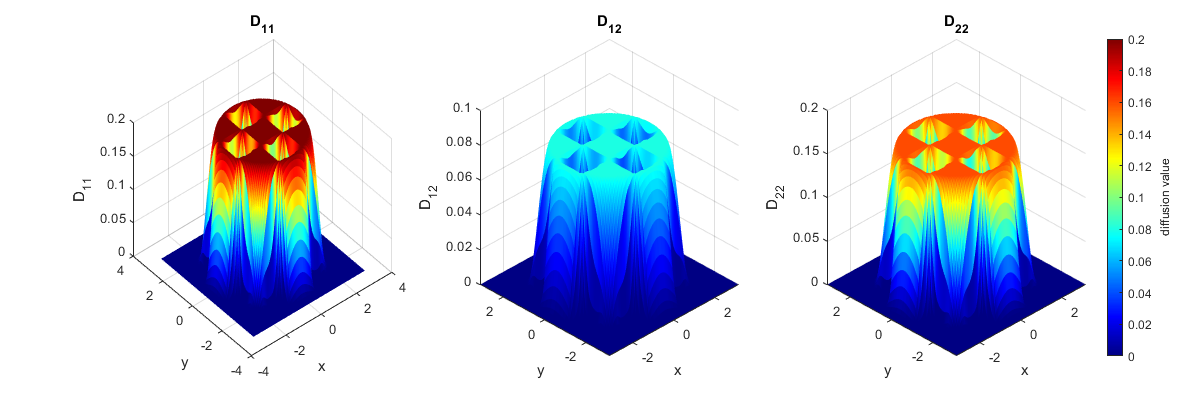}
        \label{fig:diff_case3}}

    \caption{Visual representation of the spatial distributions of the diffusion tensor components $D_{11}$, $D_{12}=D_{21}$, and $D_{22}$ across the three simulation scenarios.}
    \label{fig:tensor_profiles}
\end{figure}

\paragraph{Initial conditions:}
We prescribe initial conditions at $t = 0$ supported strictly within the active (non-degenerate) tissue subdomain $\Omega_+ := \{ (x,y) \in \Omega : D_{\min} + g(x,y) > 0 \}$. For the {MSCs of density} $c_1$, we model a localized cell seeding concentrated around the center $(0, 0)$ via a Gaussian distribution:
\begin{equation}
  c_1(x, y, 0) = K_{\text{car}} \exp\left( -\frac{x^2 + y^2}{\sigma_0^2} \right) \mathds{1}_{\Omega_+}(x,y), \quad ~\text{for all}~ \ (x,y) \in \Omega,
\end{equation}
where $K_{\text{car}} $ denotes the initial normalized seeding density and $\sigma_0^2 = 0.8$ characterizes the initial spread of the cell aggregate. 
For the chondrocytes of density $c_2$, we assume a minimal baseline presence proportional to the MSC distribution:
\begin{equation}
  c_2(x, y, 0) = 0.05 \, c_1(x, y, 0), \quad ~\text{for all}~ (x,y) \in \Omega.
\end{equation}
Finally, the ECM $\tau$ is initialized with a small random spatial perturbation to reflect baseline physiological variability within the viable tissue matrix:
\begin{equation}
  \tau(x, y, 0) = \tau_{\text{noise}}(x,y) \, \mathds{1}_{\Omega_+}(x,y), \quad \text{with } \tau_{\text{noise}}(x,y) \in [0, 10^{-3}], \quad ~\text{for all}~ (x,y) \in \Omega.
\end{equation}

%\newpage
\subsection{Computational results}

In this subsection, we present the computational results for all diffusion tensor scenarios. We directly compare the results of Case 1 and Case 2, and subsequently present the solution for the strongly degenerate case (Case 3). The physical computational domain is $1000\,\mu\text{m} \times 1000\,\mu\text{m}$, with $x, y \in [-500\,\mu\text{m}, +500\,\mu\text{m}]$. Scaling by the characteristic length $L = 158.11\,\mu\text{m}$ yields the dimensionless domain:
\begin{equation}
    \widehat{x}, \widehat{y} \in \left[ \frac{-500}{158.11}, \frac{+500}{158.11} \right] \approx [-3.16, +3.16].
\end{equation}

The physical parameters and corresponding dimensionless scaling values utilized in our computational framework are summarized in Table \ref{tab:params}. 
Time integration is carried out up to $N = 420$ discrete time steps with a uniform step size of $\Delta t = 0.01$, corresponding to $t = 4.2$ dimensionless time units, which represents an effective duration of $21\,\text{days}$.
For the spatial discretization in Cases 1 and 2, the domain is partitioned into a mesh consisting of $9{,}800$ triangular elements, yielding a total of $29{,}963$ degrees of freedom (DoFs). To adequately resolve the complex internal boundaries of the multi-hole obstacle geometry in Case 3, a finer mesh comprising $16{,}200$ elements and $49{,}323$ total DoFs is employed.

We perform the numerical simulations using the open-source finite element software FreeFem++ \cite{Hecht2012}. We use ParaView \cite{AhrensGL05} and MATLAB \cite{MATLAB:R2022b} for post-processing and visualization of the spatial fields.

\begin{table}[htbp]
\centering
\begin{tabular}{llll}
\toprule
\textbf{Parameters} & \textbf{Physical values} & \textbf{Scaling formula} & \textbf{Dimensionless value} \\
\midrule
 Carrying capacity $K$ & $10^{-3}\,\text{cells/}\mu\text{m}^3$ & Reference density & $1.0$  \\
Proliferation rate $\beta$ & $0.20\,\text{day}^{-1}$ & $T = 1/\beta = 5\,\text{days}$ & $1.0$ \\
Baseline diffusion $D_0$ & $5000\,\mu\text{m}^2/\text{day}$ & $L = \sqrt{D_0/\beta} = 158.11\,\mu\text{m}$ & $\wt{D}_{\max} = 1.0$ \\
Tissue factor $a_1$ & $0.20$ (dimensionless) & Explicit PDE prefactor & $a_1 = 0.20$ \\
Min.\ diffusion $D_{\min}$ & $0$ (or $500\,\mu\text{m}^2/\text{day}$) & $\wt{D}_{\min} = D_{\min}/D_0$ & $0.0$ (or $0.10$) \\
Support $R$ & $395 - 450\,\mu\text{m}$ & $\wt{R} = R/L$ & $2.50 - 2.85$ \\
Differentiation $\alpha_1$ & $2.4 \times 10^{-5}\,\text{cells/}(\mu\text{m}^3\text{d})$ & $\wt{\alpha}_1 = \alpha_1 / (\beta K)$ & $0.12$ \\
Dedifferentiation $\alpha_2$ & $4.0 \times 10^{-6}\,\text{cells/}(\mu\text{m}^3\text{d})$ & $\wt{\alpha}_2 = \alpha_2 / (\beta K)$ & $0.02$ \\
Haptotactic sensitivity $\chi$ & $50\,\mu\text{m}^3/\text{cells}$ & $\wt{\chi} = K \chi$ & $0.05$ \\
Signal uptake $\delta_1$ & $30\,\mu\text{m}^3/(\text{cells}\cdot\text{d})$ & $\wt{\delta}_1 = \delta_1 K / \beta$ & $0.15$ \\
Signal synthesis $\delta_2$ & $5.0 \times 10^{-5}\,\text{cells/}(\mu\text{m}^3\text{d})$ & $\wt{\delta}_2 = \delta_2 / (\beta K)$ & $0.25$ \\
Signal decay $\mu$ & $0.16\,\text{day}^{-1}$ & $\wt{\mu} = \mu / \beta$ & $0.80$ \\
\bottomrule
\end{tabular}
\caption{Parameter table for the model \eqref{model_num}; physical and dimensionless reference values are adapted from \cite{BailonPlaza2001, Campbell2019a}.}
\label{tab:params}
\end{table}

\subsubsection{Comparison of non-degenerate and radially degenerate diffusion (Cases 1 and 2)}

\begin{figure}[htbp]
	\centering
	\subfloat[$c_1(x,y,0)$]{%
		\includegraphics[width=0.32\textwidth]{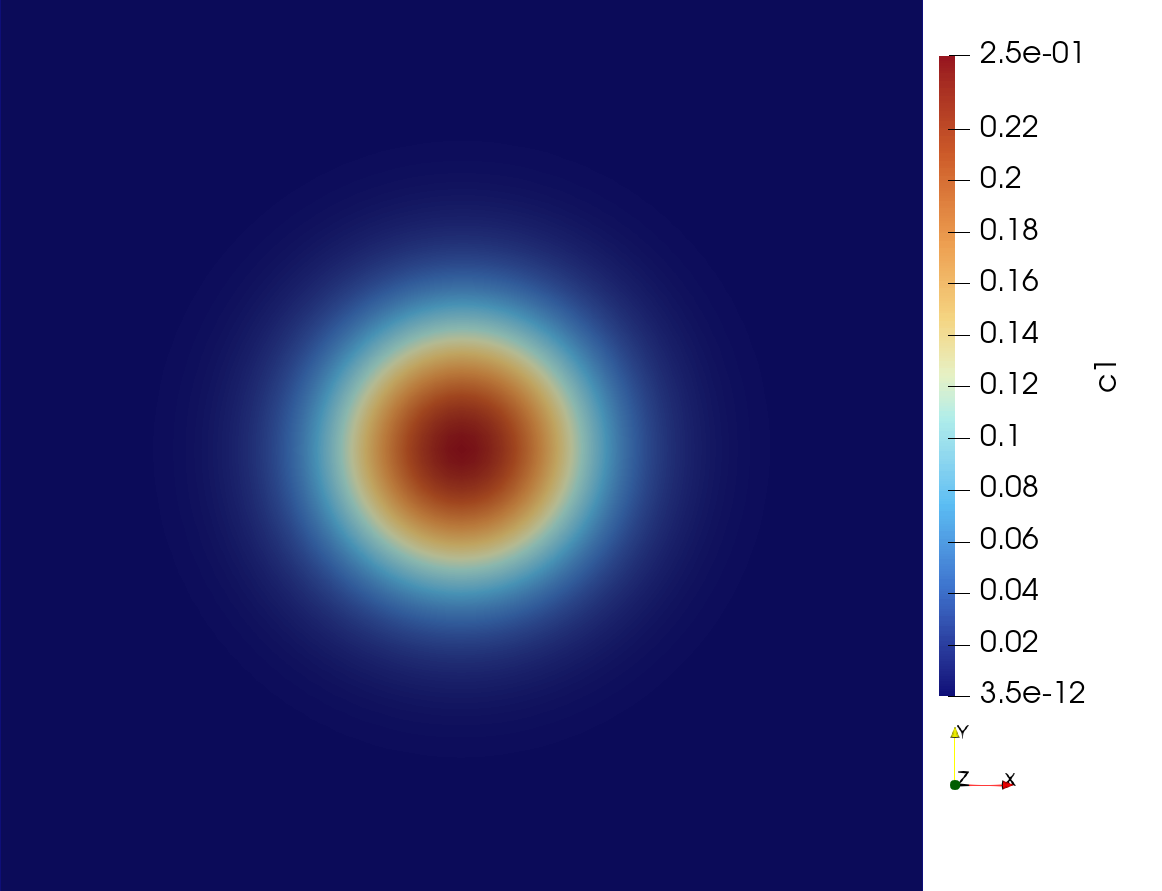}}
	\hfill
	\subfloat[$c_2(x,y,0)$]{%
		\includegraphics[width=0.32\textwidth]{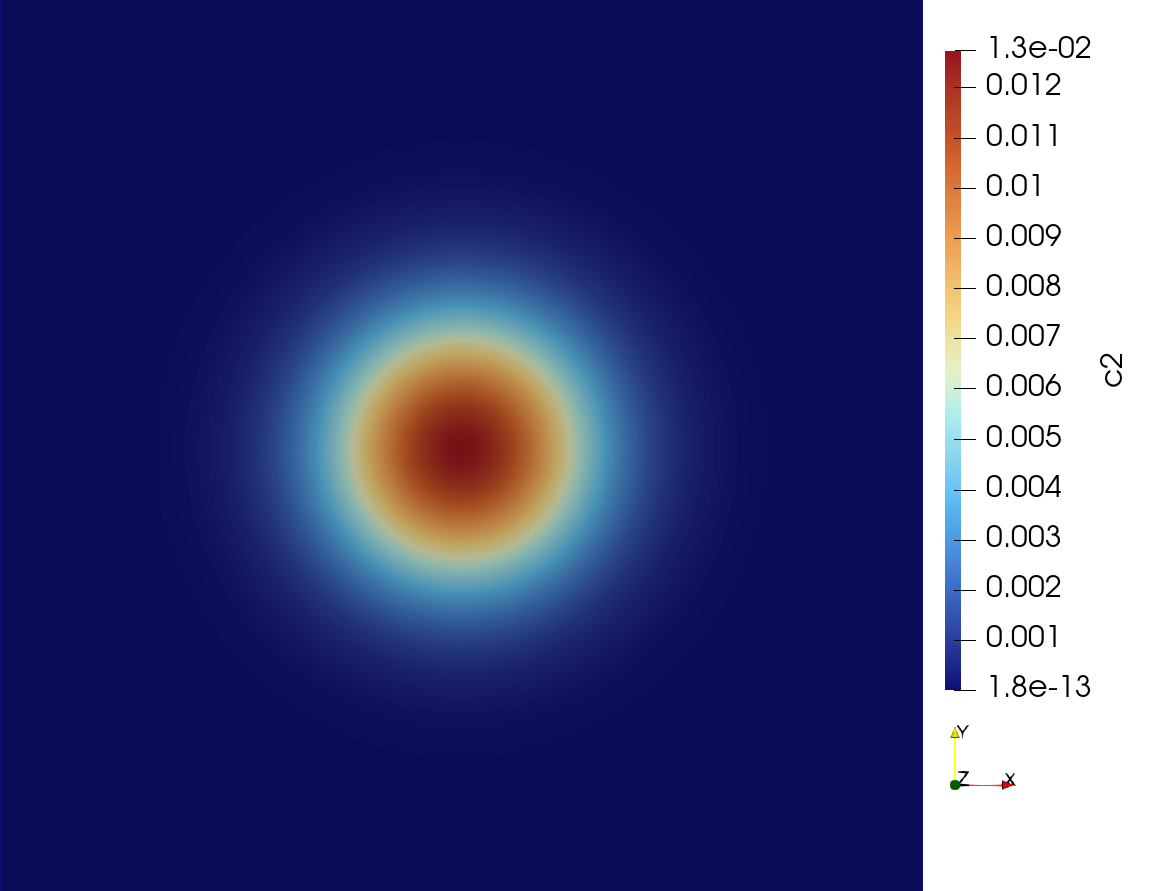}}
	\hfill
	\subfloat[$\tau(x,y,0)$]{%
		\includegraphics[width=0.32\textwidth]{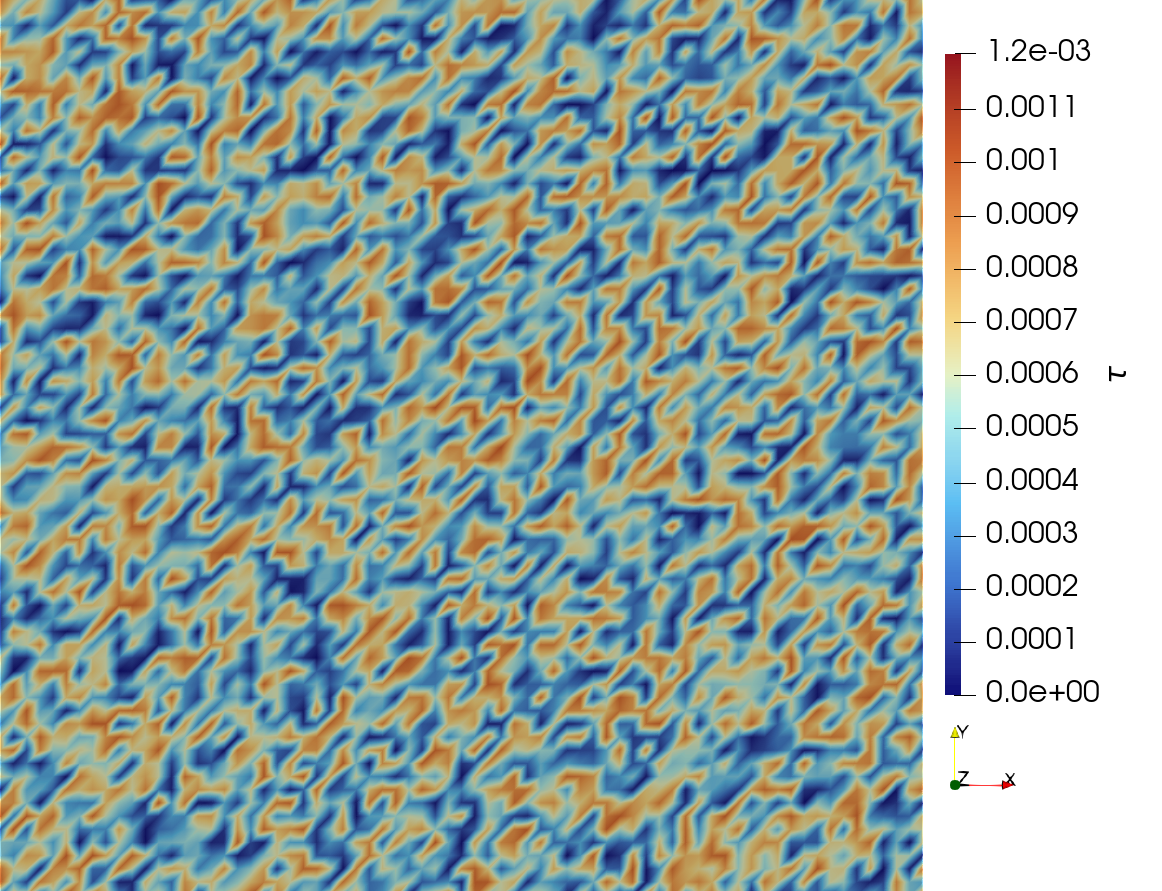}}
	\caption{Initial distributions of $c_1$, $c_2$, and $\tau$ at $t=0$ for Cases 1 and 2.}
	\label{fig:init_case12}
\end{figure}

We first examine the baseline behavior and compare it with the radially degenerate setting. Although initial cell seeding is restricted to the tissue subdomain \(\Omega_+\), we prescribe identical initial distributions throughout the computational domain \(\Omega\) for both Case 1 and Case 2. In Case 2, the degenerate outer ring effectively restricts cell presence to \(\Omega_+\). Nevertheless, prescribing the same initial distributions in both cases ensures a controlled and unbiased comparison, allowing us to isolate the effect of parabolic degeneracy on the resulting dynamics. Figure \ref{fig:init_case12} displays the prescribed initial distributions at $t = 0$.

%Figure \ref{fig:init_case12} displays the prescribed initial distributions of $c_1$, $c_2$, and $\tau$ at $t = 0$. 

\begin{figure}[htbp]
    \centering

    \subfloat[$c_1$ (non-degenerate)]{%
        \includegraphics[width=0.32\textwidth]{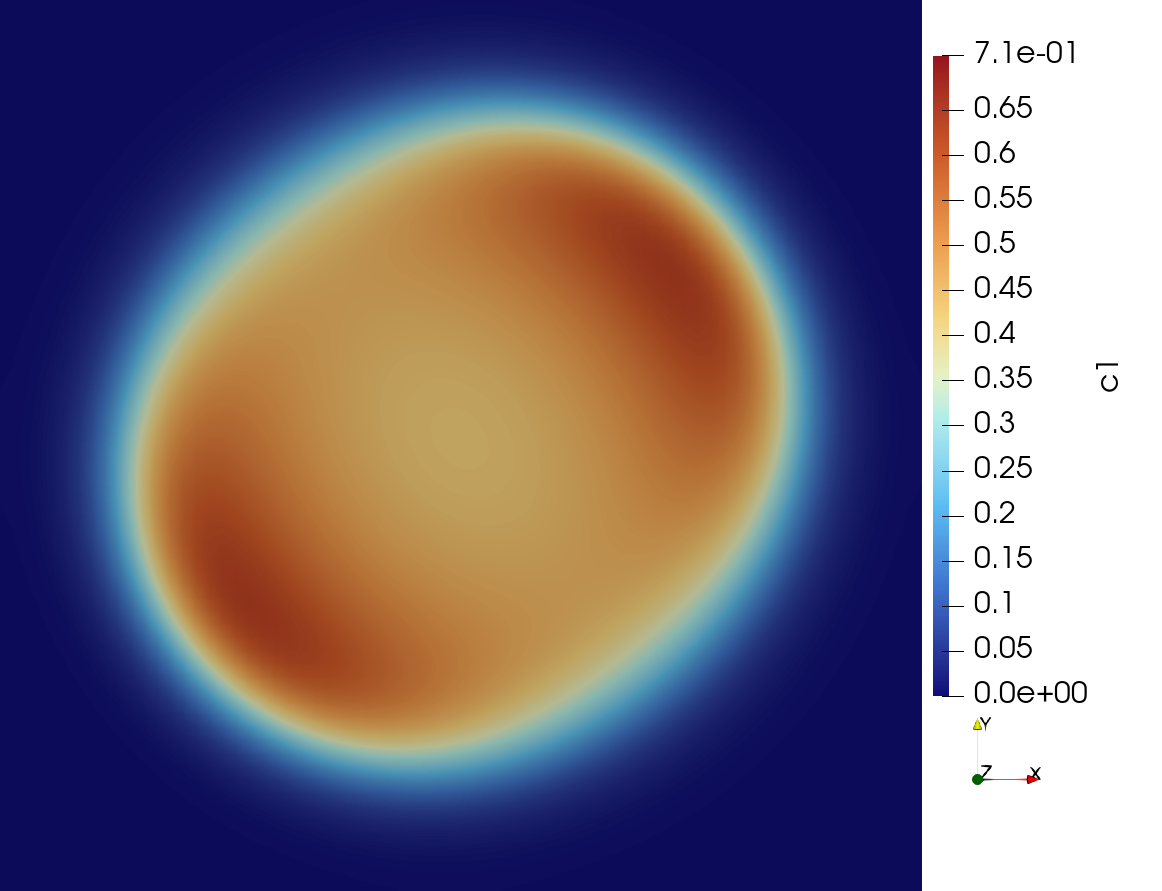}}
    \hfill
    \subfloat[$c_1$ (degenerate)]{%
        \includegraphics[width=0.32\textwidth]{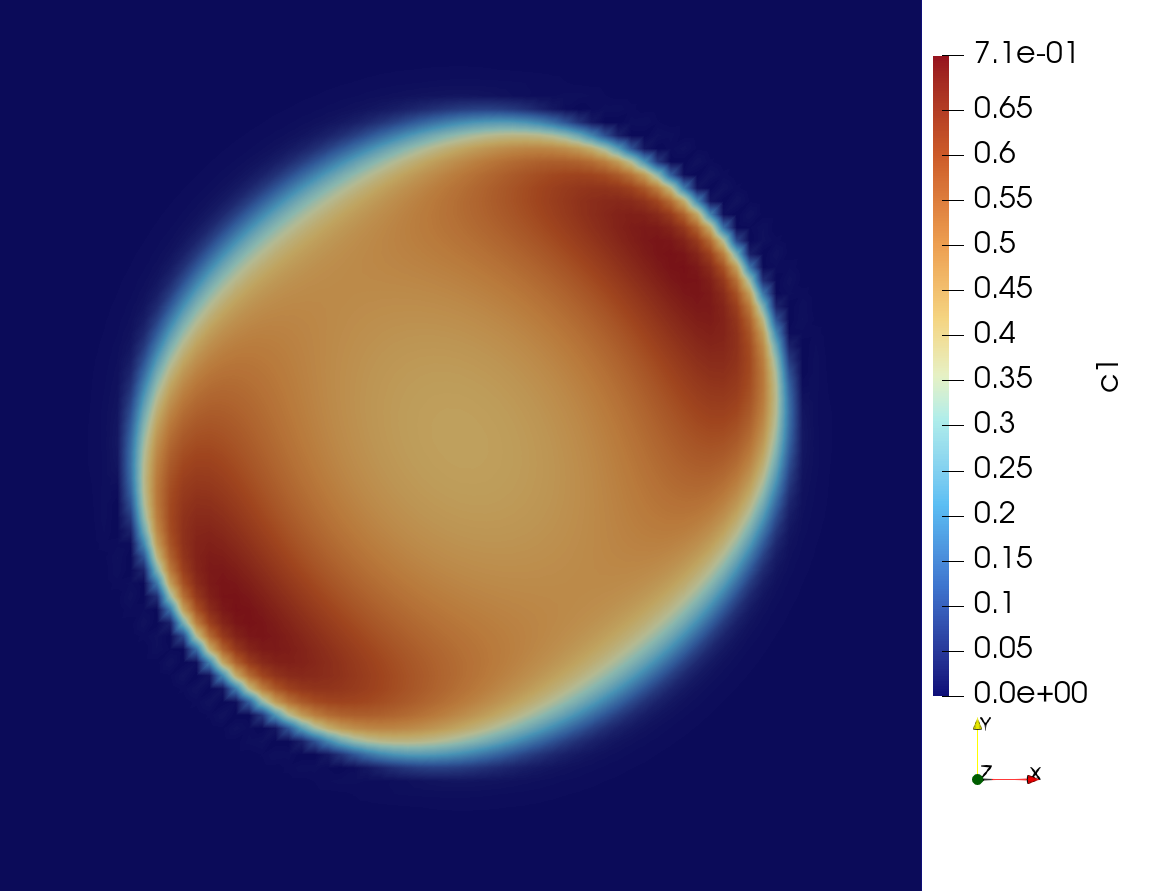}}
    \hfill
    \subfloat[$c_1$ difference]{%
        \includegraphics[width=0.32\textwidth]{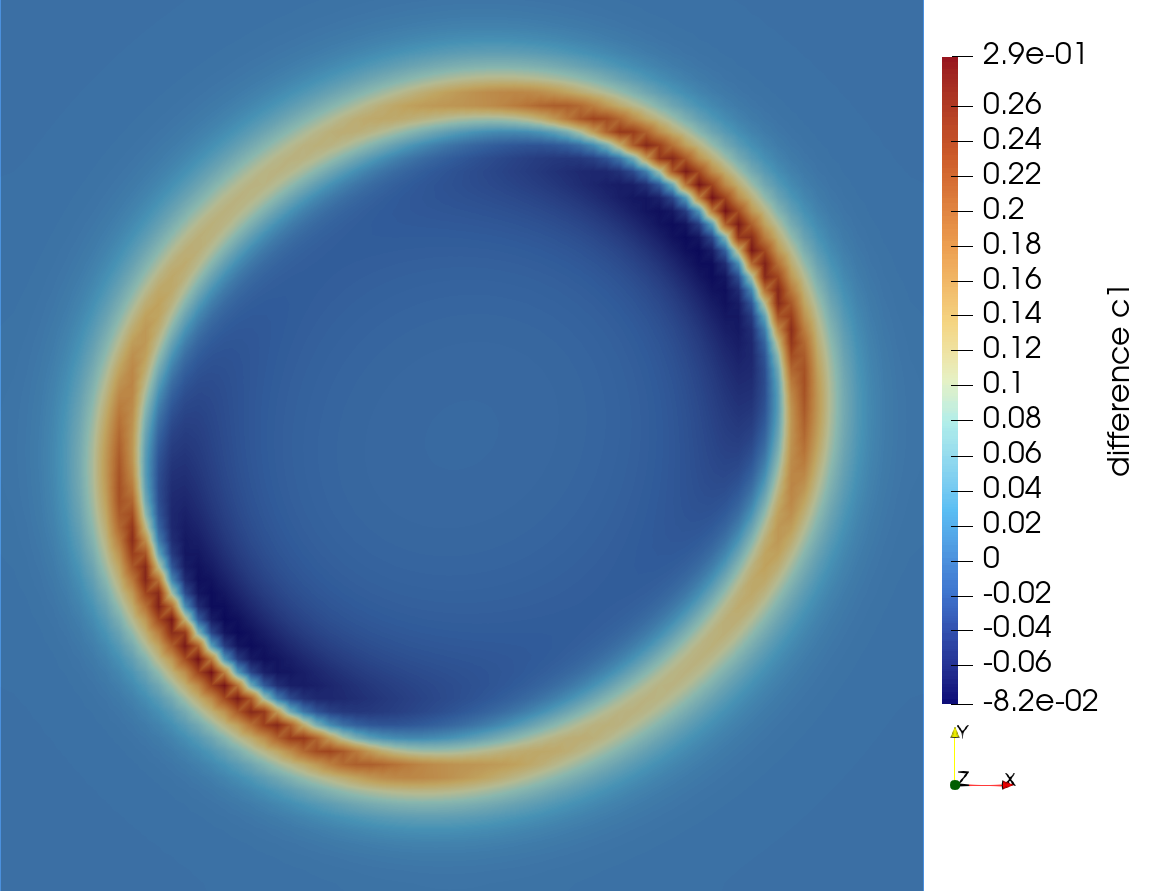}}

    \vspace{0.3cm}

    \subfloat[$c_2$ (non-degenerate)]{%
        \includegraphics[width=0.32\textwidth]{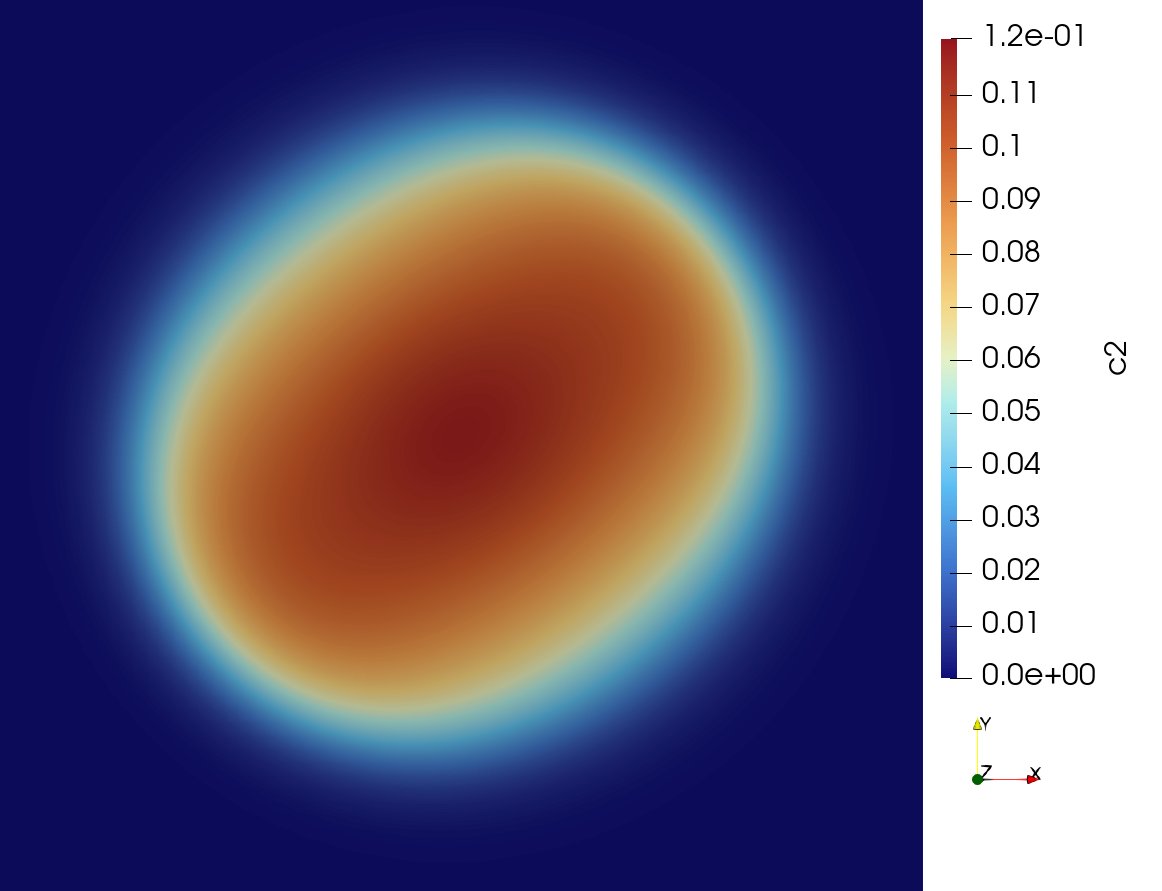}}
    \hfill
    \subfloat[$c_2$ (degenerate)]{%
        \includegraphics[width=0.32\textwidth]{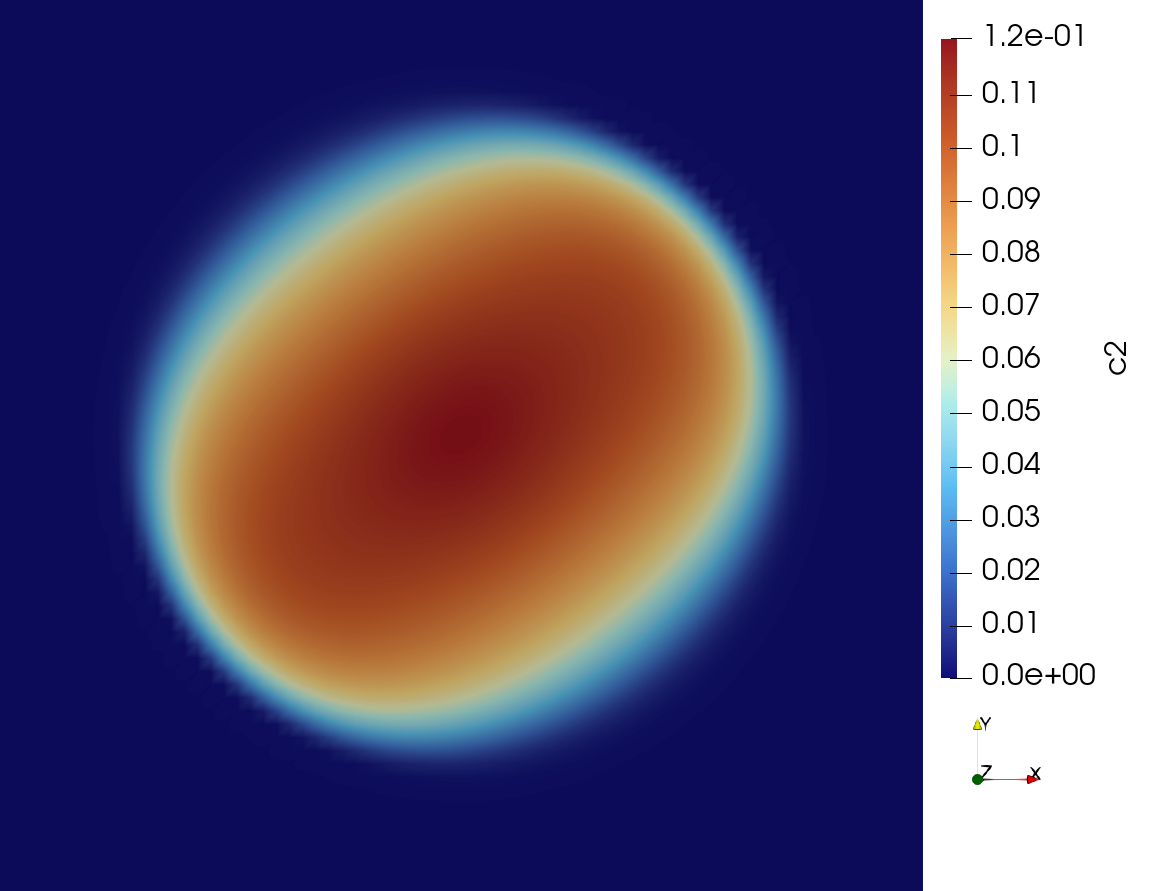}}
    \hfill
    \subfloat[$c_2$ difference]{%
        \includegraphics[width=0.32\textwidth]{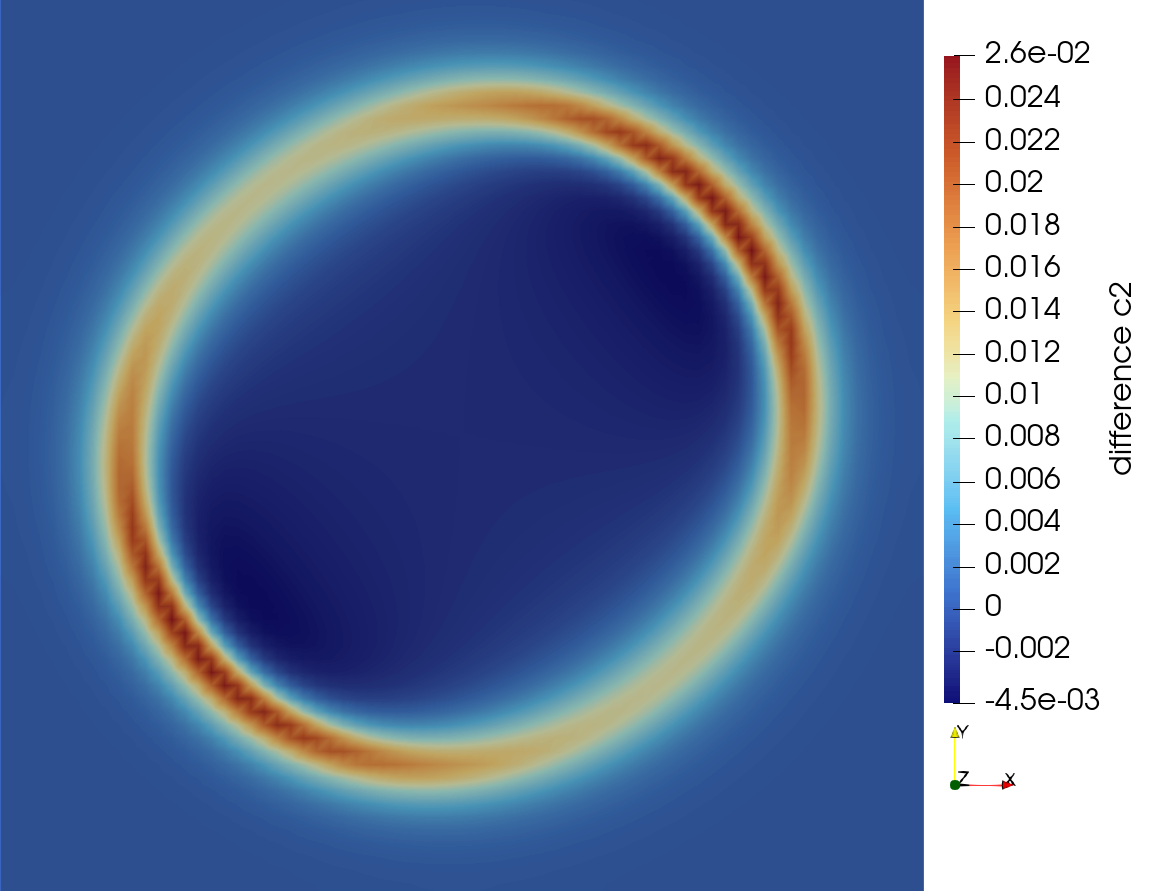}}

    \vspace{0.3cm}

    \subfloat[$\tau$ (non-degenerate)]{%
        \includegraphics[width=0.32\textwidth]{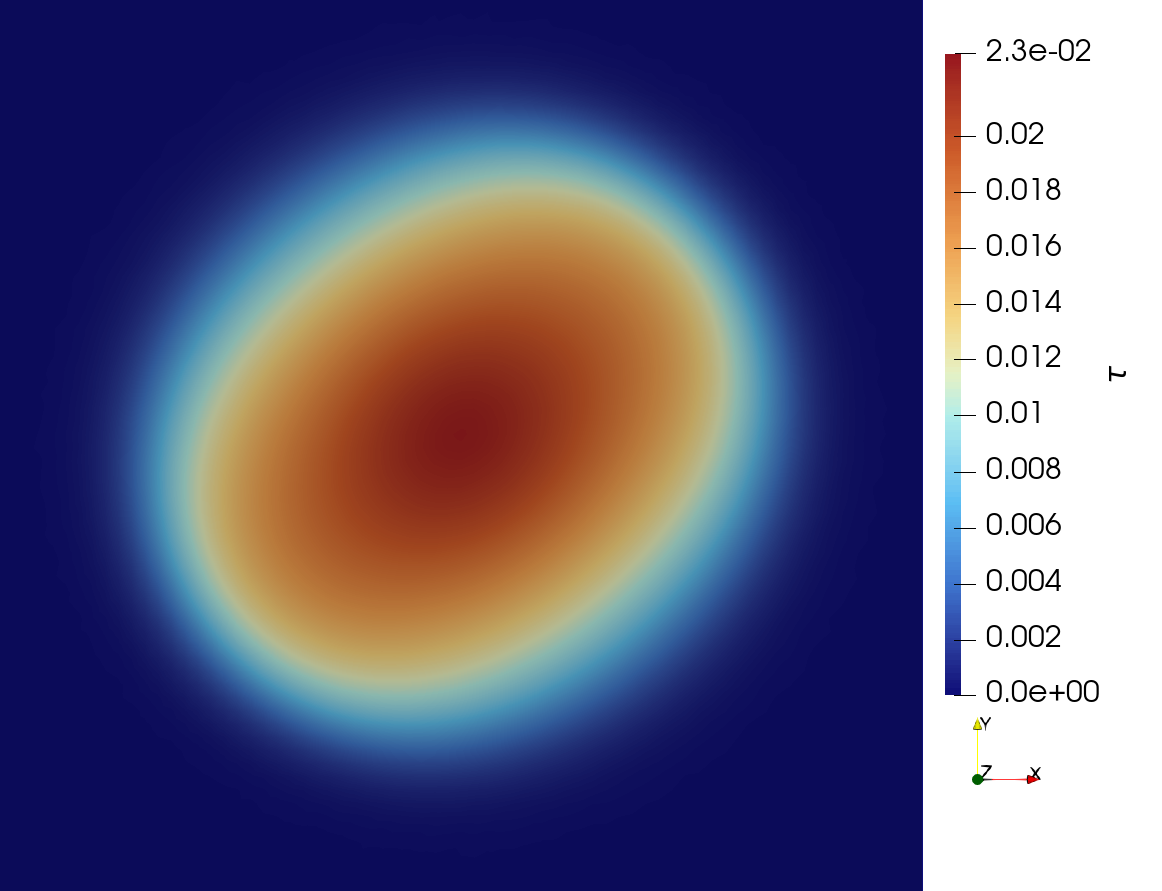}}
    \hfill
    \subfloat[$\tau$ (degenerate)]{%
        \includegraphics[width=0.32\textwidth]{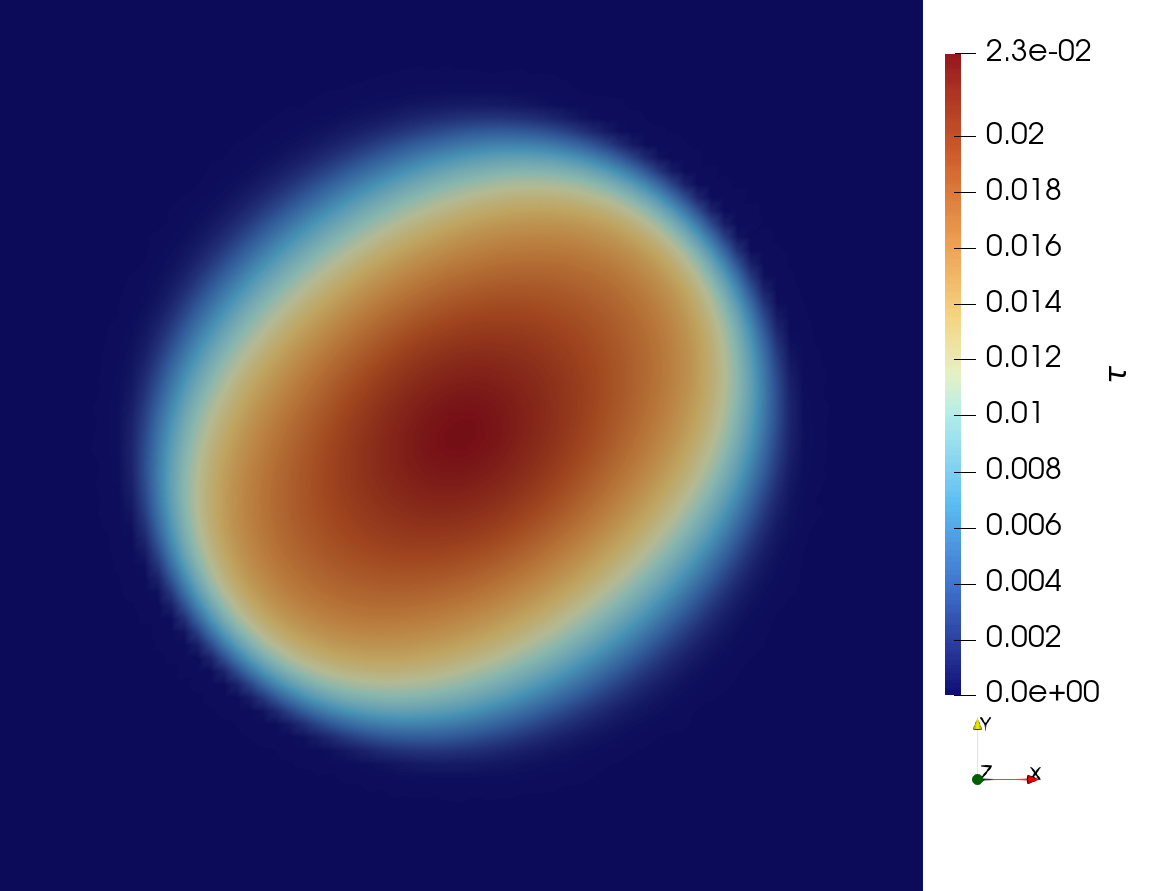}}
    \hfill
    \subfloat[$\tau$ difference]{%
        \includegraphics[width=0.32\textwidth]{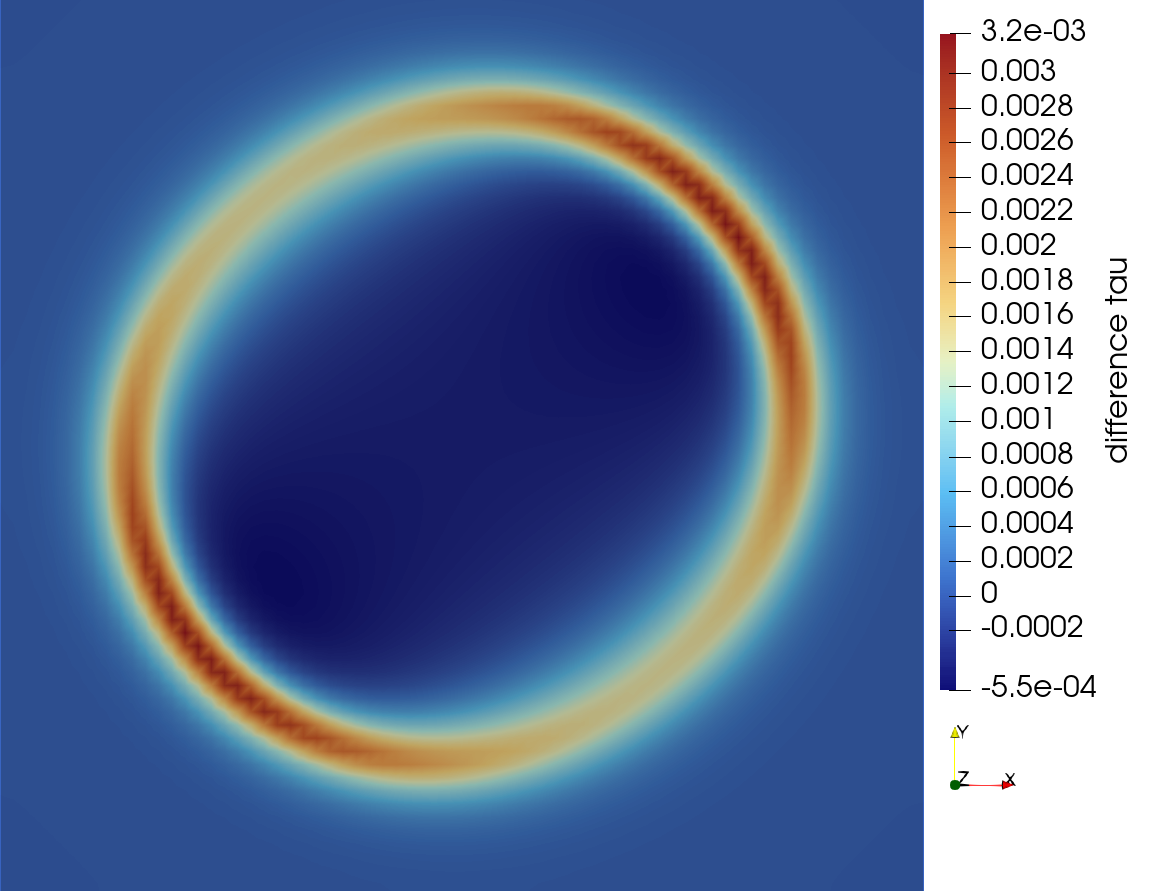}}

    \caption{Numerical solutions for the non-degenerate case (left), the radially degenerate case (middle), and their pointwise difference (right) after 420 time steps ($21\,\text{days}$).}
    \label{fig:comparison_case12}
\end{figure}

In Figure \ref{fig:comparison_case12}, we present a direct 3-column comparative evaluation after $21\,\text{days}$ ($t = 4.2$), contrasting the non-degenerate solutions (left column) with the radially degenerate results (middle column), alongside their pointwise spatial differences (right column). The non-degenerate model exhibits smooth, unconstrained diffusion profiles extending beyond the cutoff radius due to the positive background diffusion $D_{\min}$. In contrast, the degenerate case clearly demonstrates a sharp finite propagation front where cell motility drops abruptly to zero, leading to higher cell and ECM aggregates within the diffusion region and highlighting the qualitative impact of parabolic degeneracy on tissue regeneration patterns. 

\subsubsection{Strongly degenerate multi-hole diffusion (Case 3)}

\begin{figure}[htbp]
	\centering
	\subfloat[$c_1(x,y,0)$]{%
		\includegraphics[width=0.32\textwidth]{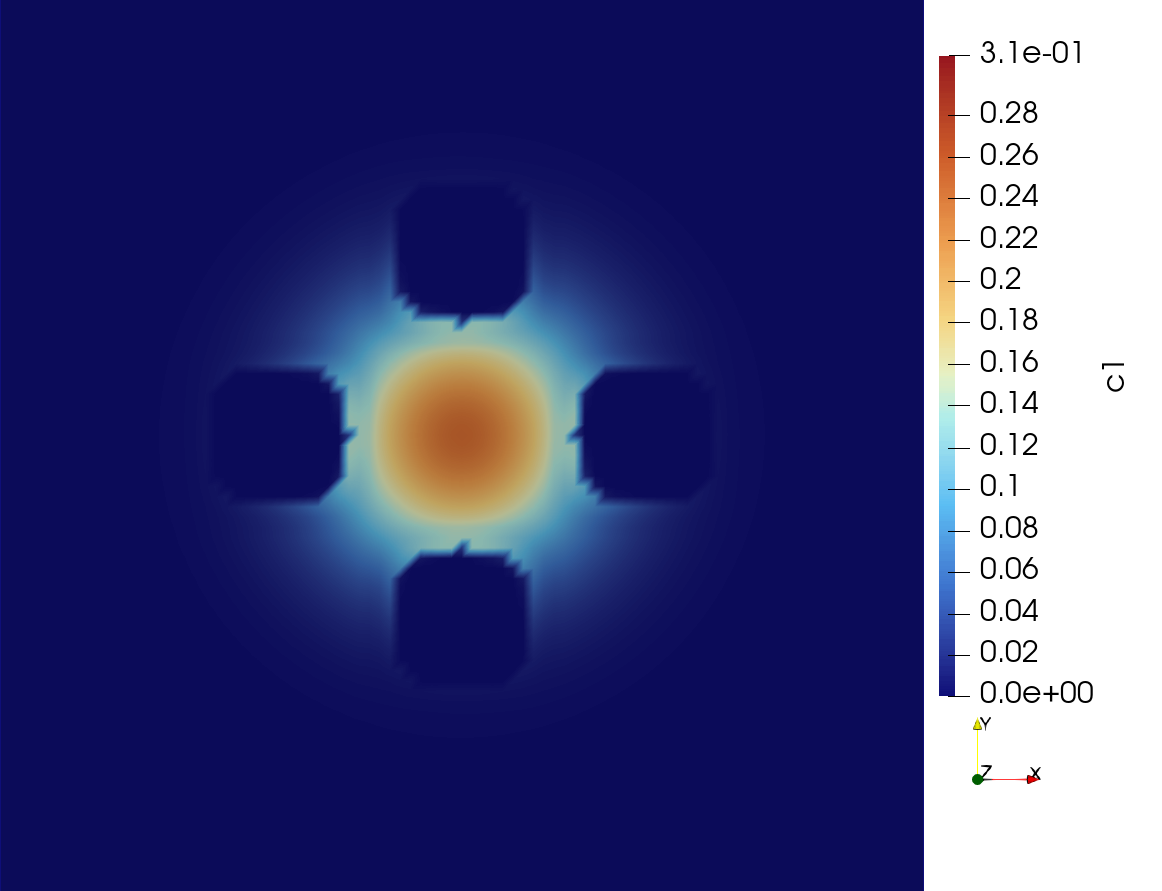}}
	\hfill
	\subfloat[$c_2(x,y,0)$]{%
		\includegraphics[width=0.32\textwidth]{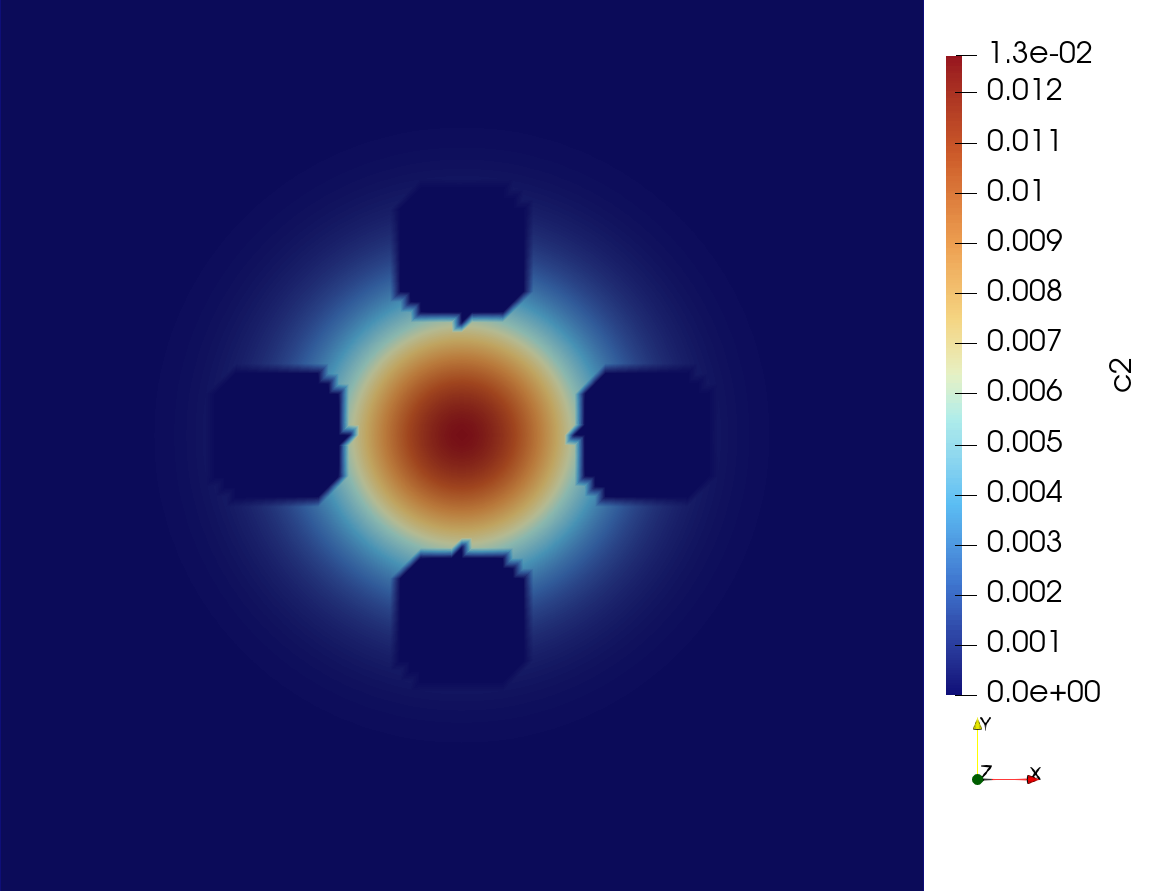}}
	\hfill
	\subfloat[$\tau(x,y,0)$]{%
		\includegraphics[width=0.32\textwidth]{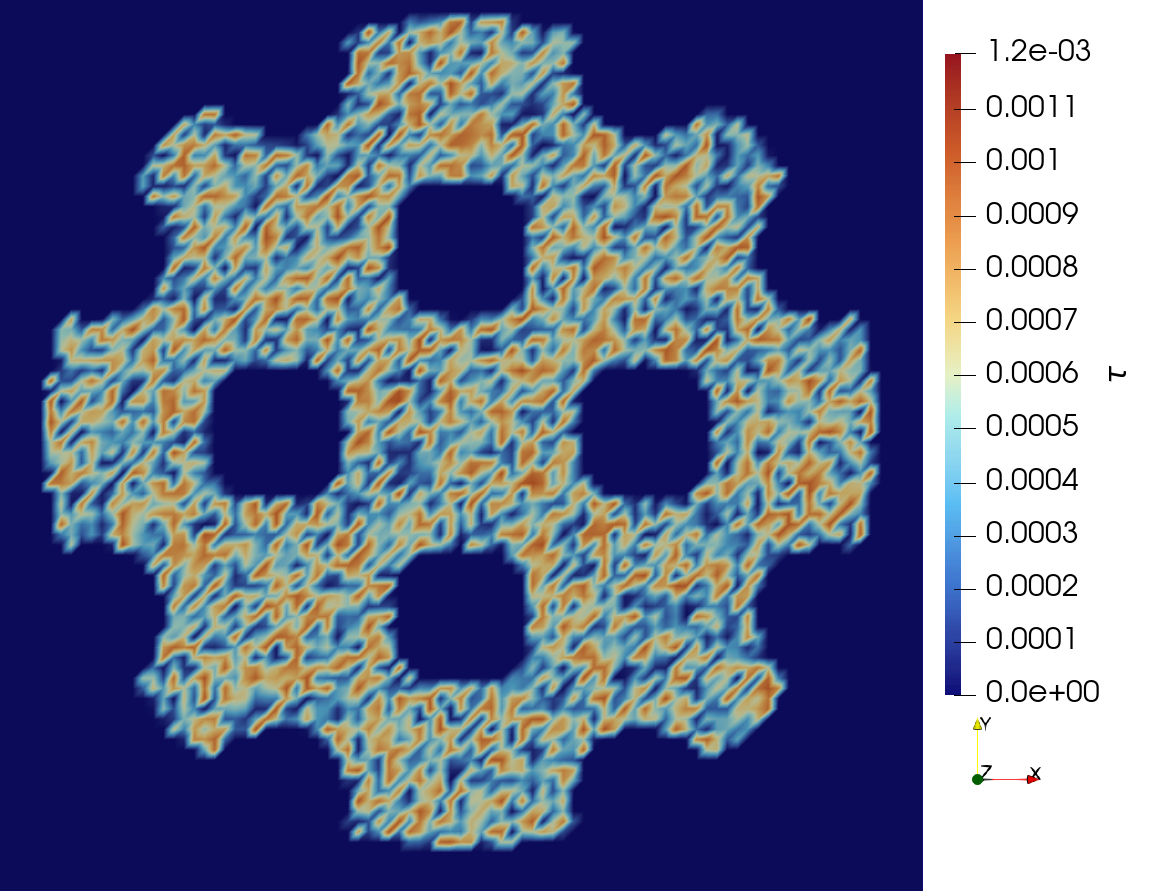}}
	\caption{Initial distributions of $c_1$, $c_2$, and $\tau$ at $t=0$ for Case 3.}
	\label{fig:init_case3}
\end{figure}

In this subsection we evaluate performance under complex topological constraints. Figure \ref{fig:init_case3} shows the modified initial conditions for Case 3, where MSCs avoid the newly introduced internal obstacle cores. 

\begin{figure}[htbp]
	\centering
	\subfloat[$c_1$]{%
		\includegraphics[width=0.32\textwidth]{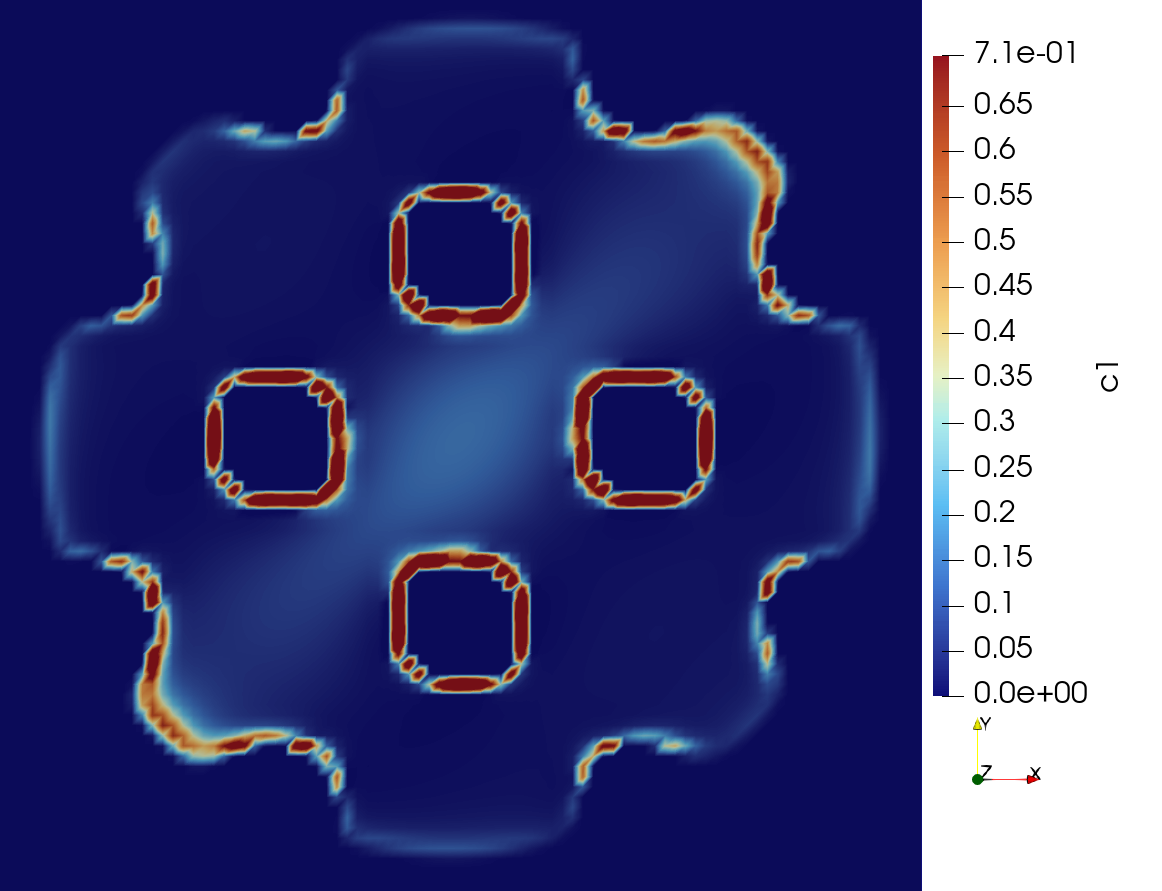}}
	\hfill
	\subfloat[$c_2$]{%
		\includegraphics[width=0.32\textwidth]{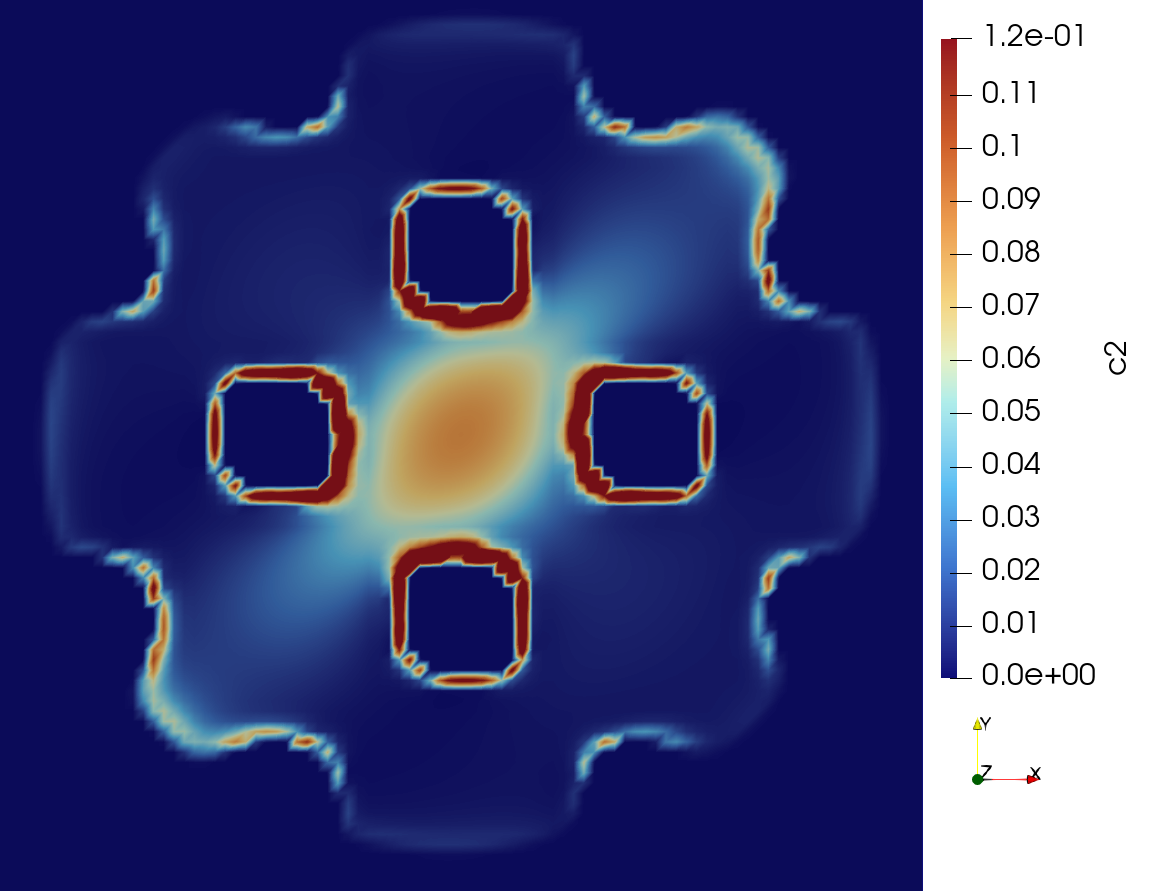}}
	\hfill
	\subfloat[$\tau$]{%
		\includegraphics[width=0.32\textwidth]{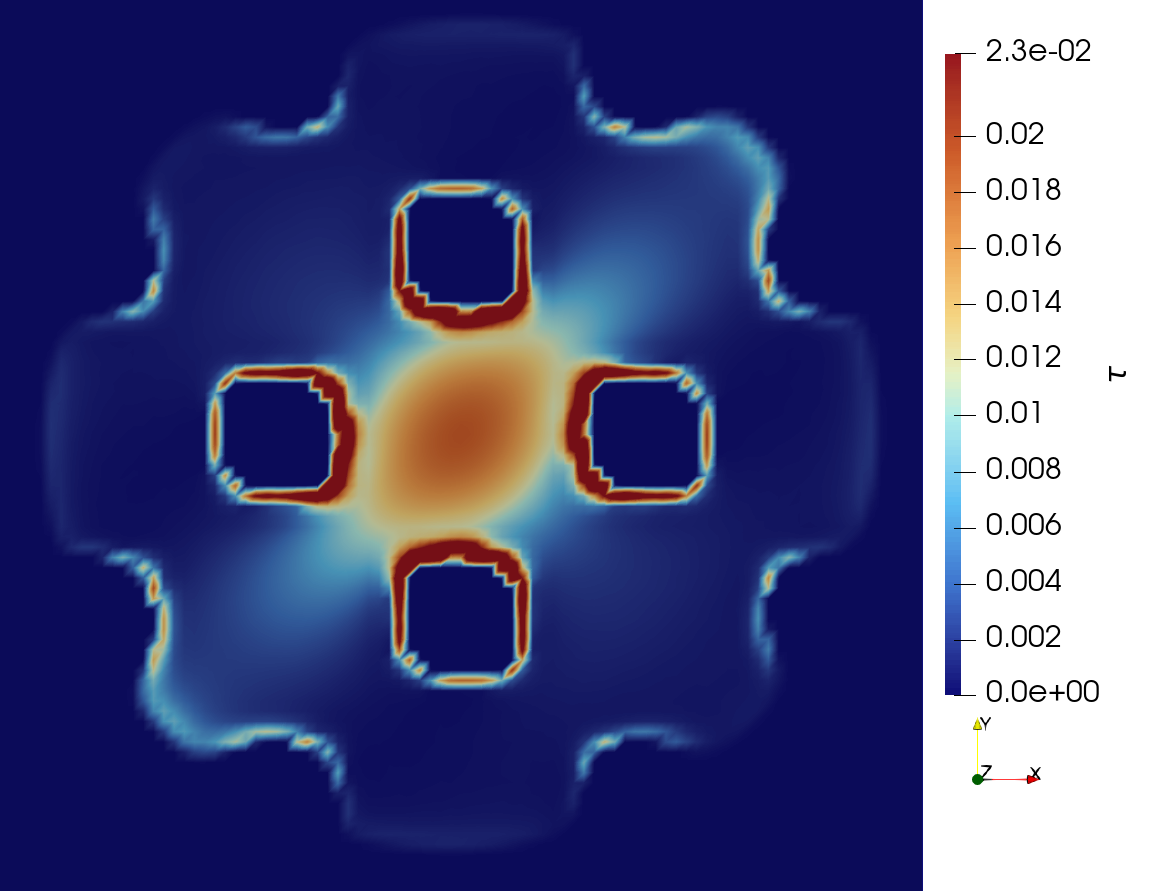}}
	\caption{Numerical solutions for $c_1$, $c_2$, and $\tau$ in the strongly degenerate multi-hole case (Case 3) after 420 time steps ($21\,\text{days}$).}
	\label{fig:results_case3}
\end{figure}

Figure \ref{fig:results_case3} illustrates the resulting distributions of $c_1$, $c_2$, and $\tau$ after $21\,\text{days}$ of simulation time under the multi-hole tensor structure introduced in \eqref{eq:g_multihole}. The numerical scheme effectively maintains zero concentration inside the impenetrable internal islands ($\mathbb{D} \equiv \mathbb O_2$). Proliferating cell fronts successfully navigate and channel around the obstacles without exhibiting unphysical artificial depletion in the core, confirming the robustness and stability of the mass-lumped finite element formulation combined with Raviart-Thomas flux projections. In this case, too, the suppressed motility due to degenrating diffusion tensor leads to accumulations of MSCs and therewith associated higher, more localized chondrocyte and ECM densities at the hole edges.

\subsection{Discussion}
{In this subsection, we summarize the numerical part of this work, in which we developed and applied a structure-preserving two-dimensional finite element framework in FreeFem++ to investigate the space-dependent, potentially degenerate PDE-ODE-ODE system. By combining continuous piecewise linear ($\mathbb{P}_1$) Lagrange finite elements for the primary variables ($c_1, c_2, \tau$) with a lowest-order Raviart-Thomas ($W_h$) projection for the auxiliary drift field $w = \nabla \cdot \mathbb{D} - \chi \mathbb{D}\nabla \tau$, we obtain a robust scheme that is applied to several test cases involving non-degenerate and degenerate diffusion tensors, including a strongly degenerate multi-hole geometry.}

\begin{figure}[htbp]
    \centering
    \subfloat[$c_1$ without mass lumping (consistent mass matrix)]{%
    \includegraphics[width=0.45\textwidth]{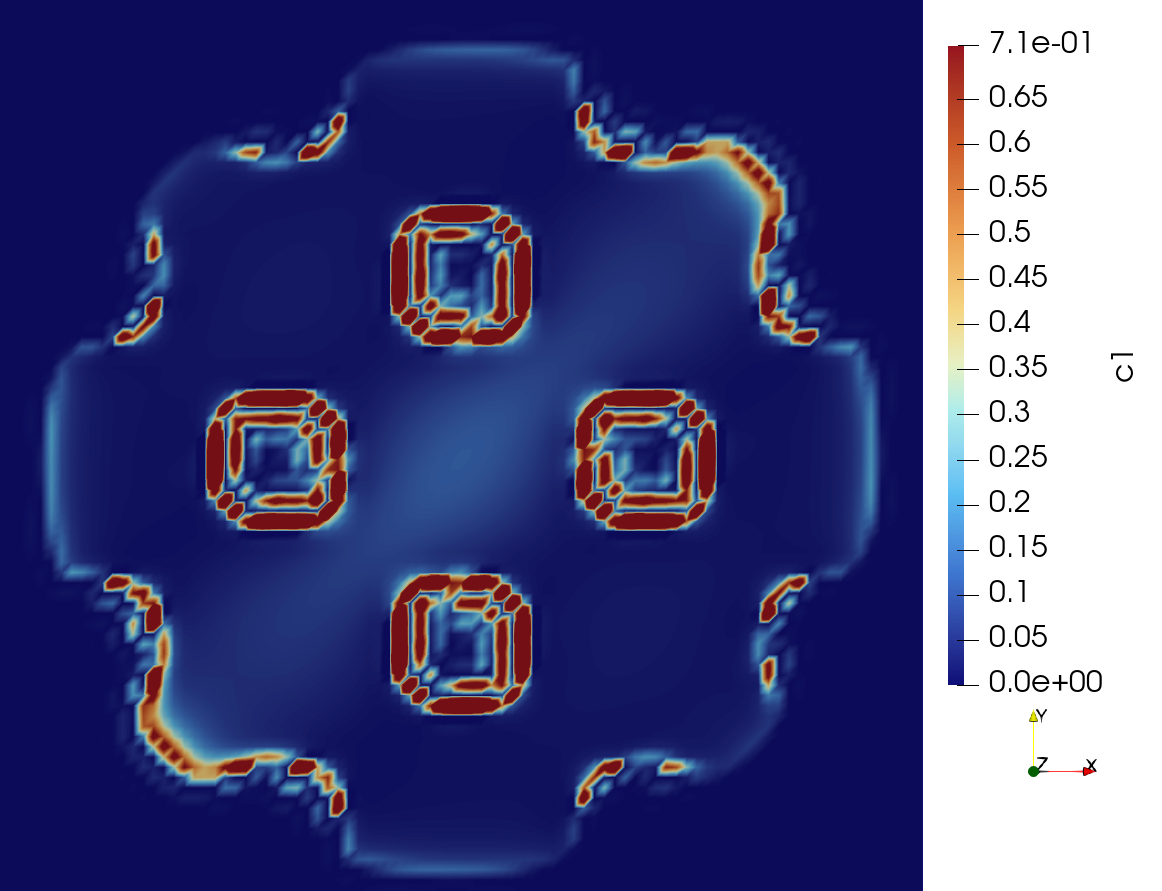}}
    \hfill
    \subfloat[$c_1$ with nodal mass lumping]{%
        \includegraphics[width=0.45\textwidth]{res/c1deg3.png}}
    \caption{Effect of mass lumping on the MSC density $c_1$ in the strongly degenerate multi-hole scenario (Case 3) at $t = 4.2$ ($21\,\text{days}$). (a) The consistent mass matrix causes spurious cross-interface leakage and subsequent proliferative growth inside the impenetrable obstacles ($\mathbb{D} \equiv \mathbb{O}_2$). (b) Nodal mass lumping strictly preserves zero concentration inside the obstacle cores.}
    \label{fig:lumping_comparison}
\end{figure}

From a numerical perspective, the simulations for the degenerate multi-hole geometry (Case 3) indicated that the consistent-mass formulation can produce nonzero numerical concentrations within the obstacle cores, where the diffusion tensor satisfies $\mathbb{D} \equiv \mathbb O_2$. This behavior is consistent with the off-diagonal entries of the consistent mass matrix, which couple neighboring nodal degrees of freedom even where the diffusion contribution vanishes. We therefore employ nodal mass lumping via vertex quadrature for the time-derivative and reaction terms, yielding a diagonal lumped mass matrix and nodal reaction contributions. As illustrated in Figure \ref{fig:lumping_comparison} for the mesenchymal stem cell density \(c_1\), the mass-lumped scheme maintains zero nodal values in the interior of the obstacle cores throughout the simulated time interval; analogous behavior is observed for \(c_2\) and \(\tau\). This comparison demonstrates the practical benefit of mass lumping for the present numerical setup, although positivity preservation of the full coupled scheme does not follow from mass lumping alone.

%The numerical simulations (also performed over a longer time span - not shown here) suggest that a degenerating diffusion tensor $\mathbb D$ does not necessarily lead to excessive accumulations of cells. Alone the presence of $\mathbb D$ in all motility terms, including haptotaxis, seems, however, not to preclude the explosion of solutions (see, \cite{Winkler2017,Winkler2018} for closely related systems in 1D). Instead, its combination with transitions between the two cell populations and a proliferation term possibly turning into degradation for large cell and ECM densities seems to keep the solution bounded. The occurence of $\mathbb D$ in the taxis term is due to the way we performed the upscaling from the KTE. In \cite{Zhigun2022} was performed a rigorous upscaling leading to a macroscopic PDE in which the cell diffusion tensor is only present in the myopic diffusion, but not in the taxis term; such deduction, even if merely informal, would lead to a model where the degeneracy of $\mathbb D$ no longer turns off the aggregation-supporting taxis.

Overall, summarizing the computational results, the experiments over the simulated $21$\ days ($t = 4.2$) highlight two key biological and physical insights into meniscus tissue regeneration:
\begin{itemize}
    \item \textbf{Consequences of parabolic degeneracy (Cases 1 and 2):} A direct comparison under identical initial conditions demonstrates how local mobility constraints alter spatial spreading. While a strictly positive background diffusion ($D_{\min} = 0.10$) permits unhindered dispersion across the domain, setting $D_{\min} = 0$ enforces a sharp motility barrier at the boundary of the support region. This confines mesenchymal stem cells (MSCs) to the active tissue core, preventing outward leakage and leading to more concentrated cellular and extracellular matrix accumulations.
    \item \textbf{Channeling and boundary aggregation in perforated media (Case 3):} In the multi-hole scenario, the heterogeneous tensor structure forces proliferating cell populations to navigate dynamically through interconnected migration channels. Furthermore, steep diffusivity gradients at the obstacle interfaces generate an inward myopic drift, trapping cells along the perimeter. This boundary-focused aggregation subsequently drives elevated, highly localized chondrogenesis and matrix synthesis around the structural pores.
\end{itemize} 

These results suggest that local variations and degeneracies in scaffold permeability and tissue anisotropy act as geometric guiding mechanisms that directly shape cell colonization and matrix synthesis.

\paragraph{Acknowledgments:}
This work was funded by the German Research Foundation (DFG) within SPP 2311 (Grant No. 465242756). We also acknowledge partial funding by the Research Initiative MSO (Modeling, Simulation and Optimisation) at the RPTU. NM gratefully acknowledges partial support by REaCH at the RPTU. AH gratefully acknowledges Dr Jeremias Arf (RPTU) for valuable discussions and helpful insights.
%%%%%%%%%%%%%%%%%%%%%%%%%%%%%%%%%%%%%%%%%%%%%%%%%%%%%%%%%%%%%%%%%%%%%%%%%%%%%%%%%%%%%%%%%%%%%%%%%%%%%%%%%%%%%%%%%%%%%%%%%%%%%%%%%%%%%%%%%%%%%%%%%%%%%%%%%%%%%%%%%%%%%%%%%%%%%%%%%%%%%%%%%%%%%%%%%%%%%%%%%%%%%%%%%%%%%%%%%%%%%%%%%
%%%%%%%%%%%%%%%%%%%%%%%%%%%%%%%%%%%%%%%%%%%%%%%%%%%%%%%%%%%%%%%%%%%%%%%%%%%%%%%%%%%%%%%%%%%%%%%%%%%%%%%%%%%%%%%%%%%%%%%%%%%%%%%%%%%%%%%%%%%%%%%%%%%%%%%%%%%%%%%%%%%%%%%%%%%%%%%%%%%%%%%%%%%%%%%%%%%%%%%%%%%%%%%%%%%%%%%%%%%%%%%%%
%%%%%%%%%%%%%%%%%%%%%%%%%%%%%%%%%%%%%%%%%%%%%%%%%%%%%%%%%%%%%%%%%%%%%%%%%%%%%%%%%%%%%%%%%%%%%%%%%%%%%%%%%%%%%%%%%%%%%%%%%%%%%%%%%%%%%%%%%%%%%%%%%%%%%%%%%%%%%%%%%%%%%%%%%%%%%%%%%%%%%%%%%%%%%%%%%%%%%%%%%%%%%%%%%%%%%%%%%%%%%%%%%

\bibliographystyle{abbrv}
\bibliography{Tensor}

\end{document}